\documentclass[11pt, reqno]{amsart}
\usepackage{amsxtra}
\usepackage{tikz-cd}

\let\oldtocsection=\tocsection

\let\oldtocsubsection=\tocsubsection

\let\oldtocsubsubsection=\tocsubsubsection

\renewcommand{\tocsection}[2]{\hspace{0em}\oldtocsection{#1}{#2}}
\renewcommand{\tocsubsection}[2]{\hspace{2em}\oldtocsubsection{#1}{#2}}
\renewcommand{\tocsubsubsection}[2]{\hspace{4.5em}\oldtocsubsubsection{#1}{#2}}

\makeatletter
\def\subsection{\@startsection{subsection}{2}
 \z@{.5\linespacing\@plus.7\linespacing}{-.5em}
 {\normalfont\bfseries}}
\def\subsubsection{\@startsection{subsubsection}{3}
 \z@{.5\linespacing\@plus.7\linespacing}{-.5em}
 {\normalfont\bfseries}}
\makeatother

\theoremstyle{plain}

\numberwithin{equation}{section}

\newtheorem{theorem}{Theorem}[section]
\newtheorem{lemma}[theorem]{Lemma}
\newtheorem{definition-theorem}[theorem]{Definition-Theorem}
\newtheorem{definition-construction}[theorem]{Definition-Construction}
\newtheorem{proposition}[theorem]{Proposition}
\newtheorem{corollary}[theorem]{Corollary}
\newtheorem{definition}[theorem]{Definition}
\newtheorem{example}[theorem]{Example}
\newtheorem{remark}[theorem]{Remark}
\newtheorem{notation}[theorem]{Notation}
\newtheorem{construction}[theorem]{Construction}
\newtheorem{assumption}[theorem]{Assumption}
\newtheorem{assumption-notation}[theorem]{Assumption-Notation}
\newtheorem{lemma-definition}[theorem]{Lemma-Definition}
\newtheorem{corollary-definition}[theorem]{Corollary-Definition}
\newtheorem{lemma-notation}[theorem]{Lemma-Notation}
\newtheorem{question}[theorem]{Question}
\newtheorem{remark-definition}[theorem]{Remarks-Definition}
\newtheorem{notation-remark}[theorem]{Notation-Remarks}
\newtheorem{conjecture}[theorem]{Conjecture}
\newcommand \bth[1] { \begin{theorem}\label{t#1} }
\newcommand \ble[1] { \begin{lemma}\label{l#1} }

\newcommand \bpr[1] { \begin{proposition}\label{p#1} }
\newcommand \bco[1] { \begin{corollary}\label{c#1} }
\newcommand \bconj[1] { \begin{conjecture}\label{co#1} }
\newcommand \bcons[1] { \begin{construction}\label{cons#1} }
\newcommand \bde[1] { \begin{definition}\label{d#1}\rm }
\newcommand \bex[1] { \begin{example}\label{e#1}\rm }
\newcommand \bre[1] { \begin{remark}\label{r#1}\rm }
\newcommand \bnota[1] {\begin{notation}\label{n#1}\rm }
\newcommand \bas[1] { \begin{assumption}\label{a#1}\rm }

\newcommand \bld[1] { \begin{lemma-definition}\label{ld#1} }

\newcommand \bqu[1] { \begin{question}\label{q#1}\rm }

\newcommand {\eth} { \end{theorem} }
\newcommand {\ele} { \end{lemma} }

\newcommand {\epr} { \end{proposition} }
\newcommand {\econj} { \end{conjecture} }
\newcommand {\eco} { \end{corollary} }
\newcommand {\ede} { \end{definition} }
\newcommand {\eex} { \end{example} }
\newcommand {\ere} { \end{remark} }
\newcommand {\enota} { \end{notation} }
\newcommand {\eas} {\end{assumption}}
\newcommand {\econs} {\end{construction}}

\newcommand {\eld}{ \end{lemma-definition} }

\newcommand {\equ} {\end{question}}
\newcommand \thref[1]{Theorem \ref{t#1}}
\newcommand \leref[1]{Lemma \ref{l#1}}
\newcommand \prref[1]{Proposition \ref{p#1}}
\newcommand \coref[1]{Corollary \ref{c#1}}
\newcommand \deref[1]{Definition \ref{d#1}}

\newcommand \exref[1]{Example \ref{e#1}}
\newcommand \reref[1]{Remark \ref{r#1}}

\newcommand \ldref[1]{Lemma-Definition \ref{ld#1}}

\def \RR {{\mathbb R}}         
\def \CC {{\mathbb C}}
\def \ZZ {{\mathbb Z}}

\def \TT {{\mathbb T}}
\def \KK {{\mathbb K}}
\def \QQ {{\mathbb Q}}
\def \PP {{\mathbb P}}

\def \calQ  {{\mathcal{Q}}}

\def \calP {{\mathcal{P}}}

\def \calS {{\mathcal{S}}}
\def \calT {{\mathcal{T}}}

\def \lan {\langle}
\def \ran {\rangle}
\def \ol {\overline}
\def \wt {\widetilde}

\def \g  {\mathfrak{g}}   

\def \t  {\mathfrak{t}}

\newcommand{\beqa}{\begin{eqnarray*}}                     
\newcommand{\eeqa}{\end{eqnarray*}}
\def \hs {\hspace{.2in}}

\def \bfu {{\bf u}}
\def \bfx {{\bf x}}
\def \bfy {{\bf y}}
\def \bfp {{\bf p}}

\def \ot {\otimes}

\def \FF {\mathbb{F}}

\def \PP {\mathbb{P}}

\def \calK {{\mathcal{K}}}
\def \calM {{\mathcal{M}}}

\def \calA {{\mathcal{A}}}
\def \calI {{\mathcal{I}}}
\def \calJ {{\mathcal{J}}}
\def \calK {{\mathcal{K}}}

\def \Acc {A^{c, c^{-1}}}

\def \Lcc {L^{c, c^{-1}}}

\def \oc {\overline{c}}

\def \wdu {\widetilde{{\bf u}}}

\def \bsig {\calS}
\def \obsig {\overline{\calS}}
\def \mathdash {-\!\!\!-}

\def \oP {\overline{P}}

\def \ot {\overline{t}}

\def \obfu {\overline{\bfu}}
\def \obfp {\overline{\bfp}}
\def \bb {\backslash\!\!\backslash}
\def \fX {{\mathfrak{X}}}
\def \bspr {\bsig^{\rm prin}}

\def \bfm {{\bf m}}
\def \ol {\overline{l}}

\def \BFZ {{\scriptscriptstyle {\rm BFZ}}}
\def \why {\widehat{y}}

\def \pbsig {\bsig^\prime}
\def \dbsig {\bsig^\dag}

\def \mF {{\mathfrak{F}}}
\def \ssq {{\scriptstyle \diamond}}
\def \calU {\mathcal{U}}
\def \calD {\mathcal{D}}

\def \ou {\overline{u}}
\def \ocalS {\overline{\calS}}
\def \ov {\overline{v}}

\def \cwc {{\bf c} {\bf w}_0({\bf c})}
\def \frakF {{\mathfrak{F}}}
\def \bfF {{\bf F}}
\def \BA {B_{\scriptstyle A}}
\def \UA {U_{\scriptstyle A}}

\begin{document}

\setlength{\baselineskip}{1.2\baselineskip}
\vspace{-.35in}
\title[Friezes of cluster algebras of geometric type]
{Friezes of cluster algebras of geometric type}
\author{Antoine de Saint Germain}
\address{
Department of Mathematics   \\
The University of Hong Kong \\
Pokfulam Road               \\
Hong Kong}
\email{adsg96@hku.hk}
\author{Min Huang}
\address{
School of Mathematics (Zhuhai)  \\
Sun Yat-sen University \\
Zhuhai              \\
China}
\email{huangm97@mail.sysu.edu.cn}
\author{Jiang-Hua Lu}
\address{
Department of Mathematics   \\
The University of Hong Kong \\
Pokfulam Road               \\
Hong Kong}
\email{jhlu@maths.hku.hk}
\date{}

\begin{abstract}
For a cluster algebra $\calA$ over $\QQ$ of geometric type, a {\it frieze} of $\calA$ is defined to be a 
$\QQ$-algebra homomorphism from $\calA$ to $\QQ$ that takes positive integer values on all cluster variables and all frozen variables. 
We present some basic facts on friezes, including frieze testing criteria, the notion of
{\it frieze points} when $\calA$ is finitely generated, and pullbacks of friezes under certain $\QQ$-algebra homomorphisms.
When the cluster algebra $\calA$ is acyclic, we define {\it frieze patterns associated to acyclic seeds of $\calA$}, generalizing the
{\it frieze patterns with coefficients of type $A$} studied by J. Propp and by M. Cuntz, T. Holm, and P. J{\o{}}rgensen, and
we give a sufficient condition for such frieze patterns to be equivalent to 
friezes. 
For the special cases when $\calA$ has an acyclic seed with 
either trivial coefficients, principal coefficients,
or what we call the {\it BFZ coefficients} (named after A. Berenstein, S. Fomin, and A. Zelevinsky), we identify
frieze points of $\calA$ 
both geometrically as certain positive integral points in explicitly described affine varieties and Lie theoretically (in the finite case) 
in terms of reduced double Bruhat cells and generalized minors on the associated semi-simple Lie groups. 
Furthermore, extending the gliding symmetry of the classical Coxeter frieze patterns of type $A$, we determine the symmetry of frieze patterns of any finite type with arbitrary coefficients.
\end{abstract}
\maketitle

\tableofcontents
\addtocontents{toc}{\protect\setcounter{tocdepth}{1}}
\section{Introduction and main results}
\subsection{Introduction}
The notion of {\it frieze patterns} was first introduced by Coxeter \cite{Cox} 
and further studied by Conway and Coxeter \cite{Con-Cox1, Con-Cox2}. Motivated by the recently discovered connections of Coxeter's frieze patterns with cluster algebras \cite{CC}, Assem, Reutenauer, and Smith introduced in \cite{ARS:friezes} frieze patterns associated to an arbitrary symmetrizable generalized Cartan matrix (called {\it frieze sequences} in \cite{Keller:Sche:linear}). 

\bde{:frieze-1}  \cite[$\S$3]{ARS:friezes} 
{\rm For an $r \times r$ symmetrizable generalized Cartan matrix $A = (a_{i,j})$, 
a {\it frieze pattern associated to $A$} is 
 a map $f: [1, r] \times \ZZ_{\geq 0} \rightarrow \ZZ_{>0}$
 satisfying
\begin{equation}\label{eq:ARS}
f(i, m)\, f(i, m+1) =  1+ \prod_{j=i+1}^r f(j, m)^{-a_{j,i}} \prod_{j = 1}^{i-1} f(j, m+1)^{-a_{j, i}}
\end{equation}
for every  $(i, m) \in [1, r] \times \ZZ_{\geq 0}$.
Here and for the rest of the paper, $\ZZ_{\geq 0}$ (resp. $\ZZ_{>0}$) denotes the set of all non-negative (resp. positive) integers, and $[1, r] = \{1, 2, \ldots, r\}$.
\hfill $\diamond$
}
\ede

The relations in \eqref{eq:ARS} are recursive with respect to the total order on $[1, r] \times \ZZ_{\geq 0}$ given by
\begin{equation}\label{eq:order-0}
(i, m) < (i', m') \hs \mbox{iff}\hs m <  m' \;\; \mbox{or} \; \;m = m' \;\; \mbox{and} \; \;i< i',
\end{equation}
and  a frieze pattern is uniquely determined by its initial values $\{f(i, 0): i \in [1, r]\}$. 
Frieze patterns associated to the standard Cartan matrix of type $A_r$ are precisely Coxeter's frieze patterns. 
 It was shown recently by Gunawan and Muller \cite{GM:finite} that 
the numbers of frieze patterns associated to Cartan matrices of finite type are finite, and their exact or conjectured numbers can be found in 
\cite{baur:friezes-higher-slk, FP:Dn}.
Other generalizations of Coxeter's frieze patterns include {\it frieze patterns with coefficients of type $A_r$} by Propp \cite{Propp:integers} and by
Cuntz, Holm, and J\o{}rgensen 
\cite{cuntz-et-al:coefficients} (see \exref{:Ptolemy}), 
and $SL_k$-friezes  for $k \geq 3$ (see \cite{baur:friezes-higher-slk, cuntz:wild-frieze}).
We refer to \cite{Sophie-M:survey} for a survey on the history and generalizations of Coxeter's frieze patterns and their relations to other fields such as linear difference equations
and integrable systems. 

Motivated by frieze patterns, the following definition has 
recently emerged (see \cite{F:non-zero, FP:Dn, GSch} and more particularly \cite{GM:finite} for cluster algebras with trivial coefficients).

\bde{:frieze-2}  
Let $\KK = \ZZ$ or a field of characteristic zero and $\calA_\KK$ be a $\KK$-cluster algebra of geometric type (see \deref{:A-0}). {\rm A {\it frieze}
of $\calA_\KK$  is 
a $\KK$-algebra homomorphism 
\[
h:\;\; \calA_\KK \longrightarrow  \KK
\]
 that  takes positive integral values on
all cluster variables and all frozen variables, collectively referred to as  extended cluster variables, 
of $\calA_\KK$. 
\hfill $\diamond$
}
\ede

To explain the connection between frieze patterns associated to Cartan matrices and friezes of cluster algebras, fix an $r \times r$ symmetrizable generalized Cartan matrix $A = (a_{i, j})$, and consider the 
skew-symmetrizable matrix 
\begin{equation}\label{eq:Bo-00}
\BA = \left(\begin{array}{cccc} 0 & -a_{1,2}  & \cdots & -a_{1,r}\\ 
a_{2,1} & 0 & \cdots & -a_{2,r}\\
\vdots & \vdots & \ddots& \vdots\\
\vspace{.08in}
a_{r,1} & a_{r,2} & \cdots &0\end{array}\right).
\end{equation}
Let $\calA_\ZZ(\BA)$  be the cluster algebra with trivial coefficients that has 
$\{u(i, 0): i \in [1, r]\}$
 as an initial cluster  and $\BA$ as the initial mutation matrix (see $\S$\ref{ss:A}). Applying the {\it knitting algorithm} 
(\cite[$\S$2.2]{Keller:06} and  \cite[$\S$2]{Keller:Sche:linear}),
one obtains cluster variables 
\[
\{u(i, m): (i, m) \in  [1, r] \times \ZZ_{\geq 0}\}
\]
in $\calA_\ZZ(\BA)$ satisfying mutation relations (see $\S$\ref{ss:belt} for detail)
\begin{equation}\label{eq:uim-intro}
u(i, m)\, u(i, m+1) =  1+ \prod_{j = i+1}^r u(j, m)^{-a_{j,i}} \prod_{j = 1}^{i-1} u(j, m+1)^{-a_{j, i}}, \hs (i, m) \in [1, r]\times \ZZ_{\geq 0}.
\end{equation}
Compairing \eqref{eq:ARS} and \eqref{eq:uim-intro}, it is clear that a frieze 
$h: \calA_\ZZ(\BA) \to \ZZ$
gives rise to a frieze pattern associated to $A$  by evaluating $h$ on the cluster variables $\{u(i, m): (i, m) \in  [1, r] \times \ZZ_{\geq 0}\}$. It is a priori not
clear, however, that the converse is true, i.e., that every frieze pattern associated to $A$ arises from a (necessarily unique) 
frieze $h: \calA_\ZZ(B_A) \to \ZZ$  via evaluation. 
The fact that this is indeed the case has been proved in \cite[Remark 6.13]{GM:finite} using a recent result of \cite{BNI20}.

\deref{:frieze-2}, valid for any cluster algebra $\calA_\KK$ of geometric type, naturally leads to geometrical as well as functorial considerations. When $\calA_\KK$ is acyclic, \deref{:frieze-2} also motivates the introduction of {\it frieze patterns} associated to acyclic seeds of $\calA_\KK$,  generally referred to as {\it frieze patterns with coefficients}. In this paper, we systematically study friezes and frieze patterns (with coefficients) in these directions. We now give more details on the main results of the paper. 

\subsection{Basics on friezes}\label{ss:main} 
We present in 
$\S$\ref{s:friezes} some basic concepts on friezes. 

The Positive Laurent Phenomenon for cluster algebras of geometric type, proved in  \cite[Proposition 11.2]{FZ:II} and \cite[Theorem 4.10]{GHKK}, 
ensures the existence of a special class of friezes, called
{\it unitary friezes}, for any cluster algebra. See \ldref{:unitary-frieze}.

When  $\calA_\CC$ is finitely generated,
we consider the affine variety $V = {\rm Spec}(\calA_\CC)$, 
and we call $v \in V$ a {\it frieze point} if all the extended cluster variables in $\calA_\CC$ 
take positive integer values at $v$.  Regarding each extended cluster of $\calA_\CC$ as a coordinate system on $V$, frieze points are then exactly those
$v \in V$ that have {\it positive integral  coordinates} in {\it every} such coordinate system. 
Existence of non-unitary friezes (see \cite{FP:Dn} for examples) is thus a rather remarkable fact from this geometric point of view. 

As a cluster algebra  typically has infinitely many cluster variables, it is natural to ask for frieze testing criteria
using a small number of extended cluster variables, which is 
desirable even when the cluster algebra is of finite type. We give in 
 \leref{:poly-0} examples of such testing criteria which cover a
rich class of polynomial cluster algebras coming from Lie theory such as coordinate rings on Schubert cells 
\cite{GLS:Kac-Moody} and double Bott-Samelson cells \cite{ShenWeng:BS}. See \exref{:poly}.
 
For a cluster algebra $\calA_\ZZ$ of geometric type, let $\calA_\ZZ^+$ be the  unital sub-semiring of $\calA_\ZZ$ generated by all the extended cluster variables of 
$\calA_\ZZ$. 
In $\S$\ref{s:pullbacks}, we introduce the category ${\bf Frieze}$, whose objects are cluster algebras $\calA_\ZZ$ of geometric type with
frozen variables non-invertible, 
and whose morphisms are ring homomorphisms $\psi: \calA_\ZZ \to \calA_\ZZ^\prime$ with {\it Property $\mF$}, in the sense
that 
\[
\psi(\calA_\ZZ^+) \subset 
(\calA_\ZZ^\prime)^+
\]
(see \deref{:mF}). Regarding $\ZZ$ as an object in ${\bf Frieze}$, i.e., the cluster algebra of rank $0$ and with trivial coefficients, 
a frieze of $\calA_\ZZ $ is nothing but a morphism from $\calA_\ZZ$ to $\ZZ$ in  {\bf Frieze}, and friezes have well-defined pullbacks under morphisms in 
${\bf Frieze}$. Our main examples of morphisms in ${\bf Frieze}$ include

(1) (\leref{:freezing-mF}) the inclusion map $\calA_\ZZ \to \calA_\ZZ^\prime$, where $\calA_\ZZ$ is obtained from $\calA_\ZZ^\prime$ by 
{\it freezing} a subset of cluster variables in a cluster of 
$\calA_\ZZ^\prime$;

(2) (\prref{:deleting}) ring homomorphisms $\calA_\ZZ \to \calA_\ZZ^\prime$, where $\calA_\ZZ^\prime$ is of finite type and is obtained from $\calA_\ZZ$ by  {\it deleting}
a subset of cluster variables in a cluster of 
$\calA_\ZZ$;

(3)  (\coref{:quasi-frakF} and \exref{:P0P0}) certain cluster algebra quasi-homomorphisms $\calA_\ZZ \to \calA_\ZZ^\prime$ in the sense of Fraser \cite{Fraser:quasi};

(4)  (\thref{:finite-mF}) all universal coefficient specializations for cluster algebras of finite type.

\subsection{Friezes and frieze patterns of acyclic cluster algebras}\label{ss:patterns-intro}
In $\S$\ref{s:patterns} we turn our attention to friezes of {\it acyclic} cluster algebras
(of geometric type), i.e. cluster algebras that have at least one acyclic 
seed \cite[Definition 1.14]{BFZ:III}. Recall (see \cite[p.36]{BFZ:III}) that
the extended mutation matrix of an acyclic seed is, up to simultaneous re-labelling of rows and columns, of the form 
\[
\widetilde{B} = \left(\begin{array}{c} \BA \\ P\end{array}\right)\, \stackrel{\rm denote}{=}\,(\BA\bb P),
\]
where $\BA$ is the skew-symmetrizable matrix given in \eqref{eq:Bo-00} for some generalized symmetrizable Cartan matrix $A$ and $P \in M_{l, r}(\ZZ)$ for any integer $l \geq 0$. 
Without loss of generality, we thus fix throughout $\S$\ref{s:patterns} an extended mutation matrix $(\BA\bb P)$ and consider the seed pattern $\calS(\BA \bb P) = \{\Sigma_t = ({\bf u}_t, \bfp, \widetilde{B}_t): t \in \TT_r\}$ 
with $\widetilde{B}_{t_0} = (\BA \bb P)$, where $(\TT_r, t_0)$ a rooted $r$-regular tree. Let $\calA_\ZZ(\BA\bb P, \emptyset)$ be the corresponding cluster algebra over $\ZZ[\bfp]$ (see $\S$\ref{ss:A} for more detail).
Consider the 
2-regular sub-tree $\TT_r^\flat$ of $\TT_r$ whose vertices are labelled 
\[ 
\cdots  \stackrel{r-1}{\mathdash} t(r,\!-2) \stackrel{r}{\mathdash} t(1, \!-1) \stackrel{1}{\mathdash} \cdots
\stackrel{r-1}{\mathdash} t(r,\!-1) \stackrel{r}{\mathdash} t(1, \!0) \stackrel{1}{\mathdash} 
\cdots \stackrel{r-1}{\mathdash} t(r,\! 0) \stackrel{r}{\mathdash} t(1,\! 1) 
\stackrel{1}{\mathdash} \cdots  
\]
where $t(1, 0) = t_0$ (see e.g. \cite[p.924]{RSW}). We call $\{\Sigma_{t(i, m)}: (i,m) \in [1, r] \times \ZZ\}$ 
the {\it acyclic belt of $\calS(\BA\bb P)$ through the seed $\Sigma_{t_0}$} (see $\S$\ref{ss:belt} for more detail). 
For $(i, m) \in [1, r] \times \ZZ$, write
\[
\widetilde{B}_{t(i, m)}= \left(\begin{array}{c} B_{t(i, m)} \\ P_{t(i, m)}\end{array} \right),
\]
let 
$p(i, m) \in \ZZ^l$ be the $i$th column of $P_{t(i, m)}$, and let $u(i, m)$ be the $i$th variable of the cluster ${\bf u}_{t(i, m)}$.
We prove in  \leref{:umi-mutation}, that for any $(i, m) \in [1, r]\times \ZZ$, the mutation 
of the seed $\Sigma_{t(i, m)}$ in direction $i$ gives the mutation relation
\begin{equation}\label{eq:ex-00}
	u(i, m)\, u(i, m+1)=  {\bf p}^{[p(i, m)]_+} + {\bf p}^{[-p(i, m)]_+}\prod_{j = i+1}^r u(j, m)^{-a_{j i}} \prod_{j =1}^{i-1} u(j, m+1)^{-a_{j i}}
\end{equation}
 in $\calA_\ZZ(\BA\bb P, \emptyset)$, which motivates the following definition (see $\S$\ref{ss:patterns} for more detail).

\bde{:pattern-intro}
{\rm 
A {\it frieze pattern associated 
to $A$ with coefficient matrix $P$} is a map 
\[
f: \;\; ([1, r] \times \ZZ) \cup \{p_1 \ldots, p_{l}\} \longrightarrow \ZZ_{>0}
\]
satisfying, for all $(i, m) \in [1, r] \times \ZZ$,
\[
f(i,m)\, f(i, m+1) =  f({\bf p})^{[p(i,m)]_+} + f({\bf p})^{[-p(i,m)]_+}\prod_{j=i+1}^r  f(j,m)^{-a_{j i}} \prod_{j = 1}^{i-1} f(j,m+1)^{-a_{j i}}.
\]
\hfill $\diamond$
}
\ede

Let $\fX(\BA \bb P)$ be the set of all cluster variables of $\calA_\ZZ(\BA\bb P, \emptyset)$. We refer to the assignment
\[
[1, r] \times \ZZ \longrightarrow \fX(\BA \bb P), \;\; (i, m) \longmapsto u(i, m)
\]
as the {\it generic frieze pattern} associated to $A$ with coefficient matrix $P$ or associated to $(\BA \bb P)$.
Every frieze  $h: \calA_\ZZ(\BA\bb P, \emptyset)\to \ZZ$ thus gives rise to a frieze pattern $f$ via 
\[
f(i, m) = h(u(i, m)), \hs (i, m) \in [1, r] \times \ZZ.
\]
 If $\calA_\ZZ(\BA\bb P, \emptyset)$ is generated over $\ZZ[\bfp]$ by a subset of 
$\{u(i, m): \,(i, m) \in [1, r] \times \ZZ\}$,
for example if $A$ is of finite type, 
we show in 
 \prref{:ARS-same} that the converse holds, i.e. that every frieze pattern associated to $(\BA \bb P)$ is obtained from a frieze of $\calA_\ZZ(\BA\bb P, \emptyset)$ via evaluation.

The set of vectors $\{p(i,m) : i \in [1,r], m \in \ZZ\}$ is completely determined
by $P$. We prove in \prref{:additive} that for every  $P \in M_{l, r}(\ZZ)$, the map
\[
    [1, r] \times \ZZ \longrightarrow \ZZ^l, \;\; (i, m) \longmapsto p(i, m),
\]
is a $\ZZ^l$-valued {\it cluster-additive function associated to $A$} in the sense of Ringel \cite{Ringel}, and that all $\ZZ^l$-valued cluster-additive functions associated to $A$ arise this way. Cluster-additive functions are further studied in 
\cite{CGL:additive} in the spirit of Fock-Goncharov duality.

In \exref{:Ptolemy}, we show that the {\it frieze patterns with coefficients} studied by Cuntz, Holm, and J{\o{}}rgensen 
\cite{cuntz-et-al:coefficients} are special cases of \deref{:pattern-intro}, where the cluster algebra is of type $A_r$ and the extended mutation matrix $(\BA\bb P)$ comes from a shell triangulation of an $(r+3)$-gon.

\subsection{Frieze points in examples}
In $\S$\ref{s:3-cases}, we study friezes of acyclic cluster algebras with three cases of coefficients. More precisely, for each symmetrizable generalized Cartan matrix $A$, we consider three extended mutation matrices, respectively given by
\[
\widetilde{B}^{\BFZ} = \left(\begin{array}{c} \BA \\ \UA\end{array}\right), \hs \widetilde{B}^{\rm prin} = \left(\begin{array}{c} \BA \\ I_r\end{array}\right) 
 \hs \mbox{and} \hs \widetilde{B}^{\rm triv} = \BA,
\]
where $I_r$ is the $r \times r$ identity matrix, and
\begin{equation}\label{eq:U0-00}
\UA= \left(\begin{array}{cccc} 1 & a_{1,2}  & \cdots & a_{1,r}\\ 
0 & 1 & \cdots & a_{2,r}\\
\vdots & \vdots & \ddots& \vdots\\
\vspace{.08in}
0 & 0 & \cdots &1\end{array}\right).
\end{equation}
Let  $\bsig^{\BFZ}$, $\bsig^{\rm prin}$, and $\bsig^{\rm triv}$ be the corresponding seed patterns,  and denote 
 by  
\[
 \calA_\CC(\bsig^{\BFZ}, \emptyset), \hs \calA_\CC(\bsig^{\rm prin}, \emptyset), \hs
\calA_\CC(\bsig^{\rm triv})
\]
 the respective cluster algebras in which the frozen variables are not invertible. 
We will refer to $\bsig^{\BFZ}$, $\bsig^{\rm prin}$, and $\bsig^{\rm triv}$ respectively  as the seed patterns associated to $A$ with {\it 
BFZ coefficients}, {\it principal coefficients}, and {\it trivial coefficients}. The term {\it BFZ coefficients} is due to the fact, to be explained in $\S$\ref{ss:BFZ}, that 
$\bsig^{\BFZ}$ gives the cluster structure on a reduced double Bruhat cell defined by A. Berenstein, S. Fomin, and A. Zelevinsky in \cite{BFZ:III}. 

The following \thref{:3-cases} summaries our geometric descriptions of friezes in the three cases of coefficients (\prref{:BFZ-test}, \prref{:prin-1}, $\S$\ref{ss:trivial}). 

\bth{:3-cases} 1)  One has ${\rm Spec} (\calA_\CC(\bsig^{\BFZ}, \emptyset)) \cong \CC^{2r}$ with frozen variables in $\calA_\CC(\bsig^{\BFZ}, \emptyset)\cong
\CC[z_1, \ldots, z_{2r}]$ given by $(p_1, \ldots, p_r)$, where
\begin{equation}\label{eq:pi-00}
p_i (z_1, \ldots, z_{2r}) = z_i z_{r+i} -\prod_{j=i+1}^r z_j^{-a_{j,i}} \prod_{j=1}^{i-1} z_{r+j}^{-a_{j,i}}, \hs i \in [1, r].
\end{equation}
The set of frieze points in $\CC^{2r}$ is 
$\{(z_1, \ldots, z_{2r}) \in (\ZZ_{>0})^{2r}: p_i(z_1, \ldots, z_{2r}) > 0, \; \forall \, i \in [1, r]\}$.

2) One has ${\rm Spec} (\calA_\CC(\bsig^{\rm prin}, \emptyset)) \cong V^{\rm prin}$, where $V^{\rm prin}$ 
 is the affine sub-variety of $\CC^{3r}$ defined in the  coordinates $(x_1, x_1^\prime, \ldots, x_r, x_r^\prime, y_1, \ldots, y_r)$ of $\CC^{3r}$ by the equations
 \[
x_ix_i^\prime = y_i \prod_{j=1}^{i-1} x_j^{-a_{j, i}} + \prod_{j=i+1}^r x_j^{-a_{j, i}},   \hs i \in [1, r].
\]
The set of frieze points in $V^{\rm prin}$ is $V^{\rm prin} \cap (\ZZ_{>0})^{3r}$;

3)  One has ${\rm Spec} (\calA_\CC(\bsig^{\rm triv})) \cong V^{\rm triv} \cong W^{\rm triv}$, where 
$V^{\rm triv} \subset \CC^{2r}$ 
 is defined in the coordinates $(z_1, \ldots, z_{2r})$ of $\CC^{2r}$ by the equations
\begin{equation}\label{eq:zzi}
z_i z_{r+i}  = 1 + \prod_{j=i+1}^r z_j^{-a_{j,i}} \prod_{j=1}^{i-1} z_{r+j}^{-a_{j,i}},  \hs i \in [1, r],
\end{equation}
and $W^{\rm triv}\subset \CC^{2r}$ is 
defined in the coordinates $(x_1, x_1^\prime, \ldots, x_r, x_r^\prime)$ of $\CC^{2r}$ by the equations
\begin{equation}\label{eq:xxi}
x_ix_i^\prime = \prod_{j=1}^{i-1} x_j^{-a_{j, i}} + \prod_{j=i+1}^r x_j^{-a_{j, i}},   \hs i \in [1, r].
\end{equation}
The sets of $\bsig^{\rm triv}$-frieze points  in  $V^{\rm triv}$ and $W^{\rm triv}$ are respectively given by 
$V^{\rm triv} \cap (\ZZ_{>0})^{2r}$ and  $ W^{\rm triv} \cap (\ZZ_{>0})^{2r}$.
\eth

When $A$ is of finite type, the recent result of 
Gunawan and Muller  \cite{GM:finite} on the finiteness of the number of frieze patterns (with trivial coefficients) implies that  the $r$ equations in \eqref{eq:zzi} have finitely many positive integral solutions, and the same holds for the $r$ equations in \eqref{eq:xxi}.

\subsection{Lie theoretical models}\label{ss:Lie-intro}
Assume the Cartan matrix $A$ is of finite type. In $\S$\ref{s:Lcc}, we let $G$ be a connected and simply connected complex Lie group of the same Cartan-Killing type as $A$, and we choose a simple root system
$\alpha_1, \ldots, \alpha_r$ 
for $G$ that gives rise to the Cartan matrix $A$. Consider the Coxeter element 
$c = s_{\alpha_1}s_{\alpha_2} \cdots s_{\alpha_r}$ in the Weyl group $W$ of $G$
 and the {\it reduced double Bruhat cell} 
 (see $\S$\ref{ss:bkgrnd-minor} for detail)
\[
\Lcc = N \oc B \cap B_- c^{-1} B_- \subset G.
\]
For $i \in [1, r]$, let $\omega_i$ be the fundamental 
weight  corresponding to $\alpha_i$. For $u, v \in W$, one then has the {\it generalized minor} $\Delta_{u\omega_i, v\omega_i} \in \CC[G]$
(see $\S$\ref{ss:bkgrnd-minor}). With $p_i \in \CC[z_1, \ldots, z_{2r}]$ given 
in \eqref{eq:pi-00} for $i \in [1, r]$, we show in  \thref{:Lcc-Af} that the map
\[
\rho: \Lcc \to \CC^{2r}, \;\;
g \longmapsto ({\Delta}_{\omega_1, \,\omega_1}(g), \, \ldots, \, {\Delta}_{\omega_r, \,\omega_r}(g), \, 
{\Delta}_{c\omega_1, \,c\omega_1}(g), \, \ldots, \, 
{\Delta}_{c\omega_r, c\omega_r}(g))
\]
gives an isomorphism of varieties from $\Lcc$ to the Zariski open subset 
\[
\CC^{2r}_{{\bf p} \neq 0} = \{(z_1, z_2, \ldots, z_{2r}) \in \CC^{2r}: \; \prod_{i=1}^r p_i(z_1, \ldots, z_{2r}) \neq 0\}
\]
of $\CC^{2r}$, and that the pullback by $\rho$ of the cluster structure $\calS^{\BFZ}$ on
$\CC^{2r}_{{\bf p} \neq 0}$ from \thref{:3-cases} is precisely the cluster structure on 
$\Lcc$ defined by A. Berenstein, S. Fomin, and A. Zelevinsky 
 \cite[Example 2.24]{BFZ:III}.

The reduced double Bruhat cell $\Lcc$ is also a geometric model for the
cluster algebra $\calA_\CC(\calS^{\rm prin})$ with {\it principal coefficients}, as shown by S.-W. Yang and A. Zelevinsky in \cite{YZ:Lcc}. In particular, all cluster variables of $\calA_\CC(\calS^{\rm prin})$ are shown to be
restrictions to $\Lcc$ of certain generalized minors.  
Set
\[
\Acc \stackrel{\rm def}{=} N \oc N \cap N_- \oc^{\, -1} N_- ,
\]
which now inherits the cluster structure $\calS^{\rm triv}$ with trivial coefficients and is isomorphic to  both cluster varieties $V^{\rm triv}$ and $W^{\rm triv}$ in  \thref{:3-cases}.
In $\S$\ref{s:Lcc}, we rephrase 
the frieze criteria in \thref{:3-cases}, and our results 
(\coref{:BFZ-in-minor} and \coref{:prin-in-minor}) are summarized as follows.

\bth{:frieze-in-minor}
1) A point $g \in \Lcc$ is an $\calS^{\BFZ}$-frieze point if and only if 
\[
\Delta_{\omega_i, \omega_i} \in \ZZ_{>0},\; \;\;\Delta_{c\omega_i, c\omega_i} \in \ZZ_{>0} \;\;\; \mbox{and} \;\; \;
\Delta_{\omega_i, c\omega_i} \in \ZZ_{>0}, \hs i \in [1, r];
\]

2) A point $g \in \Lcc$ is an $\calS^{\rm prin}$-frieze point if and only if 
\[
\Delta_{\omega_i, \,\omega_i} \in \ZZ_{>0},\; \;\;\Delta_{c\omega_i, \, s_i\omega_i} \in \ZZ_{>0} \;\;\; \mbox{and} \;\; \;
\Delta_{\nu_i, c\nu_i} \in \ZZ_{>0}, \hs i \in [1, r],
\]
where $\nu_i$ is a unique weight associated to $\alpha_i$ given in \leref{:y-ic}.

3) A point $g \in \Acc$ is an $\calS^{\rm triv}$-frieze point if and only if 
\[
	\Delta_{\omega_i, \,\omega_i}(g) \in \ZZ_{>0} \;\;\; \mbox{and} \;\;\; \Delta_{c\omega_i,c\omega_i}(g) \in \ZZ_{>0}, \hs  i \in [1,r],
\]
or, equivalently, if and only if 
\[
\Delta_{\omega_i,\, \omega_i}(g) \in \ZZ_{>0} \;\;\; \mbox{and} \;\;\; \Delta_{c\omega_i,s_i\omega_i}(g) \in \ZZ_{>0}, \hs i \in [1,r].
\]
\eth

In \cite{Lusztig1994}, G. Lusztig defined the {\it totally non-negative} part $G_{\geq 0}$ of $G$, and  it is shown in 
\cite[Theorem 3.1]{FZ:Oscillatory} that $g \in G$ lies in $G_{\geq 0}$ if and only if $\Delta(g) \geq 0$ for every generalized minor $\Delta$ on $G$. 
As part of  \thref{:frieze-pos}, which has  more related statements, we have the following result relating frieze points and Lusztig's total positivity.

\bth{:prin-frieze-intro}
A point $g \in \Lcc$ (resp. $g \in \Acc$) is a $\calS^{\rm prin}$ (resp. $\calS^{\rm triv}$)-frieze point if and only if 
$\Delta(g) \in \ZZ_{\geq 0}$
for every generalized minor $\Delta$ on $G$.
\eth

In $\S$\ref{ss:belt-finite}, we assume that the Cartan matrix $A$ is of finite type and indecomposable. For  $i \in [1, r]$, let $i^* \in [1, r]$ and $h(i; c) \geq 1$ be such that
$w_0 \omega_i = -\omega_{i^*}$ and 
$\omega_i > c\omega_i > \cdots  > c^{h(i; c)} \omega_i = -\omega_{i^*}$, where $w_0 \in W$ is the longest element. Define ${\bf F}:
[1,r] \times \ZZ \longrightarrow [1,r] \times \ZZ$ by 
\[
{\bf F}(i, m) = (i^*, \;m+h(i^*; c) +1).
\]
Using results by Yang and Zelevinsky \cite{YZ:Lcc}, we show that the generic frieze pattern associated to $A$ 
with arbitrary coefficients is $\bfF$-invariant and we determine its fundamental domain (\thref{:finite-main} and \coref{:uim-F-inv}). 
As a consequence, we prove the following statement which generalizes the gliding symmetry of Coxeter's frieze patterns
(see \reref{:ADE-F}).
\bth{:F-inv-all}
When the Cartan matrix $A$ is of finite type and indecomposable, all frieze patterns associated to $A$ with arbitrary coefficients are $\bfF$-invariant. 
\eth 
For all the cases except when $A$ is of type $A_n, D_{2n+1}, n \geq 2$, of $E_6$, one has $h(i; c) = \frac{1}{2}h$
for all $i \in [1, r]$, where $h$ is the Coxeter number of $A$ (see \reref{:F-w0=-1}). In the remaining cases, we explain in \leref{:hic-mi} how the
integers $h(i; c)$ can be computed using the quiver of the mutation matrix $\BA$.

\subsection{Notation}\label{ss:nota} For integers $m \leq n$, we use $[m, n]$ to denote the set of all integers $k$ such that $m \leq k \leq n$.
Elements in $\ZZ^n$ for any integer $n \geq 1$ are regarded as column vectors unless otherwise specified. The set of all $m \times n$ matrices with integral entries is denoted as 
$M_{m, n}(\ZZ)$.

For a real number $a$, let  $[a]_+={\rm max}(a, 0)$. 
For any matrix $M = (m_{i, j})$ with real entries, let $[M]_+ = ([m_{i,j}]_+)$, and write $M \geq 0$ if $M = [M]_+$.
The transpose of $M$ is denoted by $M^T$.

\medskip
\noindent
{\bf Notation on monomials}.
Suppose that $A$ is a commutative ring and let $A^\times$
be the group of units of $A$. For ${\bf a} = (a_1, \ldots, a_n) \in (A^\times)^n$ and
$L = (l_1, \ldots, l_n)^T \in \ZZ^n$,
we write 
\begin{equation}\label{eq:av}
{\bf a}^L = a_1^{l_1} \cdots a_n^{l_n} \in A^\times.
\end{equation}
If $L$ is an $n \times m$ integral matrix with columns $L_1, \ldots, L_m \in \ZZ^n$, we write
\begin{equation}\label{eq:aV}
{\bf a}^L = ({\bf a}^{L_1}, \ldots, {\bf a}^{L_m}) \in (A^\times)^m.
\end{equation}
For  ${\bf a} = (a_1, \ldots, a_n)$ and  ${\bf b} = (b_1, \ldots, b_n)$ in $(A^\times)^n$, we write 
\[
{\bf ab} = (a_1b_1, \ldots, a_nb_n) \in (A^\times)^n.
\]
As a non-standard notational convention, for ${\bf a} =(a_1, \ldots, a_n) \in A^n$ we also denote by ${\bf a}$ the set of its components, so 
$S \subset {\bf a}$
means that $S$ is a subset of $\{a_1, \ldots, a_n\}$.

\subsection{Acknowledgment} We would like to thank Dr. Peigen Cao  and Dr. Jie Pan for helpful discussions. We are grateful to M. Cuntz, T. Holm, and P. J{\o{}}rgensen for pointing out to us the reference \cite{Propp:integers}. The first and the third authors have been partially supported by
the Research Grants Council of the Hong Kong SAR, China (GRF 17307718 and GRF 17306621). The second author was partially supported by the Institute of Mathematical Research (IMR) of the University of Hong Kong during the academic year 2019 - 2020, and he has also been partially supported by NSF of China: 12101617 and Fundamental Research Funds for the Central Universities, Sun Yat-sen University (23qnpy51). 

 The main results from $\S$\ref{s:patterns} - $\S$\ref{s:Lcc} have been included in the PhD thesis of the first author from the University of Hong Kong \cite{Antoine:thesis}, which also extends the results in $\S$\ref{s:Lcc} to the Kac-Moody case and gives explicit descriptions of the cluster variables on the acyclic belt for the BFZ cluster structure on reduced double Bruhat cells in Kac-Moody groups.
 
\section{Preliminaries on cluster algebras}\label{s:A}

\subsection{Cluster algebras of geometric type}\label{ss:A}
We recall  the definition of cluster algebras of geometric type as in 
\cite{BFZ:III, FZ:IV}. In what follows, 
let $r \geq 1$ and $l \geq 0$ be integers, and let $\FF$ be  a field of rational functions over $\QQ$ in $r+l$ independent variables.

A {\it seed in $\FF$ of extended rank $(r, l)$} is a triple $\Sigma = ({\bf u}, {\bf p}, \widetilde{B})$, where

1) ${\bf u} = (u_1, \ldots, u_r)$ and ${\bf p} = (p_1, \ldots, p_{l})$ are such that $(u_1, \ldots, u_r, p_1, \ldots, p_{l})$ 
is a set of free generators of $\FF$ over $\QQ$;

2) $\widetilde{B}=  \left(\!\begin{array}{c} B\\ P\end{array}\!\right)$, also denoted as  $(B\bb P)$,
is an $(r+l) \times r$ integral matrix, called the {\it extended mutation matrix}, such that
$B$ formed by the first $r$ rows is skew-symmetrizable, i.e., $DB$ is skew-symmetric for some 
integral diagonal matrix $D$ with positive diagonal entries.
 The matrix $B$ is  
 called the {\it principal part} or the {\it mutation matrix} of $\widetilde{B}$ and $P$ the {\it coefficient part} or the {\it coefficient matrix} of $\widetilde{B}$.

The $r$-tuple ${\bf u}$ is called the {\it cluster} and the $(r+l)$-tuple $({\bf u}, {\bf p})$ the {\it extended cluster} of $\Sigma$. 
Elements in ${\bf u}$ are called {\it cluster variables}, elements in $({\bf u}, {\bf p})$ {\it extended cluster variables} and
elements in ${\bf p}$ {\it frozen variables} of $\Sigma$. 

Given a seed $\Sigma$ in $\FF$ as above,  for each $k \in [1, r]$, let $B_{k}$ and $P_{k}$ be the $k$th columns of $B$ and $P$ respectively, and define the
{\it $\widehat{y}$-variables} $\widehat{\bf y} = (\widehat{y}_1, \, \ldots, \, \widehat{y}_r)$ of $\Sigma$ (recall our notation on monomials in $\S$\ref{ss:nota}) by
\begin{equation}\label{eq:y-hat-de}
\widehat{\bf y} = \bfu^{B} \bfp^{P} = (\bfu^{B_{1}} \bfp^{P_{1}}, \, \ldots, \, \bfu^{B_{r}} \bfp^{P_{r}}).
\end{equation}
The mutation of $\Sigma$
in the direction of $k \in [1, r]$ is then the seed 
\[
\mu_k(\Sigma) = \left({\bf u}^\prime= (u_1^\prime, \ldots, u_r^\prime), \; {\bf p} = (p_1, \ldots, p_{l}), \; \widetilde{B}' = ({b}_{i j}^\prime)\right),
\]
where $u_j^\prime = u_j$ for $j \in [1, r]\backslash \{k\}$, 
\begin{equation}\label{eq:xk}
u_k^\prime =  \frac{1}{u_k}  \left( \bfu^{[B_{k}]_+}\bfp^{[P_{k}]_+}  +  \bfu^{[-B_{k}]_+}\bfp^{[-P_{k}]_+}\right)=  
\frac{1}{u_k} \bfu^{[-B_{k}]_+}\bfp^{[-P_{k}]_+}(1 + \widehat{y}_k),
\end{equation}
and, with $\widetilde{B} = (b_{i, j})_{i \in [1, r+l], j \in [1, r]}$, 
\begin{equation}\label{eq:B-mutation}
{b}_{i,j}^\prime = \begin{cases} -{b}_{i, j}, & \hs i = k \; \mbox{or} \; j = k, \\
{b}_{i,j} + [b_{i, k}]_+ [b_{k,j}]_+ - [-b_{i, k}]_+ [-b_{k,j}]_+, &\hs \mbox{otherwise}\end{cases}.
\end{equation}
One checks that the $\widehat{y}$-variables  $\widehat{\bf y}^{\, \prime} = (\widehat{y}_1^{\, \prime}, \, \ldots, \, \widehat{y}_r^{\, \prime})$ of the seed $\mu_k(\Sigma)$ 
are given by
\begin{equation}\label{eq:y-hat}
\widehat{y}_j^{\, \prime} = \begin{cases}\widehat{y}_k^{\, -1}, & \hs j = k, \\
\widehat{y}_j  \widehat{y}_k^{\,[b_{k, j}]_+} (1+\widehat{y}_k)^{-b_{k, j}}, & \hs j \neq k.\end{cases}
\end{equation}
In particular, the mutations of the $\widehat{y}$-variables depend only on the principal part of $\widetilde{B}$.
One has $\mu_k(\mu_k(\Sigma)) = \Sigma$ for every $k \in [1, r]$.

Let $\TT_r$ be the $r$-regular tree with the edges incident to every vertex labeled bijectively by the indices in $[1, r]$.
We write $t \in \TT_r$ for a vertex of $\TT_r$. A {\it  seed pattern} \cite{FZ:I, FZ:IV, Nakanishi-I}
with extended rank $(r, l)$ in $\FF$ is an assignment 
\[
\bsig = \{\Sigma_t =  ({\bf u}_t, {\bf p}, \widetilde{B}_t): t \in \TT_r\}
\]
of a seed of extended rank $(r, l)$ in $\FF$ to each $t \in \TT_r$ such that 
$\Sigma_{t'} = \mu_k (\Sigma_t)$ whenever 
$t \stackrel{k}{\mathdash} t'$.
A cluster (resp. an extended cluster) in any seed in $\bsig$ is called a cluster (resp. an extended cluster) of $\bsig$. 
All the seeds in $\calS$ have the same set of frozen variables, 
which we denote as\footnote{By the notational convention in $\S$\ref{ss:nota}, $\bfp$ also denotes the set $\{p_1, \ldots, p_l\}$.}
\[
{\rm Froz}(\calS) = \bfp,
\]
and we call 
elements in ${\rm Froz}(\calS)$ {\it frozen variables of $\calS$}. 
A  {\it cluster variable of $\calS$}  is by definition a cluster variable in any seed of $\calS$, and 
an {\it extended cluster variable} of $\calS$ is either  a cluster variable or a frozen variable of $\calS$.
We set
\begin{align*}
\fX(\bsig) = & \; \mbox{the set of all cluster variables of}\; \bsig,\\
\widetilde{\fX}(\bsig) =& \fX(\bsig) \sqcup {\rm Froz}(\bsig) = \; \mbox{the set of all extended cluster variables of}\; \bsig.
\end{align*}

\bde{:A-0}
{\rm 
Let $\bsig$ be a seed pattern in $\FF$.

1) For any
${\rm inv} \subset {\rm Froz}(\bsig)$, 
 the {\it $\ZZ$-cluster algebra} 
$\calA_{\ZZ}(\bsig, {\rm inv})$ is the $\ZZ$-subalgebra of $\FF$ generated by $\widetilde{\fX}(\bsig)$ and 
all  $p_j^{-1}$ for $p_j \in {\rm inv}$; 

2) For any field $\KK$ of characteristic $0$, define the {\it $\KK$-cluster algebra} 
\[
\calA_\KK (\bsig, {\rm inv})  = \calA_\mathbb{Z} (\bsig, {\rm inv}) \otimes_{\ZZ}\KK.
\]
Cluster  variables in $\bsig$ are also called cluster variables of $\calA_\KK (\bsig, {\rm inv})$, and similarly for frozen (resp. extended cluster) variables,
clusters, and extended clusters.

3) For $\KK = \ZZ$ or a field of characteristic $0$, we write $\calA_{\KK}(\bsig) = \calA_{\KK}(\bsig, {\rm inv})$ when ${\rm inv} = {\rm Froz}(\bsig)$.
\hfill $\diamond$
}
\ede

We will also frequently consider $\calA_\KK(\bsig, \emptyset)$, in which no frozen variables are invertible. 
For $\KK = \ZZ$ or a field of characteristic $0$, any cluster algebra $\calA_\KK = \calA_\KK(\bsig, {\rm inv})$ as in \deref{:A-0} 
is called a {\it cluster algebra of geometric type with extended rank $(r, l)$}.

\subsection{Tropicalization}\label{ss:trop}
Let $\bsig$ be a seed pattern of geometric type with extended rank $(r, l)$ and frozen variables ${\bf p} = (p_1, \ldots, p_{l})$. 
Using any cluster $\bfu$ of $\bsig$, define $\FF_{>0}(\bsig)$ to be the semi-field consisting of all non-zero elements in 
$\FF(\bsig) = \QQ(\bfu, \bfp)$ that have subtraction-free expressions in $(\bfu, \bfp)$. 
It is clear from the mutation rules that 
$\FF_{>0}(\bsig)$ is independent of the choice of the cluster $\bfu$. 
The semi-field $\FF_{>0}(\bsig)$ is called the {\it ambient semi-field} of the seed pattern $\bsig$.

Let $\PP = {\rm trop}(\bfp)$ be the tropical semi-field generated by $\bfp = (p_1, \ldots, p_{l})$. We denote by $1 \in \PP$ the identity element of $\PP$ with respect to the group operation.

\ble{:delta-PP}
    There exists a unique semi-field homomorphism 
    \[
        \FF_{>0}(\calS) \to \PP, \quad Z \mapsto Z\vert_\PP,
    \]
    which sends all cluster variables of $\calS$ to $1 \in \PP$ and $p_j \in \FF_{>0}(\calS)$ to $p_j \in \PP$ for each $j \in [1,l]$.
\ele
\begin{proof}
We prove existence by constructing a semi-field homomorphism with the prescribed properties. Fix any cluster ${\bf u} = (u_1, \dots, u_r)$ of $\calS$ and consider the semi-field homomorphism
\[
\FF_{>0} (\bsig) \longrightarrow {\PP}, \;\; Z \longmapsto Z|_{\bfu \to {\bf 1}, \PP},
\]
that sends $u_i \to 1 \in \PP$ for each $i \in [1, r]$ and $p_j \in \FF_{>0}(\bfu, \bfp)$ to $p_j \in {\PP}$ for each $j \in [1, l]$. Let $t_0 \in \TT_r$ be such that $\bfu$ is the cluster of $\bsig$ at $t_0$, and for $t \in \TT_r$, let
\[
\Sigma_t = (\bfu_t = (u_{t, 1}, \ldots, u_{t, r}), \; {\bf p}, \; \widetilde{B}_t)
\]
be the seed of $\bsig$ at $t \in \TT_r$. We
prove by induction on the distance of $t $ to $t_0$ that $u_{t, i}|_{\bfu \to {\bf 1}, \PP} = 1$ for all $t \in \TT_r$ and $i \in [1, r]$. 
The claim holds for $t = t_0$
by definition.
Suppose that the claim holds for
$t\in \TT_r$ and consider an edge $t \stackrel{k}{\mathdash} t'$ in $\TT_r$. By definition,
\[
u_{t', k} = \frac{1}{u_{t, k}} \left(\bfu_t^{[B_{t, k}]_+} \bfp^{[P_{t, k}]_+} + \bfu_t^{[-B_{t, k}]_+} \bfp^{[-P_{t, k}]_+}\right).
\]  
Using the induction assumption that $u_{t, i}|_{\bfu \to {\bf 1}, \PP} = 1$ for every $i \in [1, r]$,
we have
\begin{equation}\label{eq:pp}
u_{t', k}|_{\bfu \to {\bf 1}, \PP} = \bfp^{[P_{t, k}]_+} \oplus \bfp^{[-P_{t, k}]_+} = 1,
\end{equation}
where $\oplus$ is the addition in the semi-field $\PP = {\rm trop}(\bfp)$. 
Since $u_{t', i} = u_{t, i}$ for $i\neq k$, we have $u_{t', i}|_{\bfu \to {\bf 1}, \PP} = 1$ for every $i \in [1, r]$. This proves existence of the semi-field homomorphism $Z \mapsto Z\vert_\PP$. Uniqueness follows from the fact that any semi-field homomorphism 
$\FF_{>0} (\bsig) \rightarrow {\PP}$ is  uniquely
determined by its values on any extended cluster $({\bf u}, \bf {p})$.
\end{proof}

\bre{:any-PP}
{\rm The semi-field homomorphism $\FF_{>0}(\bsig) \to \PP$ in \leref{:delta-PP} 
can be similarly defined for a seed pattern $\bsig$ with coefficients in any semi-field
$\PP$, not necessarily a tropical semi-field, where \eqref{eq:pp} in the proof of \leref{:delta-PP} is due to the seeds in $\bsig$ being normalized \cite{FZ:IV}. 
\hfill $\diamond$
}
\ere

\bre{:trop-prin}
{\rm
The special case of \leref{:delta-PP} when $\bsig$ has principal coefficients at the cluster $\bfu$ 
has been observed in \cite[Proposition 4.14]{Nakanishi-I}, which, as
pointed out in \cite[Proposition 4.15]{Nakanishi-I}, does not imply that the $F$-polynomial of 
every cluster variable has constant term $1$. 
\hfill $\diamond$
}
\ere

\bco{:delta-pj}
For $j \in [1,l]$, one has a unique well-defined semi-field homomorphism
\begin{equation}\label{eq:delta-pj}
\FF_{>0}(\bsig) \longrightarrow {\rm trop}(p_j), \;\; Z \longmapsto Z|_{{\rm trop}(p_j)},
\end{equation}
that sends $p_j \in \FF_{>0}(\bsig)$ to $p_j \in {\rm trop}(p_j)$ and all the other extended cluster variables of $\bsig$ to $1 \in {\rm trop}(p_j)$.
\eco

\begin{proof} The composition of 
$\FF_{>0}(\bsig) \rightarrow \PP$ in \leref{:delta-PP} with the semi-field homomorphism
\[
\PP \longrightarrow {\rm trop}(p_j), \; \;p_k \longmapsto \begin{cases} p_j, & \hs k = j; \\
1, & \hs k \in [1, l]\backslash\{j\},\end{cases}
\]
is the necessarily unique
semi-field homomorphism from $\FF_{>0}(\bsig)$ to ${\rm trop}(p_j)$ as desired.
\end{proof}




\subsection{Positive Laurent phenomenon for cluster algebras of geometric type}\label{ss:pos-Laurent}
Let again 
\[
\bsig = \{\Sigma_t =({\bf u}_t,\, \bfp = (p_1, \ldots, p_l), \,\widetilde{B}_t): \,t \in \TT_r\}
\]
be a seed pattern of geometric type with extended rank 
$(r, l)$. Recall \cite{BFZ:III} that the upper cluster algebra of $\bsig$ is defined as
\[
\calU_\ZZ(\bsig) = \bigcap_{t \in \TT_r} \ZZ[{\bf u}_t^{\pm 1}, \bfp^{\pm 1}] \subset \FF(\bsig), 
\]
and that a non-zero $Z \in \FF(\bsig)$ is said to be {\it universally positive Laurent} if it belongs to
\[
\calU^+_\ZZ(\bsig) = \bigcap_{t \in \TT_r} \ZZ_{\geq 0} [{\bf u}_t^{\pm 1}, \bfp^{\pm 1}].
\]
For any ${\rm inv} \subset {\rm Froz}(\bsig) = \bfp$, introduce
\[
\calU_\ZZ(\calS, {\rm inv}) = \bigcap_{t \in \TT_r} \ZZ [{\bf u}_t^{\pm 1}, \bfp][p_j^{-1}: p_j \in {\rm inv}] \quad \text{and} \quad \calU^+_\ZZ(\bsig, {\rm inv}) = \bigcap_{t \in \TT_r} \ZZ_{\geq 0} [{\bf u}_t^{\pm 1}, \bfp][p_j^{-1}: p_j \in {\rm inv}].
\]
Then $\calU^+_\ZZ(\bsig, \emptyset) = \bigcap_{t \in \TT_r} \ZZ_{\geq 0} [{\bf u}_t^{\pm 1}, \bfp]$, $\calU^+_\ZZ(\bsig, \bfp) = \calU^+_\ZZ(\bsig)$,
and for every ${\rm inv} \subset \bfp = {\rm Froz}(\calS)$,
\[
 \calU^+_\ZZ(\bsig, \emptyset) \subset \calU^+_\ZZ(\bsig, {\rm inv}) 
\subset \calU^+_\ZZ(\bsig).
\]
Recall that $\fX(\bsig)$ denotes the set of all cluster variables of $\bsig$. We refer to the following \thref{:pos-Laurent} as the 
{\it Positive Laurent Phenomenon for cluster algebras of geometric type}.

\bth{:pos-Laurent} \cite[Proposition 11.2]{FZ:II} and \cite[Theorem 4.10]{GHKK}
For every seed pattern $\bsig$ of geometric type, one has $\fX(\bsig) \subset \calU^+_\ZZ(\bsig, \emptyset)$.
\eth

\bre{:1}
{\rm
By \leref{:delta-PP} one has
$u|_{\PP} = 1$ for every $u \in \fX(\bsig)$, which 
is a stronger statement than $u$ being polynomial in ${\bf p}$.
\hfill $\diamond$
}
\ere

Recall that a {\it cluster monomial} of $\bsig$ is an element of $\FF(\bsig)$ of the form ${\bf u}_t^\alpha$ for some  $t \in \TT_r$ and 
$\alpha\in (\ZZ_{\geq 0})^r$. By \thref{:pos-Laurent}, all cluster monomials of $\bsig$ are in 
$\calU^+_\ZZ(\bsig, \emptyset) \subset \calU^+_\ZZ(\bsig, {\rm inv})$ for any ${\rm inv} \subset  \bfp$.
Recall that $\bsig$ is said to be of {\it finite type} if it has finitely many cluster variables \cite{FZ:II}. The following \prref{:finite-atomic} will
be used in $\S$\ref{ss:rooted}.

\bpr{:finite-atomic}
If $\bsig$ is of finite type, then for any ${\rm inv} \subset {\bf p}$, every $Z \in \calU^+_\ZZ(\bsig, {\rm inv})$ is a unique finite linear combination of
cluster monomials with coefficients in ${\mathbb{A}}^+_{\rm inv} \stackrel{\rm def}{=}\ZZ_{\geq 0}[\bfp][p_j^{-1}: j \in {\rm inv}]$.
\epr

\begin{proof}
When ${\rm inv} = \bfp$, the statement is equivalent to the fact that cluster monomials form a basis of $\calU^+_\ZZ(\bsig)$  over $\ZZ_{\geq 0}[{\bf p}^{\pm 1}]$ (see, e.g., \cite[Theorem 10.2]{FT:orbifold-bases} or 
\cite[Corollary 3.7]{Huang:proper}). 

Let ${\rm inv} \subset \bfp$ be arbitrary, and let $Z \in \calU^+_\ZZ(\bsig, {\rm inv}) \subset \calU^+_\ZZ(\bsig)$. Then $Z$ is a unique 
sum
\[
Z = \sum_{i=1}^m c_i M_i,
\]
where each $c_i \in  \ZZ_{\geq 0} [\bfp^{\pm 1}]$ and $M_i$ a cluster monomial of $\bsig$.   
Let $j \in [1, l]$ be such that $p_j \notin {\rm inv}$. By \coref{:delta-pj}, $M_i|_{{\rm trop}(p_j)} = 1$ for every $i \in [1, m]$. Thus
\[
Z|_{{\rm trop}(p_j)} = \oplus_{i=1}^m c_i|_{{\rm trop}(p_j)},
\]
where $\oplus$ is the addition in the semi-field ${\rm trop}(p_j)$.
Writing $Z|_{{\rm trop}(p_j)} = p_j^{\delta_j(Z)}$ with $\delta_j(Z) \in \ZZ$, 
it follows from $Z \in \calU^+_\ZZ(\bsig, {\rm inv})$ that $\delta_j(Z) \geq 0$.
Thus every $c_i \in \ZZ_{\geq 0} [\bfp^{\pm 1}]$ must be polynomial in $p_j$. Since $p_j \notin {\rm inv}$ is arbitrary, 
$c_i \in {\mathbb{A}}^+_{\rm inv}$ for every $i \in [1, m]$.
\end{proof}

\section{Friezes, frieze testing criteria, and frieze points}\label{s:friezes}
\subsection{Friezes}\label{ss:friezes} 
 Let $\bsig$ be a seed pattern of geometric type, and 
recall from \deref{:A-0} the $\KK$-cluster algebra $\calA_{\KK}(\calS, {\rm inv})$ associated to $\bsig$ for any ${\rm inv} \subset
{\rm Froz}(\bsig)$, where $\KK = \ZZ$ or is a field of characteristic $0$. Recall also that
$\widetilde{\fX}(\bsig)$ denotes the set of all extended cluster variables of $\calS$. The following definition follows  \cite{F:non-zero, FP:Dn, GSch}.

\bde{:frieze-0} 
{\rm 
1) By a {\it frieze} of $\calA_\ZZ(\bsig, \emptyset)$ we mean any ring homomorphism
$h: \calA_{\ZZ}(\calS, \emptyset) \to \ZZ$ that takes positive integral values on every 
$u \in \widetilde{\fX}(\bsig)$; 

2) For a field $\KK$ of characteristic $0$, a {\it frieze} of $\calA_\KK(\bsig, {\rm inv})$ is any $\KK$-algebra homomorphism
$h: \calA_{\KK}(\calS, {\rm inv}) \to \KK$ that takes positive integral values on every $u \in \widetilde{\fX}(\bsig)$;

3) Any frieze as in 1) or 2) is called an $\bsig$-frieze.
\hfill $\diamond$
}
\ede


\bre{:ZZ-QQ}
{\rm It is clear from the definitions that every $\bsig$-frieze restricts to, and is determined by, a frieze of $\calA_\ZZ(\bsig, \emptyset)$.
The freedom of allowing arbitrary fields $\KK$ of characteristics $0$, especially for $\KK = \QQ$ or $\CC$, makes it easier to 
consider friezes when ${\rm inv} \neq \emptyset$ (see
$\S$\ref{ss:testing}) and  
 renders possible geometrical considerations (see $\S$\ref{ss:friezes-V}). 
 \hfill $\diamond$
}
\ere

The Positive Laurent Phenomenon for cluster algebras of geometric type in \thref{:pos-Laurent}
ensures the existence of a special class of $\bsig$-friezes.

\bld{:unitary-frieze} 
Let $\bsig$ be a seed pattern of geometric type with extended rank $(r, l)$, and let 
$\wdu = (u_1, \ldots, u_r, p_1, \ldots, p_{l})$ be any extended cluster in $\bsig$. For any ${\bfm} = (m_1, \ldots, m_{l}) \in (\ZZ_{>0})^{l}$,
one has a unique frieze  $h_{\wdu, {\bfm}}: \calA_\ZZ(\bsig, \emptyset) \to \ZZ$ such that
 \begin{equation}\label{eq:F-dwx}
h_{\wdu, {\bfm}}(u_i)  = 1, \;\;\; i \in [1, r]  \hs \mbox{and}\hs 
h_{\wdu, {\bfm}}(p_j) = m_j, \,\;\; j \in [1, l].
\end{equation}
Friezes obtained this way are said to be \textit{unitary}.
\eld

\begin{proof}
Let $T_{\wdu} = \ZZ[u_1^{\pm 1}, \ldots, u_r^{\pm 1}, p_1, \ldots, p_{l}]$ and 
$T_{\wdu}^+ = \ZZ_{\geq 0} [u_1^{\pm 1}, \ldots, u_r^{\pm 1}, p_1, \ldots, p_{l}]$.
Then each ${\bfm} = (m_1, \ldots, m_{l}) \in (\ZZ_{>0})^{l}$ gives rise to a 
$\ZZ$-algebra homomorphism $h_{\wdu, {\bfm}}: T_{\wdu} \to \ZZ$ via \eqref{eq:F-dwx}. 
By \thref{:pos-Laurent},  
$\widetilde{\fX}(\calS) \subset T_{\wdu}^+$, so $h_{\wdu, {\bfm}}$ restricts to a frieze of $\calA_\ZZ(\bsig, \emptyset)$.
\end{proof}

\subsection{Frieze testing criteria}\label{ss:testing}
Let again $\calS$ be a seed pattern with frozen variables $\bfp$.

\bde{:testing-criterion}
{\rm A subset $\calT$ of $\widetilde{\fX}(\bsig)$
 is called an  {\it $\bsig$-frieze testing set} if for any ${\rm inv} \subset \bfp$, a $\QQ$-algebra homomorphism 
$h: \calA_\QQ(\bsig, {\rm inv})\to \QQ$ is a frieze if and only if $h(u)  \in \ZZ_{>0}$ for every $u \in \calT$.
\hfill $\diamond$
}
\ede

In examples, when testing an $h: \calA_\QQ(\bsig,  {\rm inv}) \to \QQ$ using a subset $\calT\subset \widetilde{\fX}(\bsig)$, 
we often can modify the condition that $h(u) \in \ZZ_{>0}$ for every $u \in \calT$ by dividing $\calT$ into $\calT = \calT_1 \cup \calT_2$ 
and imposing integrality of $h$ on $\calT_1$  and positivity of $h$ on $\calT_2$. 
We refer to all such conditions as {\it frieze testing criteria}.
We illustrate this by the following simple observation in a special case.  

\ble{:poly-0}
Suppose that $\calA_\ZZ(\bsig, \emptyset)$ is generated over $\ZZ$ by a subset ${\mathcal{Z}}$ of $\widetilde{\fX}(\bsig)$
and let
$(u_1, \ldots, u_{r+l})$ be any extended cluster of $\bsig$. Then for any ${\rm inv} \subset \bfp$, a $\QQ$-algebra homomorphism
$h: \calA_\QQ(\bsig, {\rm inv}) \to \QQ$ is a frieze if and only if the following two conditions hold:

1) $h(z) \in \ZZ$ for every $z \in {\mathcal{Z}}$, and 

2) $h(u_j) >0$ for every $j \in [1, r+l]$.
\ele

\begin{proof}
Let $h: \calA_\QQ(\bsig, {\rm inv}) \to \QQ$ be a $\QQ$-algebra homomorphism satisfying 1) and 2).
Let $u\in \widetilde{\fX}(\bsig)$.
Condition 1) implies that  $h(u) \in \ZZ$, and Condition 2) implies that  $h(u) >0$ by the
positivity of the Laurent expression of $u$ in $u_1, \ldots, u_{r+l}$. Thus $h$ is a frieze.
\end{proof}

\bre{:finite-generated}
{\rm
If
$\calA_\ZZ(\bsig, \emptyset)$ is finitely generated as a ring, it is generated by 
a finite subset of $\widetilde{\fX}(\bsig)$, and \leref{:poly-0} yields a finite $\calS$-frieze testing set. 
\hfill $\diamond$
}
\ere

\bex{:poly}
{\rm
A particular case of \leref{:poly-0} is when $\calA_\ZZ(\bsig, \emptyset)$ is the polynomial ring of a subset 
$\{z_1, \ldots, z_m\}$ of $\fX(\calS)$.
A systematic class of such examples comes from the theory of symmetric Poisson CGL extensions developed by K. Goodearl and M. Yakimov \cite{GY:CGL, GY:Poi-CGL}. 
For a rich class of examples from Lie theory, see \cite{Elek-L:BS,  GLS:Kac-Moody, GSW:augmentation}.  
\hfill $\diamond$
}
\eex

\subsection{Frieze points in cluster varieties}\label{ss:friezes-V}
For cluster algebras that are finitely generated, we take a geometrical view point  on their  friezes  via the notion of cluster varieties. 

Let $V$ be an $n$-dimensional irreducible rational affine variety over $\CC$, let $\CC(V)$ be the field of rational functions on $V$, and let 
$\CC[V] \subset \CC(V)$ be the algebra of regular functions on $V$. The following definition of cluster structures on $V$ follows \cite[Proposition 3.37]{GSV:book}.

\bde{:cluster-V}
{\rm
1) A {\it cluster structure on $V$} is a seed pattern $\bsig$ in $\CC(V)$ such that 
\[
\calA_\CC (\bsig, {\rm inv}) = \CC[V]
\]
for some ${\rm inv} \subset {\rm Froz}(\bsig)$. 
A seed $\Sigma$ in $\CC(V)$ is said to {\it define a cluster structure on $V$} if it is a seed in a seed pattern that is a cluster structure on $V$;

2) For a cluster structure $\bsig$ on $V$, by an {\it $\bsig$-frieze point in $V$} we mean a point $v \in V$ such that all extended cluster variables in $\bsig$
take positive integral values at $v$.
\hfill $\diamond$
}
\ede

For a cluster structure $\calS$ on $V$, the subset ${\rm inv}$ of ${\rm Froz}(\bsig)$ such that $\calA_\CC (\bsig, {\rm inv}) = \CC[V]$
is necessarily unique, as we now explain.

\ble{:invV}
Let $\bsig$ be a seed pattern in $\CC(V)$. 

1) If there exists ${\rm inv} \subset {\rm Froz}(\bsig)$ such that  $\calA_\CC (\bsig, {\rm inv}) = \CC[V]$, then 
${\rm Froz}(\bsig) \subset \CC[V]$ and
\[
{\rm inv} = {\rm inv}(V) \stackrel{\rm def}{=} \{p \in {\rm Froz}(\bsig): p\;  \mbox{vanishes nowhere on}\; V\}.
\]

2) The seed pattern $\bsig$ is a cluster structure on $V$ if and only if 
\[
{\rm Froz}(\bsig) \subset \CC[V] \hs \mbox{and} \hs \calA_\CC (\bsig, {\rm inv}(V)) = \CC[V].
\]
\ele

\begin{proof}
1) Let ${\rm inv} \subset {\rm Froz}(\bsig)$ be such that  $\calA_\CC (\bsig, {\rm inv}) = \CC[V]$. By definition
\[
{\rm Froz}(\bsig) \cap \{p^{-1}: p \in {\rm inv}\} \subset \calA_\CC (\bsig, {\rm inv}).
\]
Thus ${\rm Froz}(\bsig) \subset \CC[V]$ and ${\rm inv} \subset {\rm inv}(V)$. 
To show that ${\rm inv}(V) \subset {\rm inv}$, let
\[
\widetilde{{\bf u}} = (u_1, \ldots, u_r, p_1, \ldots, p_{l})
\]
 be any extended cluster in $\bsig$, and assume without
loss of generality that ${\rm inv} = \{p_1, \ldots, p_a\}$ for some $1 \leq a \leq l$. By the Laurent phenomenon, 
\begin{equation}\label{eq:CV-Laurent}
\CC[V]  = \calA_\CC(\bsig, {\rm inv}) \subset \CC[u_1^{\pm 1}, \ldots, u_r^{\pm 1}, p_1^{\pm 1}, \ldots, p_{a}^{\pm 1}, p_{a+1}, \ldots, p_{l}].
\end{equation}
If $j \in [1, l]$ is such that $p_j \in {\rm inv}(V)$, then $p_j^{-1} \in \CC[V]$. It follows that 
$p_j \in  {\rm inv}$.  This shows that ${\rm inv} = {\rm inv}(V)$.
2) follows trivially from 1).
\end{proof}


Let $\calS$ be a cluster structure on $V$ with  ${\rm Froz}(\bsig) = \{p_1, \ldots, p_{l}\}$, and let 
${\rm inv}(V) = \{p_1, \ldots, p_a\}$ for some $0 \leq a \leq l$. For any extended cluster
$\wdu  = (u_1, \ldots, u_r, p_1, \ldots, p_{l})$ in $\bsig$, the inclusion in \eqref{eq:CV-Laurent}
gives rise to a biregular isomorphism
\[
\rho_{\widetilde{\bf u}}:\;\; (\CC^\times)^{r+a} \times \CC^{l-a} \longrightarrow V_{\widetilde{\bf u}}:=  \{v \in V: u_1(v) \cdots u_r(v) p_1(v) \cdots p_{a}(v) \neq 0\} \subset V,
\]
whose inverse is given by
\[
\widetilde{\bf u}:\;\; V_{\widetilde{\bf u}} \longrightarrow (\CC^\times)^{r+a} \times \CC^{l-a}, \;\; 
v \longmapsto (u_1(v), \ldots,  u_r(v),  p_1(v),  \ldots,  p_{l}(v)).
\]
We call the above map $\wdu$ 
the {\it $\bsig$-cluster coordinate chart} on $V$ defined by  $\wdu$. 
By definition, we then have the following description of $\bsig$-frieze points.

\ble{:f-points}
Given a cluster structure $\bsig$ on $V$,  the ${\bsig}$-frieze points in $V$ are precisely all those $v \in V$ that 
have positive integral coordinates in every 
$\bsig$-cluster coordinate chart in $V$. 
\ele

\bre{:inv}
{\rm When $\calA_\ZZ(\bsig, \emptyset)$
is finitely generated, we can identify $\bsig$-friezes  with
$\bsig$-frieze points in $V =  {\rm Spec}(\calA_\CC(\bsig, {\rm inv}))$ for any ${\rm inv} \subset {\rm Froz}(\bsig)$. This is 
especially helpful if we have good geometrical models for ${\rm Spec}(\calA_\CC(\bsig, {\rm inv}))$ for some choices of
${\rm inv}$. 
\hfill $\diamond$
}
\ere

\section{The category ${\bf Frieze}$ and pullbacks of friezes}\label{s:pullbacks}
\subsection{Property $\mF$ and the category ${\bf Frieze}$}\label{ss:mF}
For a seed pattern 
$\bsig$ of geometric type, 
recall again that 
$\widetilde{\mathfrak{X}}(\bsig)$ denotes the set of all extended cluster variables of $\bsig$. Let 
\[
\calA_{\ZZ}^+(\bsig, \emptyset) \subset \calA_\ZZ(\bsig, \emptyset)
\]
 be the semi-ring (with $1$ but without $0$) generated by $\widetilde{\fX}(\bsig)$.

\bde{:mF}
{\rm
For two seed patterns $\bsig$ and $\pbsig$  of geometric type, a ring homomorphism
$\psi: \calA_\ZZ(\bsig, \emptyset) \rightarrow \calA_\ZZ(\pbsig, \emptyset)$ is said to have {\it Property $\mathfrak{F}$} if
\begin{equation}\label{eq:key}
\psi\left(\widetilde{\fX}(\bsig)\right) \subset \calA_{\ZZ}^+ (\pbsig, \emptyset), \hs \mbox{equivalently} \hs
\psi\left(\calA_{\ZZ}^+ (\bsig, \emptyset)\right) \subset \calA_{\ZZ}^+ (\pbsig, \emptyset).
\end{equation}
\hfill $\diamond$
}
\ede

Property $\mF$ is clearly preserved under compositions of ring homomorphisms. 

\bde{:frieze-cat}
{\rm
By ${\bf Frieze}$ we mean the category, whose objects are $\ZZ$-cluster algebras of geometric type in which frozen variables are not invertible, 
and whose morphisms are ring homomorphisms with Property $\mF$. 
\hfill $\diamond$
}
\ede

We also regard $\ZZ$ as an object in
${\bf Frieze}$, i.e., the $\ZZ$-cluster algebra of rank $0$ with trivial coefficients.
A frieze as in \deref{:frieze-0} is nothing but a morphism to $\ZZ$ in the category ${\bf Frieze}$. 
For two objects $\calA_\ZZ$ and 
$\calA_\ZZ^\prime$ in
${\bf Frieze}$, a morphism from $\calA_\ZZ$ to $\calA_\ZZ^\prime$ in ${\bf Frieze}$ can then be regarded as an $\calA_\ZZ^\prime$-valued frieze of $\calA_\ZZ$.
We make the following explicit statement for future reference.

\ble{:pullback-obvious}
If $\psi: \calA_\ZZ(\bsig, \emptyset) \rightarrow \calA_\ZZ(\pbsig, \emptyset)$ is a ring homomorphism
with Property $\mF$, 
 then for every frieze 
$h: \calA_\ZZ(\pbsig, \,\emptyset) \to \ZZ$, the composition
$h \circ \psi: \calA_\ZZ(\bsig, \,\emptyset) \rightarrow \ZZ$ 
 is again a frieze. 
\ele

In the context of \leref{:pullback-obvious}, we will call $h \circ \psi$ the
{\it pullback of the frieze $h$ by $\psi$}.

Recall that $\calA_\QQ(\bsig)$ denotes the $\QQ$-cluster algebra defined by $\bsig$ with all frozen variables invertible.
Extending \deref{:mF}, we say that a $\QQ$-algebra homomorphism 
$\psi: \calA_\QQ(\bsig) \rightarrow \calA_\QQ(\pbsig)$
 has {\it Property $\mathfrak{F}$} if it induces a ring homomorphism $\calA_\ZZ(\bsig, \emptyset) \rightarrow \calA_\ZZ(\pbsig, \emptyset)$ 
with {Property $\mathfrak{F}$}.  
We now give examples of ring (algebra) homomorphisms between cluster algebras that have Property $\mF$.

\subsection{Freezing and deletion of cluster variables}\label{ss:rooted}  Consider a seed pattern
\[
\bsig = \{\Sigma_t = (\bfu_t = (u_{t, 1}, \ldots, u_{t, r}), \,\bfp, \,\widetilde{B}_t): t \in \TT_r\}
\]
of geometric type  with extended rank $(r, l)$. Fix $t_0 \in \TT_r$
and write ${\bf u}_{t_0} = (u_1, \ldots, u_r)$.

For an integer $0 \leq s < r$, consider first a seed pattern 
\begin{equation}\label{eq:ssq}
\bsig^\ssq = \{\Sigma^\ssq_t = (\bfu^\ssq_t, \,\bfp^\ssq, \,\widetilde{B}^\ssq_t): t \in \TT_s\}
\end{equation}
of extended rank $(s, \,r-s+l)$ obtained by {\it freezing} 
$r-s$ cluster variables of $\Sigma_{t_0}$, which, for simplicity, are assumed to be
$u_{s+1}, \ldots, u_{r}$. More precisely, identifying $\TT_s$ with the sub-graph of
$\TT_r$ generated by  $t_0$ and all the edges labeled by elements in $[1, s]$, we have $\FF(\bsig^\ssq) = \FF(\bsig) = \QQ({\bf u}_{t_0}, \bfp)$,  
\[
{\bf u}^\ssq_{t_0} =(u_{1}, \, \ldots, \, u_{s}), \hs \bfp^\ssq = (u_{s+1}, \, \ldots, \, u_{r}, \, p_1, \ldots, p_l),
\]
and $\widetilde{B}^\ssq_{t_0}$ is obtained from $\widetilde{B}_{t_0}$ by deleting the {\it columns} of $\widetilde{{B}}_{t_0}$ indexed by $j \in [s+1, r]$.
Then $\widetilde{\fX}(\bsig^\ssq) \subset \widetilde{\fX}(\bsig)$ and $\calA_\ZZ(\bsig^\ssq, \emptyset)$ is a sub-ring of $\calA_\ZZ(\bsig, \emptyset)$. Consequently, we have

\ble{:freezing-mF}
With notation as above, the inclusion
$\calA_\ZZ(\bsig^\ssq, \emptyset) \hookrightarrow \calA_\ZZ(\bsig, \emptyset)$
 has Property $\mF$.
\ele

Again for an integer $0 \leq s < r$, we consider next the seed pattern 
\[
\bsig^\dag = \{\Sigma^\dag_t = (\bfu^\dag_t, \,\bfp^\dag, \,\widetilde{B}^\dag_t): t \in \TT_s\}
\]
of extended rank $(s, l)$ obtained by {\it deleting} $r-s$ cluster variables of $\Sigma_{t_0}$, again assumed to be
$u_{s+1}, \ldots, u_{r}$. More precisely, the mutation matrix $\widetilde{B}^\dag_{t_0}$ is obtained from 
$\widetilde{B}_{t_0}$ by deleting the {\it rows and columns} of $\widetilde{B}_{t_0}$ indexed by $j \in [s+1, r]$, and we write the
extended cluster of $\Sigma^\dag_{t_0}$ as
\[
{\bf u}^\dag_{t_0} =(u^\dag_{1}, \, \ldots, \, u^\dag_{s}), \hs \bfp^\dag = (p^\dag_1, \, \ldots,\,  p^\dag_l).
\]
Let $\psi: \ZZ[u_1^{\pm 1}, \ldots, u_r^{\pm 1},p_1, \ldots, p_l] \to \ZZ[(u_1^\dag)^{\pm 1}, \ldots, (u_s^\dag)^{\pm 1}, p_1^\dag, \ldots, p_l^\dag]$ be the ring homomorphism
given by
\[
\psi(u_{i}) = u^\dag_{i}, \hs \psi(u_j) = 1, \hs \psi(p_{j'}) = p^\dag_{j'}, \hs i \in [1, s], \, j \in [s+1, r], \, j' \in [1, l].
\]
Let $\bsig^\ssq$ be the seed pattern obtained from $\bsig$ by freezing the variables $u_{s+1}, \ldots, u_{r}$ as in \eqref{eq:ssq}.
Then $\psi$ restricts to a ring homomorphism 
$\calA_\ZZ(\bsig^\ssq, \emptyset) \longrightarrow \calA_\ZZ(\dbsig, \emptyset)$
 which induces a bijection 
from $\fX(\bsig^\ssq)$ to $\fX(\dbsig)$ 
In particular, if $\dbsig$ is of finite type, so is $\bsig^\ssq$. The following proposition 
is adapted from \cite[Proposition 2.9]{GM:finite} to include cluster algebras with coefficients. 

\bpr{:deleting}
With notation as above, if $\dbsig$ is of finite type, then the ring homomorphism \[
\psi: \ZZ[u_1^{\pm 1}, \ldots, u_r^{\pm 1},p_1, \ldots, p_l] \to \ZZ[(u_1^\dag)^{\pm 1}, \ldots, (u_s^\dag)^{\pm 1}, p_1^\dag, \ldots, p_l^\dag]
\]
restricts to a ring homomorphism
$\psi: \calA_\ZZ(\bsig, \emptyset) \to \calA_\ZZ(\dbsig, \emptyset)$ that has Property $\mF$.
\epr

\begin{proof}
It suffices to show that $\psi(\fX(\bsig)) \subset \calA^+_{\ZZ}(\dbsig, \emptyset)$. 

Consider again the seed pattern $\bsig^\ssq$ in $\FF(\bsig^\ssq) = \FF(\bsig)$ given in \eqref{eq:ssq} and let
\[
{\rm inv}^\ssq = (u_{s+1}, \ldots, u_r) 
\subset \bfp^\ssq = (u_{s+1}, \ldots, u_r, p_1, \ldots, p_l).
\]
 In the notation of $\S$\ref{ss:pos-Laurent}, we then have
\[
\calU^+_{\ZZ}(\bsig^\ssq, {\rm inv}^\ssq) = \bigcap_{t \in \TT_s} \ZZ_{\geq 0}[{\bf u}_t^{\pm 1}, \bfp].
\]
By the positive
Laurent phenomenon for seed patterns of geometric type in \thref{:pos-Laurent}, 
\[
\fX(\bsig) \subset \bigcap_{t \in \TT_r} \ZZ_{\geq 0}[{\bf u}_t^{\pm 1}, \bfp] \subset \calU^+_{\ZZ}(\bsig^\ssq, {\rm inv}^\ssq).
\]
By \prref{:finite-atomic} applied to $\calU^+_{\ZZ}(\bsig^\ssq, {\rm inv}^\ssq)$, every 
 $u \in \fX(\bsig)$ is a unique sum
\[
u = \sum_{i=1}^m c_i M_i^\ssq,
\]
where each $c_i \in \ZZ_{\geq 0} [u_{s+1}^{\pm 1}, \ldots, u_r^{\pm 1}, \bfp]$, $c_i \neq 0$,  and $M_i^\ssq$ a cluster monomial of $\bsig^\ssq$.
Thus
\[
\psi(u) = \sum_{i=1}^m  \psi(c_i) \psi(M_i^\ssq) \in \calA^+_{\ZZ}(\dbsig, \emptyset).
\]
\end{proof}


\bre{:rooted}
{\rm
In the terminology of \cite{ADS:rooted}, the ring homomorphism $\psi: \calA_\ZZ(\bsig, \emptyset) \to \calA_\ZZ(\dbsig, \emptyset)$ 
in \prref{:deleting} is a {\it specialization} and is, 
by \cite[Proposition 6.9]{ADS:rooted}, an example of a surjective {\it ideal rooted cluster morphism}.
In general,  $\psi: \calA_\ZZ(\bsig, \emptyset) \to \calA_\ZZ(\bsig^\dag, \emptyset)$ in \prref{:deleting} 
sends cluster monomials of $\bsig$ to linear combinations of cluster monomials of $\bsig^\dag$. Indeed, in the 
simplest example of $\bsig$ having trivial coefficients and
\[
\Sigma_{t_0} = \left((u_1, u_2), \, B =\left(\begin{array}{cc} 0 & 1 \\ -1 & 0\end{array}\right)\right)
\]
and $s = 1$, so that $\Sigma_{t_0}^\dag = (u_1^\dag, (0))$, one has
\[
\psi\left(u_1, \; u_2, \; \frac{1+u_2}{u_1}, \; \frac{1+u_1}{u_2}, \; \frac{1+u_1+u_2}{u_1u_2}\right) = 
\left(u_1^\dag, \; 1,  \; \frac{2}{u_1^\dag}, \; 1 + u_1^\dag, \; 1 + \frac{2}{u_1^\dag}\right).
\]
Examples when $\bsig$ is of type $D_4$ are given in \cite[Example 6.19]{ADS:rooted}. 
\hfill $\diamond$
}
\ere

\subsection{Homomorphisms of coefficient rescaling type}\label{ss:rescaling}
For an arbitrary seed pattern $\calS$ of geometric type with extended rank $(r, l)$ and frozen variables $\bfp  = (p_1, \ldots, p_l)$, 
we say that two non-zero elements $z, z' \in \FF(\calS)$ are {\it proportional} if 
there exists $v \in \ZZ^l$ such that $z = z' \bfp^v$.

Consider now
two seed patterns   $\bsig$ and $\obsig$
of geometric type with  respective extended ranks $(r, l)$ and $(r, \ol)$ and frozen variables $\bfp$ and $\obfp$.
Recall  again that $\calA_\QQ(\bsig)$ and $\calA_\QQ(\obsig)$ are the respective
cluster algebras in which all frozen variables are invertible. 

\bde{:rescaling-type-00}
{\rm 
A $\QQ$-algebra homomorphism $\psi: \calA_\QQ(\bsig) \to \calA_\QQ(\obsig)$ is said to be of {\it coefficient rescaling type} if there 
exists $E \in M_{\ol, l}(\ZZ)$ such that $\psi(\bfp) = \obfp^E$, and if for every $ u\in \fX(\calS)$, there exists $\overline{u} \in \fX({\ocalS})$ such that
$\psi(u)$ and $\ou$ are proportional.
In such a case, for every 
$u \in \fX(\calS)$, the cluster variable $\overline{u} \in \fX({\ocalS})$  is uniquely determined by $u$ and
$\psi$ via
\[
\overline{u} =\frac{\psi(u)}{\psi(u)|_{\overline{\PP}}}.
\]
Consequently one has the  well-defined map 
\[
\fX_{\psi}: \;\; \fX(\calS) \longrightarrow \fX(\overline{\calS}), \;\; u \longmapsto \overline{u}.
\]
For
$u \in \fX(\bsig)$, the unique $R_u \in \ZZ^{\overline{l}}$ such that
$\psi(u)|_{\overline{\PP}} = \obfp^{R_u}
$
is called the {\it coefficient rescaling vector} of $u$. We say that $\psi$ rescales $u$ {\it trivially} if $R_u = 0$.
\hfill $\diamond$}
\ede

Recall that for an integral matrix $M$, we write $M \geq 0$ if all of its entries are non-negative.
\ble{:rescaling-mF}
Let  $\psi: \calA_\QQ(\bsig) \to \calA_\QQ(\obsig)$ be a $\QQ$-algebra homomorphism of coefficient rescaling type such that
$\psi(\bfp) = \obfp^E$, where $E \in M_{\ol, l}(\ZZ)$. 

1) $\psi$ has Property $\mF$ 
if and only if $E \geq 0$ and $R_u \geq 0$ for every $u \in \fX(\bsig)$;

2) If $E \geq 0$ and if there exists a cluster ${\bf u}$ of $\calS$ such that $R_u = 0$ for every cluster variable in ${\bf u}$, 
then $R_u \geq 0$ for every $u \in \fX(\bsig)$, and thus $\psi$ has  Property $\mF$.
\ele

\begin{proof}
1) follows from the definitions. 
To prove 2), 
write  $\psi(\bfu) = \obfu$ (not assuming that $\obfu$ is a cluster in $\obsig$).
Let $u \in \fX(\bsig)$ be arbitrary and write
$u = U(\bfu, \bfp) \in \ZZ_{\geq 0}[\bfu^{\pm 1}, \bfp]$. Then 
\[
\psi(u) = U(\obfu, \obfp^E).
\]
 Since $E \geq 0$, $\psi(u)$ is polynomial in $\obfp$. As $\obfp^{R_u} = \psi(u)|_{\overline{\PP}}$, we have 
$R_u \geq 0$. 
\end{proof}

In the remainder of $\S$\ref{ss:rescaling}, we present our main examples of  homomorphisms of coefficient rescaling type, namely 
{\it cluster algebra quasi-homomorphisms} (for cluster algebras of geometric type), and we give examples of such homomorphisms with Property $\mF$.

Consider two seed patterns 
\begin{equation}\label{eq:bsig-obsig}
\bsig = \{\Sigma_t = (\bfu_t, \,\bfp, \,(B_t \bb P_t)): t \in \TT_r\} \hs \mbox{and} \hs 
\obsig = \{\overline{S}_t = (\overline{\bfu}_t,\, \overline{\bfp}, \, (\overline{B}_t \bb \oP_t)): t \in \TT_r\}
\end{equation}
of geometric type  with  respective extended ranks $(r, l)$ and $(r, \ol)$. 
Let $I_r \in M_{r, r}(\ZZ)$ be the identity matrix. 
The following definition of {\it cluster algebra quasi-homomorphisms}, which is more suitable for our study of friezes, 
is based on \cite[Remark 3.1 and Corollary 4.5]{Fraser:quasi}. 

\bde{:psi-quasi}
{\rm A {\it cluster algebra quasi-homomorphism from $\calA_\QQ(\bsig)$ to $\calA_\QQ(\obsig)$} is a $\QQ$-algebra homomorphism
$\psi: \calA_\QQ(\bsig) \rightarrow \calA_\QQ(\obsig)$ such that 
\begin{equation}\label{eq:Psi-00}
\psi(\bfu_{t_0}, \, \bfp) =( \obfu_{\ot_0} \, \obfp^R, \, \obfp^E)
\end{equation}
for some
$t_0, \ot_0 \in \TT_t$, $R \in M_{\ol, r}(\ZZ)$ and $E \in M_{\ol, l}(\ZZ)$ satisfying 
$B_{t_0} = \overline{B}_{\ot_0}$ and 
\begin{equation}\label{eq:RE}
RB_{t_0} + EP_{t_0} = \overline{P}_{\ot_0}.
\end{equation}
\hfill $\diamond$
}
\ede

In what follows, by a {\it tree automorphism} of 
$\TT_r$ we mean
a bijection $t \mapsto \ot$ on the set of vertices of $\TT_r$ such that $\ot \stackrel{k}{\mathdash} \overline{t'}$ 
whenever $t \stackrel{k}{\mathdash} t'$.
Note that for any given $t_0, \ot_0 \in \TT_r$ there is a unique tree automorphism of $\TT_r$ sending $t_0$ to $\ot_0$. 
If $\{A_t: t \in \TT_r\}$ is a collection of matrices in $M_{n, r}(\ZZ)$ for some integer $n$, we denote by
$A_{t, k}$ the $k$th column of $A_t$ for every $t \in \TT_r$ and $k \in [1, r]$.

\bpr{:psi-quasi-all}
1) For any $t_0, \ot_0 \in \TT_0$ such that $B_{t_0}=\overline{B}_{\ot_0}$, 
every solution $(R, E)$ to
\eqref{eq:RE} 
gives rise to a unique cluster algebra quasi-homomorphism $\psi:  \calA_\QQ(\bsig) \to \calA_\QQ(\obsig)$ satisfying \eqref{eq:Psi-00}.

2) For any cluster algebra quasi-homomorphism $\psi: \calA_\QQ(\bsig) \rightarrow \calA_\QQ(\obsig)$, there exist $E \in M_{\ol, l}(\ZZ)$,
a tree automorphism $\TT_r \to \TT_r, t \mapsto \ot$, and $\{R_t \in M_{\ol, r}(\ZZ): t \in \TT_r\}$, such that for every $t \in \TT_r$,
\begin{equation}\label{eq:B-psi}
B_t = \overline{B}_{\ot}, \hs R_t B_t + E P_t = \overline{P}_{\ot}, \hs \mbox{and} \hs \psi({\bf u}_t, \, \bfp) = (\overline{\bf u}_{\ot} \,\obfp^{R_t}, \, \obfp^E).
\end{equation}
Furthermore, for
$k \in [1, r]$ and $t \stackrel{k}{\mathdash} t'$ in $\TT_r$, one has
$R_{t, j} = R_{t', j}$ for $j \neq k$, and 
\begin{equation}\label{eq:Rt-k-00}
R_{t, k} + R_{t', k} = R_t [B_{t, k}]_+ +E [P_{t, k}]_+ -[\overline{P}_{\ot,k}]_+ = R_t [-B_{t, k}]_+ +E [-P_{t, k}]_+ -[-\overline{P}_{\ot,k}]_+.
\end{equation}
\epr

\begin{proof}
1) Let $t_0, \ot_0 \in \TT_r$ and $(R, E)$ be as given, and let 
$\Psi: \FF_{>0}(\calS) \to \FF_{>0}(\overline{\calS})$ be the semi-field homomorphism uniquely defined by 
$\Psi(\bfu_{t_0}, \, \bfp) =(\obfu_{\ot_0} \, \obfp^R, \, \obfp^E)$. 
Let $\TT_r\to \TT_r, t \mapsto \ot,$ be the tree automorphism sending $t_0$ to $\ot_0$. 
It follows from $B_{t_0} = \overline{B}_{\ot_0}$ that $B_{t} = \overline{B}_{\ot}$
for every $t \in \TT_r$. Recall that the $\widehat{y}$-variables of $\Sigma_t$ and $\overline{\Sigma}_{\ot}$
are respectively given by
 $\widehat{{\bf y}}_t = \bfu_t^{B_t} \bfp^{P_t}$ and $\widehat{\overline{{\bf y}}}_{\ot} = 
\obfu_t^{B_t} \obfp^{\overline{P}_{\ot}}$. 
The condition $B_{t_0} = \overline{B}_{\ot_0}$ and \eqref{eq:RE} imply that   $\Psi(\widehat{{\bf y}}_{t_0})  = \widehat{\overline{{\bf y}}}_{\ot_0}$. By
the mutation rules of the $\widehat{y}$-variables, one has 
$\Psi(\widehat{{\bf y}}_t) = \widehat{\overline{{\bf y}}}_{\ot}$ for every $t \in \TT_r$.
It then follows from 
the mutation rule \eqref{eq:xk} that for every $t \in \TT_r$ one has $\Psi(\bfu_t) = \obfu_{\ot} \, \obfp^{R_t}$ for a unique
$R_t \in M_{\overline{l}, r}(\ZZ)$. 
It also follows from $\Psi(\widehat{{\bf y}}_t) = \widehat{\overline{{\bf y}}}_{\ot}$ that 
$R_t B_t + E P_t = \overline{P}_{\ot}$ for every $t \in \TT_r$.

Consider the unique $\QQ$-algebra homomorphism
\[
\psi_{t_0}: \;\; \QQ[\bfu_{t_0}^{\pm 1}, \bfp^{\pm 1}]  \longrightarrow \QQ[\obfu_{\ot_0}^{\,\pm 1},  \obfp^{\, \pm 1}]
\]
determined by 
$\psi_{t_0}(\bfu_{t_0}, \bfp) = \Psi(\bfu_{t_0}, \bfp) = (\obfu_{\ot_0}\,\obfp^R, \, \obfp^E)$.
Then $\psi_{t_0}$ and $\Psi$ agree on
$\ZZ_{\geq 0}[\bfu_{t_0}^{\pm 1}, \bfp]$. In particular, for any cluster variable $u$ of $\bsig$ one has
$\psi_{t_0}(u) = \Psi(u) \in \calA_\QQ(\obsig)$. 
Set
\[
\psi: \;\; \calA_\QQ(\bsig) \longrightarrow \calA_\QQ(\obsig), \; u \longmapsto \psi_{t_0}(u).
\]
Then $\psi$ is the unique cluster algebra quasi-homomorphism 
satisfying \eqref{eq:Psi-00}.

2) Let $\psi: \calA_\QQ(\bsig) \rightarrow \calA_\QQ(\obsig)$ be a cluster algebra quasi-homomorphism with $t_0, \ot_0$ and $(R, E)$ as in \deref{:psi-quasi}, so we have the tree automorphism $\TT_r \to \TT_r, t \mapsto \ot,$ sending $t_0$ to $\ot_0$. Let
$\Psi: \FF_{>0}(\calS) \to \FF_{>0}(\overline{\calS})$ be again the semi-field homomorphism defined by 
$\Psi(\bfu_{t_0}, \, \bfp) =(\obfu_{\ot_0} \, \obfp^R, \, \obfp^E)$. We have seen in the proof of 1) that 
there exists a unique $R_t \in M_{\overline{l}, r}(\ZZ)$ for each $t \in \TT_r$ such that 
\eqref{eq:B-psi} holds and that $\psi(u) = \Psi(u)$ for very cluster variable $u$ of $\calS$. Let $k \in [1, r]$ and $t \stackrel{k}{\mathdash} t'$ in $\TT_r$. Then $R_{t', j} = R_{t, j}$ for $j \in [1, r]\backslash\{k\}$.
Applying $\Psi$  to the mutation relation
\begin{equation}\label{eq:uu-0}
u_{t, k} u_{t', k} =\bfu_t^{[B_{t, k}]_+} \bfp^{[P_{t, k}]_+} \left(1+ \why_{t, k}^{\,-1}\right),
\end{equation}
and using the facts that $B_t = \overline{B}_{\ot}$ and $\Psi(\why_{t, k}) = \widehat{\overline{y}}_{\ot, k}$, one has
\begin{align*}
\Psi(u_{t', k}) &=\overline{u}_{\ot, k}^{\,-1} \,\obfu_{\ot}^{[B_{t, k}]_+} \obfp^{\,-R_{t, k} + R_t{[B_{t, k}]_+} + E[P_{t, k}]_+} 
\left(1+\widehat{\overline{y}}_{\ot, k}^{\,-1}\right)\\
& = \obfp^{\,-R_{t, k} +R_t[B_{t, k}]_+ + E[P_{t, k}]_+-[\overline{P}_{\ot,k}]_+} \,
\overline{u}_{\ot, k}^{\,-1}\, \obfu_{\ot}^{[\overline{B}_{\ot, k}]_+}\obfp^{[\overline{P}_{\ot,k}]_+} \left(1+\widehat{\overline{y}}_{\ot, k}^{\,-1}\right)\\
& =\overline{u}_{\overline{t'}, k} \,  \obfp^{\,-R_{t, k} +R_t[B_{t, k}]_+ + E[P_{t, k}]_+-[\overline{P}_{\ot,k}]_+},
\end{align*}
where in the last step we used that fact that $\ot \stackrel{k}{\mathdash} \overline{t'}$. Thus
\[
R_{t, k} + R_{t', k} = R_t[B_{t, k}]_+ + E[P_{t, k}]_+-[\overline{P}_{\ot,k}]_+.
\]
 Since $R_t B_t + E P_t = \overline{P}_{\ot}$, one also the second 
identity in \eqref{eq:Rt-k-00}. 
\end{proof}


\bco{:quasi-rescaling}
Every cluster algebra quasi-homomorphism 
is of coefficient rescaling type.
\eco

We now have the following  construction of cluster algebra quasi-homomorphisms with Property $\frakF$ as a consequence of  \leref{:rescaling-mF}.

\bco{:quasi-frakF}
Let $t_0, \ot_0 \in \TT_0$ be such that $B_{t_0}=\overline{B}_{\ot_0}$, and assume that  $(0, E)$ is a solution to \eqref{eq:RE}, i.e.,
$EP_{t_0} = \overline{P}_{\ot_0}$. Then the cluster algebra quasi-homomorphism $\psi:  \calA_\QQ(\bsig) \to \calA_\QQ(\obsig)$ determined by 
$\psi(\bfu_{t_0}, \, \bfp) =( \obfu_{\ot_0}, \, \obfp^E)$ has Property $\frakF$ if and only if $E \geq 0$.
\eco

 Denote the set of all 
integral row vectors of size $r$ by $\ZZ^r_{\rm row}$.
Note then that there exists $E \geq 0$ such that $EP_{t_0} = \overline{P}_{\ot_0}$ if and only if
every row of $\overline{P}_{\ot_0}$ lies in the cone in $\ZZ^r_{\rm row}$ generated by the rows of $P_{t_0}$.

\bex{:P0P0}
{\rm
Let $\bsig = \{{\Sigma}_t = ({\bfu}_t,\, {\bfp}, \, ({B}_t \bb P_t)): t \in \TT_r\}$ 
be any seed pattern of geometric type with extended rank $(r, l)$ and let $t_0 \in \TT_r$ be arbitrary. Let 
$P_0 \in GL(r, \ZZ)$, and let $\bsig^{P_0}$ be the unique seed pattern of extended rank $(r, r)$ with seed
$({\bf x}_{t_0}, {\bf q}, ({B}_{t_0} \bb P_0))$ at $t_0$.
Taking $\ot_0 = t_0$, and $E = {P}_{t_0} P_{0}^{-1}$ in \coref{:quasi-frakF}, if 
${P}_{t_0} P_{0}^{-1} \geq 0$,  then the  
cluster algebra quasi-homomorphism $\psi: \calA_\QQ(\bsig^{P_0}) \to \calA_\QQ(\bsig)$ determined by 
$\psi({\bf x}_{t_0}, \, {\bf q}) = ({\bf u}_{t_0}, \, \bfp^{{P}_{t_0} P_{0}^{-1}})$
has Property $\mF$.

When $P_{0} = I_r$,  we denote  $\bsig^{P_0}$ by $\bsig^{\rm prin}$ for {\it principal coefficients} and denote the cluster of $\bsig^{\rm prin}$ at
$t \in \TT_r$ by ${\bf x}_t = (x_{t, 1}, \ldots, x_{t, r})$. For the cluster algebra quasi-homomorphism
$\psi: \calA_\QQ(\bsig^{\rm prin}) \to \calA_\QQ(\bsig)$ defined by $R = 0$ and $E = {P}_{t_0}$, one has 
\[
{u}_{t, i} = \frac{\psi(x_{t, i})}{\psi(x_{t, i})|_{{\PP}}},\hs t \in \TT_r, \;\; i \in [1, r],
\]
which is precisely the Separation Formula \cite[Corollary 6.3]{FZ:IV} (this fact was already observed in \cite[Remark 4.6]{Fraser:quasi}).
 If ${P}_{t_0} \geq 0$, then
$\psi: \calA_\QQ(\bsig^{\rm prin}) \to \calA_\QQ(\bsig)$ has Property $\mF$, so $\psi$ pulls back every frieze of $\calA_\QQ(\bsig)$ to one of
$\calA_\QQ(\bsig^{\rm prin})$. 
\hfill $\diamond$
}
\eex

We finish this section with a discussion on coefficient specializations for cluster algebras of geometric type. We follow the definition 
given by N. Reading \cite{Reading:univ-coeff}.

\bde{:specialization} \cite[Definition 3.1]{Reading:univ-coeff} 
{\rm
Let $\calS$ and $\overline{\calS}$ be as in \eqref{eq:bsig-obsig}. 
A $\QQ$-algebra homomorphism $\psi: \calA_\QQ(\bsig)\to \calA_\QQ(\obsig)$ is called a {\it coefficient specialization} if 

1) $\psi(\bfp) = \obfp^E$ for some $E \in M_{\ol, l}(\ZZ)$;

2) there exist $t_0, \ot_0 \in \TT_r$ such that $B_{t_0} = \overline{B}_{\ot_0}$, $\phi({\bf u}_{t_0}) = \overline{{\bf u}}_{\ot_0}$, and
\begin{equation}\label{eq:tt-2}
E P_t = \overline{P}_{\ot} \hs \mbox{and} \hs E [P_t]_+ = [\overline{P}_{\ot}]_+
\end{equation}
for all $t \in \TT_r$, where $\TT_r \to \TT_r, t \mapsto \ot,$ is the tree automorphism that sends $t_0$ to $\ot_0$.
\hfill $\diamond$
}
\ede

\bre{:specialization}
{\rm By \eqref{eq:Rt-k-00},  a $\QQ$-algebra homomorphism $\psi: \calA_\QQ(\bsig)\to \calA_\QQ(\obsig)$ 
is a coefficient specialization if and only if $\psi$ 
is a cluster algebra quasi-homomorphism that rescales every cluster variable of $\bsig$ trivially. 
As noted in \cite[Remark 3.4]{Reading:univ-coeff}, \deref{:specialization}  is stronger than that given in 
\cite[Definition 12.1]{FZ:IV}.
\hfill $\diamond$
}
\ere

We now turn to {\it universal coefficient specializations} for cluster algebras of finite  type.

Fix $t_0 \in \TT_r$. For a  skew-symmetrizable $B \in M_{r, r}(\ZZ)$ and any $P \in M_{l, r}(\ZZ)$, let $\bsig(B\bb P)$ be a seed pattern 
with extended mutation matrix $(B \bb P)$ at $t_0$, and let $\calA_\QQ(B\bb P) = \calA_\QQ(\bsig(B\bb P))$.

Assume that $B \in M_{r, r}(\ZZ)$ is a skew-symmetrizable matrix of  {\it finite type}, so that $\calA_\QQ(B\bb P)$ has finitely many cluster variables
 for  every $l \in \ZZ_{\geq 0}$ and every $P \in M_{l, r}(\ZZ)$. Let $({\bf u}_{t_0}, \, \bfp)$ denote the extended cluster of $\bsig(B \bb P)$ at $t_0$. 
Consider the transpose $B^T$ of $B$ which is also of finite type, and denote the set of all cluster variables of $\calA_\QQ(B^T \bb I_r)$ by
\[
\fX= \fX(\bsig(B^T\bb I_r)).
\] 
For $x \in \fX$, let $g_x\in \ZZ^r$ be the $g$-vector of $x$ with respect to the seed of $\bsig(B^T\bb I_r)$ at $t_0$ as defined in \cite[$\S$6]{FZ:IV}. 
Since $g_x$ is a column vector by definition, we set $\rho_x = g_x^T \in \ZZ^r_{\rm row}$.

Let $P^{\rm univ} \in M_{\fX, r}(\ZZ)$ be the matrix with the row $\rho_x$ for 
$x \in \fX$, and consider $\widetilde{B}^{\rm univ} = (B\bb P^{\rm univ})$ as an extended mutation matrix. Here and in what follows, for two finite sets $I$ and $J$ we
denote by   $M_{I, J}(\ZZ)$ the set of all matrices whose rows and columns are respectively indexed by $I$ and $J$, and we also set 
$M_{I, r}(\ZZ) = M_{I, [1, r]}(\ZZ)$ and $M_{l, J} (\ZZ)= M_{[1, l], J}(\ZZ)$. 
Multiplication of such matrices is carried out in the natural way.\footnote{One could identify $\fX$ with $[1, |\fX|]$ using any listing of
elements in $\fX$, but we prefer not to do so as there is no such natural listing.} Let $({\bf u}_{t_0}^{\rm univ}, \, \bfp^{\rm univ})$ denote the extended cluster of $\bsig(B \bb P^{\rm univ})$ at $t_0$. 
By \cite[Theorem 10.12]{Reading:univ-coeff} and \cite[Corollary 5.5]{GHKK}, for any $P \in M_{l, r}(\ZZ)$ there exists a unique coefficient specialization
\begin{equation}\label{eq:psi-univ}
\psi: \;\; \calA_\QQ(B \bb P^{\rm univ}) \longrightarrow \calA_\QQ(B \bb P),
\end{equation}
called the {\it universal coefficient specialization}, such that $\psi({\bf u}_{t_0}^{\rm univ}) = {\bf u}_{t_0}$.

\bth{:finite-mF} Assume that the mutation matrix $B$ is of finite type. The
 universal coefficient specialization $\psi$ in \eqref{eq:psi-univ} has Property $\mF$ for every $P \in M_{l, r}(\ZZ)$. Thus $\psi$ pulls back friezes of $\calA_\QQ(B \bb P)$ to friezes of
 $\calA_\QQ(B \bb P^{\rm univ})$.
\eth

\begin{proof} Let $E\in M_{l, \fX}(\ZZ)$ be such that $\psi({\bf p}^{\rm univ}) = {\bf p}^E$. By \coref{:quasi-frakF} it is enough to show that $E \geq 0$. As $B^T$ is of finite type, the integral $g$-vector fan in $\ZZ^r$ defined by $B^T$ is complete. This fact is reformulated in 
\cite[Theorem 10.12 and Definition 6.1]{Reading:univ-coeff} as saying that 
the set $\{\rho_x: x \in  \fX\}$ is a {\it positive $\ZZ$-basis for $B$}, in the sense that for any $\rho \in \ZZ^r_{\rm row}$ there exists a unique 
non-negative $E_\rho \in M_{1, \fX}(\ZZ)$ such that 
\[
E_\rho P^{\rm univ} = \rho \hs \mbox{and} \hs E_\rho [P^{\rm univ}]_+ = [\rho]_+.
\]
If   $\rho_1, \ldots, \rho_l$ are the rows of $P$, \cite[Remark 4.5]{Reading:univ-coeff} implies that $E$ has rows $E_{\rho_1}, \ldots, E_{\rho_l}$. Since $E_{\rho_i} \geq 0$ for every $i \in [1, r]$,  we deduce that $E \geq 0$. 
\end{proof}

\section{Frieze patterns with coefficients}\label{s:patterns}
\subsection{The setup}\label{ss:setup} Throughout $\S$\ref{s:patterns}, we fix an $r \times r$  symmetrizable Cartan matrix $A =(a_{i, j})$, i.e., 

1) $a_{i, i} = 2$ and $a_{i, j} \in \ZZ$ for all $i, j \in [1, r]$;

2) $a_{i, j} \leq 0$ and $a_{i, j} = 0$ if and only if $a_{j, i} =0$ for all $i, j \in [1, r]$ and $i \neq j$;

3) there exists a diagonal matrix $D$ with positive integers on the diagonal such that $DA$ is symmetric. 
When $DA$  is positive definite, $A$ is called a Cartan matrix (of finite type).

Recall the associated skew-symmetrizable matrix $\BA$ associated to $A$ given in
\eqref{eq:Bo-00}.
 We also fix $P \in M_{l, r}(\ZZ)$ where $l \geq 0$ is
arbitrary. 
Let $t_0 \in \TT_r$, and let $\bfp = (p_1, \ldots, p_{l})$.  In this section, we consider the seed pattern
\begin{equation}\label{eq:sig}
\bsig (\BA\bb P)= \{\Sigma_t = (\bfu_t,\; \bfp, \; \widetilde{B}_t = (B_t\bb P_t)): \; t \in \TT_r\}
\end{equation}
with $\widetilde{B}_{t_0} = (\BA\bb P)$. 
We introduce {\it frieze patterns associated to $A$ with coefficient matrix $P$} and discuss their relations with friezes of the cluster algebra $\calA_\QQ(\BA\bb P) \stackrel{\rm def}{=} \calA_\QQ(\bsig(\BA\bb P))$.


\subsection{The acyclic belt and cluster additive functions}\label{ss:belt} Recall that we have introduced in $\S$\ref{ss:patterns-intro} the  2-regular sub-tree $\TT_r^\flat$
of $\TT_r$ given as
\[ 
\cdots  \stackrel{r-1}{\mathdash} t(r,\!-2) \stackrel{r}{\mathdash} t(1, \!-1) \stackrel{1}{\mathdash} \cdots
\stackrel{r-1}{\mathdash} t(r,\!-1) \stackrel{r}{\mathdash} t(1, \!0) \stackrel{1}{\mathdash} 
\cdots \stackrel{r-1}{\mathdash} t(r,\! 0) \stackrel{r}{\mathdash} t(1,\! 1),
\stackrel{1}{\mathdash} \cdots  
\]
where $t(1, 0) = t_0$. Correspondingly we have  the seeds
\[
\Sigma_{t(i, m)} = ({\bf u}_{t(i, m)}, \, \bfp, \, (B_{t(i, m)}\bb P_{t(i, m)})), \hs (i,m) \in [1, r] \times \ZZ.
\]
Note that for every $(i,m) \in [1, r] \times \ZZ$, one has
\begin{equation}\label{eq:mui}
\mu_i \Sigma_{t(i,m)} =\begin{cases} \Sigma_{t(i+1, m)}, & \hs i \in [1, r-1], \\
\Sigma_{t(1, m+1)}, & \hs i = r,\end{cases}
\end{equation}
so $\Sigma_{t(i, m+1)} = \mu_{i-1} \cdots \mu_2 \mu_1 \mu_r \cdots \mu_{i+1} \mu_i \Sigma_{t(i, m)}$.  The following \leref{:B-im-column} is a direct consequence of \eqref{eq:mui} and the mutation rule \eqref{eq:B-mutation}. In particular, each
$\mu_i$ in \eqref{eq:mui} is a source mutation.

\ble{:B-im-column}
For every $(i, m) \in [1, r] \times \ZZ$, 

1) $B_{t(i, m+1)} = B_{t(i, m)}$; In particular, $B_{t(1, m)} = B$ for every $m \in \ZZ$;

2) the $i$th column of  
$B_{t(i, m)}$ 
is $(a_{1,i}, \ldots, a_{i-1,i}, 0, a_{i+1,i}, \ldots, a_{r,i})^T\leq 0$.
\ele

\bre{:acyclic-belt}
    {\rm Every seed on the acyclic belt through $\Sigma_{t(1,0)}$ is acyclic, in the sense of \cite[Definition 1.14]{BFZ:III}. In general however, not every acyclic seed in $\calS(B_A \bb P)$ appears on the acyclic belt through $\Sigma_{t(1,0)}$ (this can already be seen when $A$ is of type $A_3$). 
    \hfill $\diamond$
    }
\ere

\bnota{:uim}
{\rm
For $(i, m) \in [1, r] \times \ZZ$, let $u(i, m)$ be the $i$th variable of the
cluster ${\bf u}_{t(i, m)}$.
\hfill $\diamond$
}
\enota

\ble{:umi-mutation}
For all $(i, m) \in [1, r] \times \ZZ$ one has
\begin{equation}\label{eq:u-m-0}
{\bf u}_{t(i,m)} = (u(1, m+1), \ldots, u(i-1, m+1), \; u(i,m),\; u(i+1, m), \ldots , u(r, m)),
\end{equation}
and the mutation of $\Sigma_{t(i, m)}$ in direction $i$ gives the mutation relation
\begin{equation}\label{eq:ex-0}
	u(i,m)\, u(m+1,i) =  {\bf p}^{[p(i,m)]_+} + {\bf p}^{[-p(i,m)]_+}\prod_{j = i+1}^r u(m,j)^{-a_{j i}} \prod_{j = 1}^{i-1} u(m+1,j)^{-a_{j i}},
\end{equation}
where $p(i, m)$ is the $i$th column of $P_{t(i, m)}$. 
\ele

\begin{proof}
\eqref{eq:u-m-0} 
follows directly from \eqref{eq:mui}, and \eqref{eq:ex-0} follows from 2) of \leref{:B-im-column}.
\end{proof}

\bre{:knitting}
{\rm
The above process of assigning the cluster variable $u(i,m)$ to  $(i, m) \in [1, r] \times \ZZ$ is called the {\it knitting algorithm} 
(\cite[$\S$2.2]{Keller:06} and  \cite[$\S$2]{Keller:Sche:linear}).
\hfill $\diamond$
}
\ere

We now prove a distinguished  property of the
vectors $\{p(i, m): (i, m) \in [1, r] \times \ZZ\}$ that appear in \eqref{eq:ex-0}. 
We first recall a notion 
introduced by C. M. Ringel \cite{Ringel}. 

\bde{:additive}
{\rm 
For any integer $l\geq 1$, a {\it $\ZZ^{l}$-valued cluster additive  function associated to the symmetrizable generalized Cartan matrix $A$} is any map $q:
[1, r] \times \ZZ \to \ZZ^l$ satisfying
\begin{equation}\label{eq:pmi-additive}
q(i,m) + q(i, m+1) = \sum_{j=i+1}^r (-a_{j,i})[q(j,m)]_+ + \sum_{j=1}^{i-1}(-a_{j,i}) [q(j, m+1)]_+, \hs  (i, m) \in [1, r] \times \ZZ.
\end{equation}
\hfill $\diamond$
}
\ede

\bpr{:additive} For every $P \in M_{l, r}(\ZZ)$, the map 
\[
p: \;\; [1, r] \times \ZZ \longrightarrow \ZZ^{l}, \;\; (i, m) \longmapsto p(i, m) \in \ZZ^l
\]
that appears in \eqref{eq:ex-0} is a $\ZZ^{l}$-valued additive cluster function  associated to $A$. Every 
every $\ZZ^{l}$-valued cluster additive function associated to $A$ arises this way for a unique $P \in M_{l, r}(\ZZ)$.
\epr
\begin{proof}
Fix $m \in \ZZ$.
For $i\in [1, r]$, let $P(i, m)$ be the $i$th column of $P_{t(1, m)}$. 
  Using the facts that $B_{t(1, m)} = \BA$ and
$\Sigma_{t(i, m)} = \mu_{i-1} \cdots \mu_1 \Sigma_{t(1,m)}$, and using  the mutation formula \eqref{eq:B-mutation}, one computes  $p(i, m)$ directly and obtains
\begin{equation}\label{eq:pmi-1}
P(i,m) = p(i,m) + \sum_{j=1}^{i-1} a_{j,i} [p(j,m)]_+.
\end{equation}
Similarly, using 
$\Sigma_{t(1,m+1)} =\mu_r \cdots  \mu_{i+1} \mu_i \Sigma_{t(i, m)}$, one obtains
\begin{equation}\label{eq:pmi-2}
P(i, m+1) = -p(i,m) + \sum_{j=i+1}^r (-a_{j,i}) [p(j,m)]_+.
\end{equation}
Replacing $m$ by $m+1$ in \eqref{eq:pmi-1} and comparing with \eqref{eq:pmi-2}, one  arrives at 
\[
p(i,m) + p(i, m+1) = \sum_{j=i+1}^r (-a_{j,i})[p(j,m)]_+ + \sum_{j=1}^{i-1}(-a_{j,i}) [p(j, m+1)]_+.
\]
Any $\ZZ^{l}$-valued cluster additive cluster function associated to $A$ is determined by its values on $\{(i, 0): i \in [1, r]\}$, which corresponds to a unique
$P \in M_{l, r}(\ZZ)$ via \eqref{eq:pmi-1} for $m = 0$.
\end{proof}

\bre{:additive}
{\rm
It is clear from the definition that a $\ZZ^l$-valued function on $[1, r] \times \ZZ$ is cluster additive if and only if every one of its component functions is a
$\ZZ$-valued cluster additive function, which is what is defined in \cite{Ringel}. Cluster additive functions are studied in more detail in
\cite{CGL:additive} from the point view of tropical points and the Fock-Goncharov Duality.  
\hfill $\diamond$
}
\ere

\subsection{Frieze patterns with coefficients}\label{ss:patterns} Continuing with the setup in $\S$\ref{ss:belt},
let
$\calA_\ZZ(\BA\bb P, \emptyset) = \calA_\ZZ(\calS(\BA\bb P), \emptyset)$.
Suppose that $h:\calA_\ZZ(\BA\bb P, \emptyset) \to \ZZ$ is a frieze. Setting
\begin{equation}\label{eq:F-f}
f(i, m) = h(u(i, m)), \hs (i, m) \in [1, r] \times \ZZ,
\end{equation}
and applying $h$ to both sides of \eqref{eq:ex-0}, one gets, for all $(i, m) \in [1, r] \times \ZZ$,
\begin{equation}\label{eq:ex-f}
	f(i,m)\, f(i, m+1) =  f({\bf p})^{[p(i,m)]_+} + f({\bf p})^{[-p(i,m)]_+}\prod_{j = i+1}^r f(j, m)^{-a_{j,i}} \prod_{j = 1}^{i-1} f(j, m+1)^{-a_{j,i}}.
\end{equation}
Recall from \deref{:pattern-intro} that we call any map 
$f:  ([1, r] \times \ZZ)   \cup \{p_1, \ldots, p_{l}\}  \rightarrow \ZZ_{>0}$
satisfying \eqref{eq:ex-f} for all $(i, m) \in [1, r] \times \ZZ$ a {\it frieze pattern associated to $A$ with coefficient matrix $P$}. 
We also call such an $f$ a {\it frieze pattern associated to $(\BA \bb P)$}. 

\bre{:repetition-quiver}
{\rm
Let $Q_{\BA}$ be the valued quiver associated to $\BA$ (\cite[$\S$3.3]{Keller:derived}). 
Frieze patterns associated to $(\BA \bb P)$ can be interpreted as $\ZZ_{>0}$-valued
{\it weighted multiplicative functions} on the repetition quiver $\ZZ Q_{\BA}$.
See \cite[$\S$2]{Sophie-M:survey} and \cite[$\S$6.2 and 6.3]{GM:finite} for the case when $P = 0$.
\hfill $\diamond$
}
\ere

\bde{:generic-frieze}
{\rm 
We call both maps
\begin{align}\label{eq:generic-variable}
&[1, r] \times \ZZ \longrightarrow \fX(\calS(\BA\bb P)), \;\; (i, m) \longmapsto u(i, m),\\
 \label{eq:generic-algebra}
&([1, r] \times \ZZ)    \cup \{p_1, \ldots, p_{l}\} \longrightarrow \calA_\ZZ(\BA\bb P, \emptyset), \;\; (i, m) \longmapsto u(i, m), \; p_j \longmapsto p_j,
\end{align}
 the {\it generic frieze pattern associated to $A$ with coefficient matrix $P$ or associated to $(\BA\bb P)$}. 
\hfill $\diamond$
}
\ede

\bpr{:ARS-same}
Suppose that $\calA_\ZZ(\BA\bb P, \emptyset)$ is generated over $\ZZ[\bfp]$ by a subset ${\mathcal{Z}}$ of
$\{u(i, m): (i, m) \in  [1, r] \times \ZZ\}$. Then
 every frieze pattern $f$ associated to $(\BA \bb P)$ comes from a unique frieze $h: \calA_\ZZ(\BA \bb P, \emptyset) \to \ZZ$  via evaluation as in \eqref{eq:F-f}.
\epr

\begin{proof}
Let  $f: ([1, r] \times \ZZ)   \cup \{p_1, \ldots, p_{l}\}  \to \ZZ_{>0}$ be a frieze pattern associated to $(\BA \bb P)$. Let
$\widetilde{\bfu} = (u(1,0), \ldots, u(r,0), p_1, \ldots, p_{l})$, and let
$T_{\widetilde{\bfu}} = \QQ[u(1,0)^{\pm 1}, \ldots, u(r,0)^{\pm 1}, p_1, \ldots, p_{l}]$.
Let $h: T_{\widetilde{\bfu}} \to \QQ$ be the $\QQ$-algebra homomorphism uniquely determined by 
\[
h(u(i, 0)) = f(i, 0), \;\;i \in [1, r], \hs \mbox{and} \hs  
h(p_j) = f(p_j), \;\; j \in [1, l].
\]
Denote also by $h$ the ring homomorphism 
$h: \calA_\ZZ(\BA \bb P, \emptyset) \to \QQ$ via restriction. By the recursive nature of \eqref{eq:ex-0}, one sees that \eqref{eq:F-f} holds for all $(i, m) \in [1, r] \times \ZZ$. Since $\calA_\ZZ(\BA \bb P, \emptyset)$ is generated 
over $\ZZ[\bfp]$ by 
${\mathcal{Z}}$, $h(a) \in \ZZ$ for every $a \in \calA_\ZZ(\BA \bb P,\emptyset)$, and one thus has the ring homomorphism
$h: \calA_\ZZ(\BA \bb P,\emptyset) \to \ZZ$. Applying \leref{:poly-0} to the $\calS(\BA \bb P)$-frieze testing set
\[
\calT = {\mathcal{Z}} \cup \{u(1,0), \ldots, u(r,0), p_1, \ldots, p_{l}\},
\]
one sees that $h: \calA_\ZZ(\BA \bb P, \emptyset) \to \ZZ$ is a frieze. Since both $f$ and $h$ are determined by their values on $\widetilde{\bfu}$, 
there can be only one $h$ satisfying \eqref{eq:F-f}.
\end{proof}

When $A$ is of finite type, it is well-known, see, e.g., \cite[Theorem 3.1]{Keller:06} or \coref{:uim-F-inv} in this paper, that for all $P \in M_{l, r}(\ZZ)$ every
cluster variable of $\calA_\ZZ(\BA \bb P, \emptyset)$ appears as $u(i, m)$ for some $(i, m) \in [1, r]\times \ZZ$. We thus have the following consequence of 
\prref{:ARS-same}.
 
\bco{:finite-equivalent}
When $A$ is of finite type,  frieze patterns associated to $A$ with any coefficient matrix $P$
are equivalent to friezes of $\calA_\ZZ(\BA \bb P, \emptyset)$ via evaluation.
\eco

In \coref{:BFZ-same},  we will give another class of examples to which \prref{:ARS-same} applies.

\bex{:Ptolemy}
{\rm
For $r \geq 2$, the Ptolemy cluster algebra of rank $r$ \cite{FZ:II,cuntz-et-al:coefficients, Schiffler:ptolemy} is the one of type $A_r$ whose  clusters are in bijection with the triangulations of an
$(r+3)$-gon. More specifically, the {\it shell triangulation} of the $(r+3)$-gon, i.e. the triangulation whose $r$ diagonals have a common vertex, gives rise to 
 the seed 
\[
\calS^{\rm tol} = ({\bf u} = (u_1, \ldots, u_r), \, {\bf p} = (p_1, \ldots, p_{r+3}), \, \widetilde{B} = (B \bb P^{\rm tol}))
\]
in $\FF = \QQ({\bf u}, \bfp)$ of extended rank $(r, r+3)$, where
\[
B = \begin{pmatrix} 0 & 1 & 0 &\cdots &0\\
-1 & 0 & 1 & \cdots &0\\
0& -1 & 0 &\cdots &0\\
\cdots& \cdots  & \cdots &\cdots & \cdots\\
0&\cdots&\cdots&0&1 \\
0& \cdots&\cdots&-1 & 0
\end{pmatrix} \quad \text{ and }  \quad P^{\text{tol}} = \begin{pmatrix} 
1& -1 & 0 &\cdots &0\\
0 & 1 & -1 & \cdots&0\\
0& 0 & 1 &\cdots& 0\\
\cdots& \cdots  & \cdots &\cdots & \cdots\\
0&\cdots&\cdots&1&-1 \\
0&0 &\cdots&0 & 1 \\
0&0& \cdots&0 & -1 \\
1& 0&\cdots&0 & 0 \\
-1&0 &\cdots&0 & 0 \\
\end{pmatrix} .
\]
Note that  $B = \BA$ for the standard Cartan matrix $A$ of type $A_r$. 
Frieze patterns with coefficients defined 
in \cite{cuntz-et-al:coefficients} are the same as frieze patterns associated to $(B \bb P^{\rm tol})$ per our \deref{:pattern-intro}
(see also \cite[Remark 3.4]{cuntz-et-al:coefficients}), and they  
can be arranged so as to generalize Coxeter's original definition of frieze patterns.  For $r = 3$, we have the generic frieze pattern,

\footnotesize
\[
\begin{tikzcd}[column sep=0em, row sep=0.25em]
	&{\bf p_2}	&	&{\bf p_3}&	&{\bf p_4}	&	&{\bf p_5}&	&{\bf p_6}&	&{\bf p_1}\\
{\bf u_1}	& p_1p_3	&{\bf u_4}&p_2p_4	&{\bf u_7}	&p_3p_5	&{\bf u_9}&p_4p_6	&{\bf u_3}	&p_5p_1&{\bf u_6}	&p_2p_6&{\bf u_1}	\\
	&{\bf u_2}	& p_1p_4	&{\bf u_5}	&p_2p_5&{\bf u_{8}}&	p_3p_6&{\bf u_2}&	p_4p_1&{\bf u_5}	&p_2p_5&{\bf u_8}&	p_3p_6&{\bf u_2}\\
	&	&{\bf u_3}	& p_1p_5	&{\bf u_6}&p_2p_6&{\bf u_1}&	p_3p_1&{\bf u_4}&p_2p_4	&{\bf u_7}&p_3p_5&{\bf u_9}&	p_4p_6&{\bf u_3}\\
&	&	&{\bf p_6}&	&{\bf p_1}&	&{\bf p_2}&	&{\bf p_3}&	&{\bf p_4}&	&{\bf p_5}&	&,
\end{tikzcd}
\]
\normalsize
where $\{{\bf u}_k: k \in [1, 9]\}$ are cluster variables, $\{{\bf p}_j=p_j: j \in [1, 6]\}$ are the frozen variables, and for any four adjacent boldfaced entries arranged in a diamond {\scriptsize $\begin{matrix}
			&{\bf b}&\\
			{\bf a}&\square&{\bf d}\\
			&{\bf c}&
		\end{matrix},
	$}
we have the {\it generalized diamond rule} 
${\bf a} {\bf d} - {\bf b} {\bf c}  = \square$ with $\square$ being the product of two frozen variables at the center of the diamond. 
We refer to  
\cite{cuntz-et-al:coefficients} for details.

	}
  \eex

\section{Frieze points in examples}\label{s:3-cases}
Fix the $r\times r$ symmetrizable generalized Cartan matrix
$A = (a_{i, j})$, and let
$\BA$ be as in \eqref{eq:Bo-00}. Recall from 
$\S$\ref{ss:patterns-intro} that we introduced three extended mutation matrices
\begin{equation}\label{eq:MM}
 \widetilde{B}^{\BFZ} = \left(\begin{array}{c} \BA \\ \UA\end{array}\right), \hs \widetilde{B}^{\rm prin} = \left(\begin{array}{c} \BA \\ I_r\end{array}\right), \hs
 \hs \mbox{and} \hs \widetilde{B}^{\rm triv} = \BA,
\end{equation}
where $\UA$ is given in \eqref{eq:U0-00}. 
In this section, we prove the combined statements in \thref{:3-cases} on frieze points associated to these three mutation matrices. 

\subsection{The case of BFZ coefficients}\label{ss:BFZ-6} Consider the seed 
\begin{equation}\label{eq:Sigma-BFZ}
\Sigma^{\BFZ} = ({\bf z} = (z_1, \, \ldots, \, z_r),\, {\bf p} = (p_1, \, \ldots, \, p_r), \, \widetilde{B}^{\BFZ})
\end{equation}
in $\FF^{\BFZ} = \QQ(z_1, \ldots, z_r, z_{r+1}, \ldots, z_{2r})$, where  $(z_1, \ldots, z_{2r})$ are independent variables
and 
\begin{equation}\label{eq:pi}
p_i = z_i z_{r+i} -\prod_{j=i+1}^r z_j^{-a_{j,i}} \prod_{j=1}^{i-1} z_{r+j}^{-a_{j,i}}, \hs \hs i \in [1, r].
\end{equation}
Let $\bsig^{\BFZ}$ be the seed pattern in $\FF^{\BFZ}$ with the seed $\Sigma^{\BFZ}$ at $t_0$, and let ${\bf z}[1]$ be the cluster of the seed $\mu_r \cdots \mu_2 \mu_1 \Sigma^{\BFZ}$.

\bpr{:BFZ-1} \cite[(7) and Proposition 7.4]{GLS:factorial} One has ${\bf z}[1] = 
(z_{r+1}, \ldots, z_{2r})$, and
\[
\calA_\QQ(\bsig^{\BFZ}, \emptyset) = \ZZ[z_1, \ldots, z_r, z_{r+1}, \ldots, z_{2r}].
\]
\epr

We now have the following immediate consequence of \prref{:BFZ-1} and \leref{:poly-0}.

\bpr{:BFZ-test} One has ${\rm Spec}_\CC(\calA_\CC(\bsig^{\BFZ}, \emptyset)) \cong \CC^{2r}$, and 
$(z_1, \ldots, z_{2r}) \in \CC^{2r}$ is an $\bsig^{\BFZ}$-frieze point if and only if 
$(z_1, \ldots, z_{2r}) \in (\ZZ_{>0})^{2r}$ and 
$p_i(z_1, \ldots, z_{2r}) > 0$ for every $ i \in [1, r]$.
\epr

In the notation of $\S$\ref{ss:belt},  ${\bf z}$ and ${\bf z}[1]$ are respectively the clusters of the seeds at the vertices $t(1,0)$ and $t(1, 1)$ on the acyclic belt through $\Sigma^{\BFZ}$. 
We thus have the following consequence of \prref{:ARS-same}.

\bco{:BFZ-same}
For any $r \times r$ symmetrizable generalized Cartan matrix $A$, frieze patterns associated to $A$  with coefficient matrix in $\UA$ 
are,  via evaluation, in bijection with  friezes of the polynomial ring 
$\ZZ[z_1, \ldots, z_r, z_{r+1}, \ldots, z_{2r}]$ equipped with the $\bsig^{\BFZ}$-cluster structure.
 
\eco

Solving $z_{r+1}, \ldots, z_{2r}$ from the equations in \eqref{eq:pi} and writing them as
\begin{equation}\label{eq:z:ri}
z_{r+i} = Q_i(z_1, \ldots, z_r, p_1, \ldots, p_r) \in \ZZ_{\geq 0} [z_1^{\pm 1}, \ldots, z_r^{\pm r}, p_1, \ldots, p_r], \hs i \in [1, r],
\end{equation}
we see that $\bsig^\BFZ$-frieze points in $\CC^{2r}$ are also in bijection with the set
\[
 \{(z_1, \ldots, z_r, p_1, \ldots, p_r) \in (\ZZ_{>0})^{2r}: \, Q_i(z_1, \ldots, z_r, p_1, \ldots, p_r) \in \ZZ, \; \forall \; i \in [1, r]\}.
 \]


\bex{:G2-1}
{\rm Let $A$ be of type $G_2$, so that
\[
 \widetilde{B}^{\BFZ} = \left(\begin{array}{cc} 0 & 3 \\ -1 & 0\\ 1 & -3 \\ 0 & 1\end{array}\right).
\]
One has $p_1 = z_1z_3 -z_2$ and $p_2 = z_2z_4 - z_3^3$, and the other four cluster variables are 
\begin{align*}
z_5 & = z_1z_4-z_3^2, \hs z_6 = z_1^3z_4^2-3z_1^2z_3^2z_4+3z_1z_3^4+z_2^2z_4-2z_2z_3^3,\\
z_7 & = z_1^2z_4-3z_1z_3^2+3z_2z_3+z_1z_3^2-2z_2z_3, \hs z_8 = z_1^3z_4-3z_1^2z_3^2+3z_1z_2z_3-z_2^2.
\end{align*}
Moreover, one has 
\begin{equation}\label{eq:z3z4-G2}
z_3 = \frac{p_1+z_2}{z_1}, \hs z_4 = \frac{p_2z_1^3 + (p_1 + z_2)^3}{z_1^3z_2}.
\end{equation}
Consequently, $\bsig^{\BFZ}$-frieze points in $\CC^4$ of type $G_2$ are all $(z_1, z_2, z_3, z_4) \in (\ZZ_{>0})^4$ such that
$z_1z_3 > z_2$ and $z_2z_4 > z_3^3$. Alternatively, such frieze points are parametrized by 
all $(z_1, z_2, p_1, p_2) \in (\ZZ_{>0})^4$ such that $z_3$ and $z_4$ given in \eqref{eq:z3z4-G2} are integers. 
\hfill $\diamond$
}
\eex

\subsection{The case of principal coefficients}\label{ss:prin} Consider  
the matrix $\widetilde{B}^{\rm prin}$ in \eqref{eq:MM} and the seed 
\[
\Sigma^{\rm prin} = ({\bfx =(x_1, \ldots, x_r),   \; \bfy = (y_1, \ldots, y_r}), \; \widetilde{B}^{\rm prin})
\]
in $\FF^{\rm prin} = \QQ(\bfx, \bfy)$.  Let
$\bsig^{\rm prin}$ be the seed pattern in $\FF^{\rm prin}$ with seed $\Sigma^{\rm prin}$ at $t_0 \in \TT_r$. Correspondingly we have 
the cluster algebras 
$\calA_{\ZZ}(\bsig^{\rm prin}, \emptyset)$ and $\calA_{\ZZ}(\bsig^{\rm prin}) = \calA_{\ZZ}(\bsig^{\rm prin}, \emptyset)[\bfy^{\pm 1}]$.

The mutation of $\Sigma^{\rm prin}$ at each $i \in [1, r]$ gives the {\it lower bound cluster variable} $x_i^\prime$ via
\begin{equation}\label{eq:ppi}
x_ix_i^\prime = y_i q_i^+ + q_i^-, \hs \mbox{where}\;\;\; q_i^+ = \prod_{j=1}^{i-1} x_j^{-a_{j, i}} \;\;\; \mbox{and} \;\;\; q_i^- =\prod_{j=i+1}^r x_j^{-a_{j, i}}.
\end{equation}
Following  \cite{BFZ:III}, for ${\bfm} = (m_1, \ldots, m_r)^T\in \ZZ^r$, define $\bfx^{\lan {\bfm}\ran} = x_1^{\lan m_1 \ran} \cdots x_r^{\lan m_r\ran}$, 
\[
 x_j^{\lan m_j \ran} = \begin{cases} x_j^{m_j}, & \hs \mbox{if} \; m_j \geq 0, \\ (x_j^\prime)^{-m_j},  & \hs \mbox{if} \; m_j \leq 0,\end{cases} \hs j \in [1, r].
\]
Each $\bfx^{\lan {\bfm}\ran}$, called {\it a standard monomial} \cite{BFZ:III}, is an element in $\calA_{\ZZ}(\bsig^{\rm prin}, \emptyset)$. 
 By \cite[Corollary 1.21]{BFZ:III}, the set
${\mathcal{M}} = \{\bfx^{\lan {\bfm}\ran}: \; {\bfm} \in \ZZ^r\}$ is a basis of $\calA_{\ZZ}(\bsig^{\rm prin})$ as a module over the Laurent polynomial ring $\ZZ[\bfy^{\pm 1}]$. 
The following  \prref{:calM} is proved in the Appendix.

\bpr{:calM} The set ${\mathcal{M}}$ of all standard monomials is a basis of $\calA_{\ZZ}(\bsig^{\rm prin}, \emptyset)$ 
as a module over the polynomial ring $\ZZ[\bfy]$.
\epr

\bre{:GLS-wrong}
{\rm 
Contrary to what is stated in \cite[Proposition 7.5]{GLS:factorial}, the analog of \prref{:calM} may not hold for non-principal coefficients.
Consider, for example,  the case of $A_2$ with initial seed $((z_1, z_2, p_1, p_2), \widetilde{B}^{\BFZ})$, where
\[
\widetilde{B}^{\BFZ} = \left(\begin{array}{cc} 0 & 1 \\ -1 & 0\\ 1 & -1 \\ 0 & 1\end{array}\right).
\]
Let z$_3 = z_1^\prime$ and $z_5 = z_2^\prime$ be the two new cluster variables obtained by mutating the initial cluster in directions $1$ and $2$, i.e., 
$z_1z_1^\prime = p_1 + z_2$ and $z_2z_2^\prime  = p_1 + z_1 p_2$.
Let $z_4$ be the only cluster variable different from $z_1, z_2, z_3, z_5$.  
One checks directly that 
\[
z_4 = \frac{z_1'z_2' - p_2 }{p_1}, 
\]
which is not in the $\ZZ[p_1, p_2]$-span of standard cluster monomials formed by $\{z_1, z_1^\prime, z_2, z_2^\prime\}$. 
\hfill $\diamond$
}
\ere

Consider now the polynomial ring $\ZZ[x_1, x_1^\prime, \ldots, x_r, x_r^\prime, y_1, \ldots, y_r]$ and introduce
\[
\phi_i = y_i \prod_{j=1}^{i-1} x_j^{-a_{j, i}} + \prod_{j=i+1}^r x_j^{-a_{j, i}},   \hs i \in [1, r].
\]
Let $\calI$ be the ideal of $\ZZ[x_1, x_1^\prime, \ldots, x_r, x_r^\prime, y_1, \ldots, y_r]$ generated by the elements $x_ix_i^\prime - \phi_i$, $i \in [1, r]$. Let
$V^{\rm prin}$ be the affine sub-variety of $\CC^{3r}$ defined by the equations
\[
x_ix_i^\prime  -\phi_i = 0, \hs i \in [1, r].
\]

\bpr{:prin-1}
One has  $\calA_\QQ(\bsig^{\rm prin}, \emptyset) \cong \ZZ[x_1, x_1^\prime, \ldots, x_r, x_r^\prime, y_1, \ldots, y_r]/\calI$.
Consequently, 
\begin{equation}\label{eq:iso}
{\rm Spec}(\calA_\CC(\bsig^{\rm prin}, \emptyset)) \cong V^{\rm prin} \subset \CC^{3r}.
\end{equation}
The set of $\bsig^{\rm prin}$-frieze points in $V^{\rm prin}$ is precisely  $(\ZZ_{>0})^{3r} \cap V^{\rm prin}$.
\epr

\begin{proof}
By \prref{:calM},  one has a surjective ring homomorphism
\[
H: \;\; \ZZ[x_1, x_1^\prime, \ldots, x_r, x_r^\prime, y_1, \ldots, y_r]/\calI \longrightarrow \calA_\QQ(\bsig^{\rm prin}, \emptyset).
\]
By \cite[Corollary 1.21]{BFZ:III}, the polynomials $x_ix_i^\prime-\phi_i$, for $i \in [1, r]$, generate the ideal of relations
among the variables $x_1, \ldots, x_r, x_1^\prime, \ldots, x_r^\prime$ in $\calA_\QQ(\bsig^{\rm prin}, \emptyset)$. Thus $H$ is injective.

By \leref{:poly-0}, $\{x_1, x_1^\prime, \ldots, x_r, x_r^\prime, y_1, \ldots, y_r\}$ is a testing set for $\bsig^{\rm prin}$-friezes. Thus 
the set of $\bsig^{\rm prin}$-frieze points in $V^{\rm prin}$ is precisely  $(\ZZ_{>0})^{3r} \cap V^{\rm prin}$.
\end{proof}


\subsection{The case of trivial coefficients}\label{ss:trivial}
Turning to the seed pattern $\calS^{\rm triv}$ with trivial coefficients defined by $\BA$.  
We have then two coefficient specializations \cite[Definition 12.1]{FZ:IV}
\begin{align}
\label{eq:psi-BFZ}
&\psi^{\BFZ}: \; \calA_\CC(\bsig^{\BFZ}, \emptyset) \longrightarrow \calA_\CC(\bsig^{\rm triv}),\\
\label{eq:psi-prin}
&\psi^{\rm prin}: \; \calA_\CC(\bsig^{\rm prin}, \emptyset) \longrightarrow \calA_\CC(\bsig^{\rm triv})
\end{align}
such that $\psi^{\BFZ}(p_i) = \psi^{\rm prin}(y_i) = 1$ for every $i \in [1, r]$. Consequently,
\[
\calA_\CC(\bsig^{\rm triv}) \cong \calA_\CC(\bsig^{\BFZ}, \emptyset)/\calJ \cong \calA_\CC(\bsig^{\rm prin}, \emptyset)/\calK,
\]
where $\calJ$ is the ideal of  $\calA_\CC(\bsig^{\BFZ}, \emptyset)$ generated by $\{p_i-1: i \in [1, r]\}$, and $\calK$ is the ideal of 
$\calA_\CC(\bsig^{\rm prin}, \emptyset)$ generated by $\{y_i-1: i \in [1, r]\}$.
By \prref{:BFZ-test}  and  \prref{:prin-1}, 
one has
\[
{\rm Spec} (\calA_\CC(\bsig^{\rm triv})) \cong W^{\rm triv} \cong V^{\rm triv},
\]
where $W^{\rm triv}\subset \CC^{2r}$ is 
defined in the coordinates $(x_1, x_1^\prime, \ldots, x_r, x_r^\prime)$ of $\CC^{2r}$ by the equations
\begin{equation}\label{eq:xxi-1}
x_ix_i^\prime = \prod_{j=1}^{i-1} x_j^{-a_{j, i}} + \prod_{j=i+1}^r x_j^{-a_{j, i}},   \hs i \in [1, r],
\end{equation}
and $V^{\rm triv} \subset \CC^{2r}$ 
 is defined in the coordinates $(z_1, \ldots, z_{2r})$ of $\CC^{2r}$ by the equations
\begin{equation}\label{eq:zzi-1}
z_i z_{r+i}  = 1 + \prod_{j=i+1}^r z_j^{-a_{j,i}} \prod_{j=1}^{i-1} z_{r+j}^{-a_{j,i}},  \hs i \in [1, r].
\end{equation}
It again follows from \prref{:BFZ-test} and \prref{:prin-1} that 
the set of $\bsig^{\rm triv}$-frieze points in $W^{\rm triv}$ is $ W^{\rm triv} \cap (\ZZ_{>0})^{2r}$, and the set of $\bsig^{\rm triv}$-frieze 
points  in  $V^{\rm triv}$ is $V^{\rm triv} \cap (\ZZ_{>0})^{2r}$.

\section{Frieze points in the reduced double Bruhat cell $\Lcc$}\label{s:Lcc}
In this section, we  assume that $A$ is a Cartan matrix of finite type, and we describe in Lie theoretical terms the frieze points in  the 
three cases discussed in $\S$\ref{s:3-cases} as well as the periodicity of frieze patterns associated to $A$ with any coefficients.

\subsection{Notation and preliminaries on generalized minors}\label{ss:bkgrnd-minor}
Throughout this section, let $G$ be a connected and simply connected complex semi-simple Lie group, 
 let $(B,B_-)$ be a pair of opposite Borel subgroups of $G$,
 and let $T = B\cap B_-$, a maximal torus of $G$. Let $\t$ be the Lie algebra of $T$, and let $\Phi \subset \t^*$ and $\Phi_+ \subset \Phi$  be respectively the
 sets of roots and positive roots determined by $(B, T)$, and let $\Phi^\vee \subset \t$ be the set of co-roots. We label the simple roots and simple co-roots respectively 
 as $\{\alpha_1, \alpha_2, \cdots , \alpha_r\} \subset \t^*$ and $\{\alpha_1^\vee, \alpha_2^\vee, \cdots , \alpha_r^\vee\} \subset \t$. One then has the Cartan matrix $
A = (a_{i, j})_{i, j \in [1, r]}$, where  
\begin{equation}
	a_{i,j} = \alpha_i^\vee(\alpha_j), \hs i, j \in [1, r].
\end{equation}

Let $W$ be the Weyl group of $(G,T)$. For each $i \in [1,r]$, let $s_i \in W$ be the simple reflection defined by $\alpha_i$. 
For $i \in [1,r]$, we fix root vectors $e_i$ for $\alpha_i$ and $e_{-i}$ for $-\alpha_i$ such that $[e_i,e_{-i}] = \alpha_i^\vee\in \t$, and let 
$\theta_i: \text{SL}_2(\CC) \to G$ be the Lie group homomorphism determined by $\{e_i,e_{-i},\alpha_i^\vee\}$. By abuse of notation,  we also let
 $e_i: \CC \to N$, $e_{-i}: \CC \to N_-$ and $\alpha_i^\vee: \CC^* \to T$ denote the one-parameter subgroups respectively given by 
\begin{equation}
	e_i(z) = \theta_i\begin{pmatrix}1&z\\0&1 \end{pmatrix}, \qquad e_{-i}(z) = \theta_i\begin{pmatrix}1&0\\z&1 \end{pmatrix}, 
  \qquad \alpha_i^\vee(z) = \theta_i\begin{pmatrix}z&0\\0&z^{-1} \end{pmatrix}.
\end{equation}
Here, $N$ and $N_-$ denote respectively the unipotent radicals of $B$ and $B_-$. For $i \in [1, r]$, let
\begin{equation}\label{eq:si}
	\overline{s}_i = \theta_i\begin{pmatrix}0&-1\\1&0 \end{pmatrix} = e_i(-1)e_{-i}(1)e_i(-1)\in N_G(T),
\end{equation}
where $N_G(T)$ is the normalizer of $T$ in $G$.
Each $w \in W$ then \cite[$\S$1.4]{FZ:double} has the
well-defined  representative $\overline{w} \in N_G(T)$, uniquely determined by the condition $\overline{w_1w_2} = \overline{w}_1\overline{w}_2$ 
whenever $\ell(w_1w_2) = \ell(w_1)+\ell(w_2)$, where $\ell: W \to {\mathbb{Z}}_{\geq 0}$ is the length function on $W$.


For $t \in T$ and $w \in W$, set $t^w =\overline{w}^{\,-1} t \overline{w} \in T$, 
so  that
\[
T \times W \longrightarrow T, \;\; (t, w) \longmapsto  t^w,
\]
is a right action of $W$ on $T$. For  $i \in [1,r]$, let $\omega_i \in\t^*$ to be the fundamental weight corresponding to the simple root $\alpha_i$. 
We identify $\calP = \sum_{i=1}^r\ZZ \omega_i$ with the lattice of algebraic characters of $T$ such that each $\gamma \in \calP$ gives rise to the morphism
$T \to \CC^*, t \mapsto t^\gamma$, determined by
\begin{equation*}
	\alpha_i^\vee(z)^\gamma = z^{\alpha_i^\vee(\gamma)}, \quad  i \in [1,r],\; z \in \CC^*. 
\end{equation*}
One then has
$(t^w)^\gamma = t^{w \gamma}$ for $t \in T, w \in W$.

We now recall some basic facts on generalized minors from \cite{FZ:double}. For  $i \in [1, r]$ and  $u, v \in W$, the {\it generalized minor} $\Delta_{u\omega_i, v\omega_i}$ is the regular function on $G$ uniquely determined by 
\[
\Delta_{u\omega_i, v\omega_i}(g)= \Delta_{\omega_i, \omega_i}(\overline{u}^{\, -1} g \overline{v}) = 
([\overline{u}^{\, -1} g \overline{v}]_0)^{\omega_i}, \quad g \in \overline{u} B_- B \overline{v}^{\, -1},
\]
where  for $h \in B_-B$ we write 
$h = [h]_- [h]_0 [h]_+$ with $[h]_- \in N_-, [h]_0 \in T$ and $[h]_+ \in N$.
Moreover, $\Delta_{u\omega_i, v\omega_i}$ depends only on the weights $u\omega_i$ and $v\omega_i$ and not on the choices of $u, v \in W$. 
Note from the definitions that for any $i \in [1, r]$ and $\delta, \gamma \in W\omega_i$, one has
\begin{equation}\label{eq:Delta-tt}
\Delta_{\gamma, \delta}(t_1gt_2) = t_1^{\gamma}t_2^{\delta} \Delta_{\omega_i, \, \omega_i}(g), \hs t_1, t_2 \in T, g \in G.
\end{equation}
We will use the following identity in $\CC[G]$ from \cite[Theorem 1.17]{FZ:double}.

\ble{le:master}
If $u, v \in W$ and $i \in [1, r]$ are such that
$l(us_i) = l(u) +1$ and $l(vs_i) = l(v)+1$, then
\begin{equation}\label{eq:master-equation}
\Delta_{u\omega_i, \, v\omega_i} \Delta_{us_i\omega_i, \, vs_i\omega_i} =
\Delta_{us_i\omega_i, \, v\omega_i} \Delta_{u\omega_i, \, vs_i\omega_i} +\prod_{j \neq i} \Delta_{u\omega_j, \, v\omega_j}^{-a_{j, i}}.
\end{equation} 
\ele


 For $i \in [1, r]$ and $z \in \CC$, let
\begin{equation}\label{eq:ga}
q_i(z) = e_i(z) \overline{s_i} =\theta_{i}\left(\begin{array}{cc} z & -1 \\ 1 & 0 \end{array}\right) \in G.
\end{equation}
The following facts are well-known (see \cite[Proposition 2.2]{FZ:double} and \cite[Remark 7.5]{BZ:tensor}).

\ble{:gig} Let $i \in [1, r]$ and $\delta \in W\omega_i$. Then for 
any $z \in \CC$ and $g \in G$,
\begin{align}\label{eq:gq2}
\Delta_{\omega_i, \, \delta}(q_{i}(z)g) &= z \Delta_{\omega_i, \, \delta}(g) - \Delta_{s_i\omega_i, \, \delta}(g),\hs
\Delta_{\omega_i, \, \delta}(e_{i}(z)g)  =z \Delta_{s_i\omega_i, \, \delta}(g) +\Delta_{\omega_i, \, \delta}(g),\\
\label{eq:gq1}
\Delta_{\delta, \,\omega_i}(g q_{i}(z)) &= z \Delta_{\delta, \,\omega_i, \, }(g) + \Delta_{\delta, \,s_i\omega_i}(g),\hs
\Delta_{\delta, \,\omega_i}(g e_{-i}(z)) = z \Delta_{\delta, \,s_i\omega_i}(g) + \Delta_{\delta, \,\omega_i}(g).
\end{align}
Furthermore, for any $i' \in [1, r]$ and $i' \neq i$, and for any $a \in SL(2, \CC)$ and $g \in G$,
\[
\Delta_{\omega_i, \, \delta}(\theta_{i'}(a)g) =\Delta_{\omega_i, \, \delta}(g), \hs   
\Delta_{\delta, \,\omega_i}(g \theta_{i'}(a)) = \Delta_{\delta, \, \omega_i}(g).
\]
\ele

For a sequence ${\bf i} = (i_1, \ldots, i_l)$ in the alphabet $\{-r, \ldots, -1, 1, \ldots, r\}$, define 
\[
e_{\bf i}: \;\; \CC^l \longrightarrow G,\;\; e_{\bf i}(z_1, z_2, \ldots, z_l)  = e_{i_1}(z_1) e_{i_2}(z_2) \cdots e_{i_l}(z_l).
\]
The following  \leref{:minor-ui} is a generalization of \cite[Theorem 5.8]{BZ:tensor} and is proved the same way. We give an outline of the proof 
of \leref{:minor-ui} for completeness.

\ble{:minor-ui}
Let ${\bf i} = (i_1, \ldots, i_l)$ be any sequence in the alphabet $\{-r, \ldots, -1, 1, \ldots, r\}$. Then for any $i \in [1, r]$ and 
$\gamma, \delta \in W\omega_i$, $\Delta_{\gamma, \delta}(e_{\bf i}(z_1, \ldots,z_l))$ is a polynomial in $(z_1, \ldots, z_l)$ 
with non-negative integral coefficients. 
\ele

\begin{proof}
Note first that if  $f \in \CC[G]$ and $\gamma_1, \gamma_2 \in \calP$ are such that  $f(e) \neq 0$ and 
\begin{equation}\label{eq:fggt}
f(t_1gt_2) = t_1^{\gamma_1} t_2^{\gamma_2} f(g), \hs \forall \;\; t_1, t_1 \in T, \; g \in G,
\end{equation}
 then  $f(e) = t^{\gamma_1-\gamma_2} f(e)$ for any $t \in T$ and thus $\gamma_1 =\gamma_2$.

Consider now $\CC[G]$ with the left $G$-action given by $(g_1 \cdot f)(g) = f(gg_1)$ for $g, g_1 \in G$, and denote the induced
$U(\g)$-module structure on $\CC[G]$ by $(a, f) \mapsto a \cdot f$ for $a \in U(\g)$. Then 
\[
\Delta_{\gamma, \delta}(e_{\bf i}(z_1, \ldots, z_l))  = ((e_{i_1}(z_1), \ldots, e_{i_l}(z_l)) \cdot \Delta_{\gamma, \delta})(e).
\]
The sub-module $U(\g) \Delta_{\gamma, \omega_i} \subset \CC[G]$ has highest weight $\omega_i$ and highest weight vector
$\Delta_{\gamma, \omega_i}$, and $\Delta_{\gamma, \delta} \in U(\g) \Delta_{\gamma, \omega_i}$ is an extremal weight vector
with extremal weight $\delta$.
For $k \in [1, r]$ and $m \in \ZZ_{\geq 0}$, let
$e_{\pm k}^{(m)} \in U(\g)$ be the divided power  given by $e_{\pm k}^{(m)}:= e_{\pm k}^m/m!$. Then
\[
\Delta_{\gamma, \delta}(e_{\bf i}(z_1, \ldots, z_l))=\sum_{k_1=0}^\infty\cdots \sum_{k_l=0}^\infty z_1^{k_1} \cdots z_l^{k_l} 
\left(e_{i_1}^{(k_1)} \cdots e_{i_l}^{(k_l)} \cdot \Delta_{\gamma, \delta}\right)(e),
\]
where the  sum has only finitely many non-zero terms. For $j \in [1,l]$, let ${\rm sign}(i_j) = -1$ if $i_j \in \{-r, \ldots, -1\}$ and 
${\rm sign}(i_j) = 1$ if $i_j \in [1, r]$. For 
 $(k_1, \ldots, k_l) \in (\ZZ_{\geq 0})^l$, let 
\[
b_{k_1, \ldots, k_l} = \left(e_{i_1}^{(k_1)} \cdots e_{i_l}^{(k_l)} \cdot \Delta_{\gamma, \delta}\right)(e) \in \CC \;\;\;\;
\mbox{and} \;\;\;\; \delta_{k_1, \ldots, k_l} = \delta +\sum_{j=1}^l k_j ({\rm sign}(i_j) \alpha_{|i_j|}) \in \calP.
\]
Suppose that $b_{k_1, \ldots, k_l} \neq 0$ and let $f_{k_1, \ldots, k_l} = e_{i_1}^{(k_1)} \cdots e_{i_l}^{(k_l)} \cdot \Delta_{\gamma, \delta} \in \CC[G]$.
Then $f_{k_1, \ldots, k_l} (e) \neq 0$, and $f_{k_1, \ldots, k_l} $ satisfies
\eqref{eq:fggt} with $\gamma_1 = \gamma$ an $\gamma_2 = \delta_{k_1, \ldots, k_l}$.  Thus 
$\delta_{k_1, \ldots, k_l} = \gamma$.  Consequently, both
$f_{k_1, \ldots, k_l} $ and $\Delta_{\gamma, \gamma}$ span the
$1$-dimensional weight sub-space of
$U(\g) \Delta_{\gamma, \omega_i}$ with extremal weight $\gamma$. Since $\Delta_{\gamma, \gamma}(e) = 1$, one thus has
\begin{equation}\label{eq:ee-b}
e_{i_1}^{(k_1)} \cdots e_{i_l}^{(k_l)} \cdot \Delta_{\gamma, \delta} = b_{k_1, \ldots, k_l} \Delta_{\gamma, \gamma} \in \CC[G].
\end{equation}
Let $\gamma = u\omega_i$ for $u \in W$. Evaluating both sides of \eqref{eq:ee-b} at $\overline{u} g$ for all $g \in G$, one gets
\[
e_{i_1}^{(k_1)} \cdots e_{i_l}^{(k_l)} \cdot \Delta_{\omega_i, \delta} = b_{k_1, \ldots, k_l} \Delta_{\omega_i, \gamma} \in \CC[G].
\]
By the proof of \cite[Theorem 5.8]{BZ:tensor}, there exist a monomial $a$ in the divided powers
of  $e_1, \ldots, e_r$  and a monomial $a'$ in the divided powers of $e_{-1}, \ldots, e_{-r}$ such that
$a \cdot \Delta_{\omega_i, \gamma} = \Delta_{\omega_i, \omega_i}$ and $a' \cdot \Delta_{\omega_i, \omega_i} = \Delta_{\omega_i, \delta}$.  
It follows that 
\[
\left(a e_{i_1}^{(k_1)} \cdots e_{i_l}^{(k_l)} a'\right)\cdot \Delta_{\omega_i, \, \omega_i} = b_{k_1, \ldots, k_l} \Delta_{\omega_i, \, \omega_i}.
\]
By  \cite[Lemma 7.4]{BZ:tensor}, 
$b_{k_1, \ldots, k_l}$ is a positive integer.
\end{proof}

For $v \in W$, let
$N_v = N \cap (\ov N_- \ov^{\, -1})$.
The following facts will be used in $\S$\ref{ss:Lie-prin}. 
\ble{le:Nww} 1) If $v, v_1, v_2 \in W$ are such that $v = v_1v_2$ and $l(v) = l(v_1) + l(v_2)$, then one has the unique decomposition
$N_v \ov= (N_{v_1} \overline{v_1}) (N_{v_2} \overline{v_2})$.

2) Suppose that $v_1, \ldots, v_p \in W$ are such that $l(v_1\cdots v_p) = l(v_1) + \cdots + l(v_p)$, and let $1 \leq p_1 \leq p_2 \leq p$.
Then for any $i \in [1, r]$ and $q_j \in N_{v_j}\overline{v_j}$ for $j \in [1, p]$, one has
\[
\Delta_{v_1 \cdots v_{p_1}\omega_i, \, v_1 \cdots v_{p_2}\omega_i}(q_{1} \cdots q_{p} \overline{v_1\cdots v_p}^{\, -1})
=\Delta_{\omega_i, \, \omega_i}(q_{{p_1+1}} \cdots q_{{p_2}}).
\]
\ele

\begin{proof}
1) follows from \cite[Proposition 2.11]{FZ:double}. For 2), since $q_{p_2+1}\cdots q_p \overline{v_{p_2+1} \cdots v_p}^{\, -1} \in N$ and
$\overline{v_1\cdots v_{p_1}}^{\, -1} q_1 \cdots q_{p_1} \in N_-$, 
 one has
\begin{align*}
\Delta_{v_1 \cdots v_{p_1}\omega_i, \, v_1 \cdots v_{p_2}\omega_i}(q_{1} \cdots q_{p} \overline{v_1\cdots v_p}^{\, -1}) & = 
\Delta_{v_1 \cdots v_{p_1}\omega_i, \, \omega_i}(q_{1} \cdots q_{p_2} q_{p_2+1}\cdots q_p \overline{v_{p_2+1} \cdots v_p}^{\, -1}) \\
&=
\Delta_{\omega_i,\, \omega_i}(\overline{v_1\cdots v_{p_1}}^{\, -1} q_1 \cdots q_{p_1} q_{p_1+1} \cdots q_{p_2}) \\
&=\Delta_{\omega_i, \, \omega_i}(q_{{p_1+1}} \cdots q_{{p_2}}).
\end{align*}
\end{proof}

\subsection{The cluster structure on $\Lcc$ with BFZ coefficients}\label{ss:BFZ}
Recall that for any $u, v \in W$, the double reduced cell $G^{u, v}$ is defined as
\[
G^{u, v} = BuB \cap B_- v B_- \subset G.
\]
By \cite[$\S$2]{BFZ:III}, each double reduced word of $(u, v)$ gives rise to a seed in 
$\CC(G^{u, v})$, and  it is proved in \cite[Page 2 and Theorem 1.2]{ShenWeng:BS} that all such seeds are mutation equivalent and that they
define a cluster structure on $G^{u, v}$ which we will refer to as the {\it BFZ cluster structure on $G^{u, v}$}. Recall that as a subset of $G^{u, v}$, 
the reduced double Bruhat cell 
\[
L^{u, v} = N \overline{u} N \cap B_- v B_-
\]
is  defined by setting $\Delta_{u \omega_i, \omega_i} = 1$ for all $i \in [1, r]$. Since $\{\Delta_{u\omega_i, \omega_i}: i \in [1, r]\}$ is a subset of the set of all frozen variables
for the BFZ cluster structure on $G^{u, v}$, 
one obtains by coefficient specialization \cite[Proposition 4.10]{YZ:Lcc} a cluster structure on $L^{u, v}$, which we also call the {\it BFZ cluster structure on $L^{u, v}$}.

Consider now the Coxeter element 
$c=s_1s_2 \cdots s_r \in W$ and the BFZ cluster structure, denoted by $\calS^{\BFZ}(\Lcc)$, on the reduced double Bruhat cell 
\[
\Lcc = N \overline{c} N \cap B_- c^{-1}B_-.
\]
Let
$\Sigma^{\BFZ}(\Lcc)$ be the seed in $\CC(\Lcc)$ defined by the double reduced word 
\[
(r, \,\ldots, \,2,\, 1,\, -1,\, -2, \,\ldots,\, -r)
\]
of $(c, c^{-1})$ via \cite[Theorem 2.10]{BFZ:III}. For $i \in [1, r]$,  set
\[
x_{\omega_i} = \Delta_{\omega_i, \omega_i}|_{\Lcc}, \hs x_{c\omega_i} = \Delta_{c\omega_i, c\omega_i}|_{\Lcc}, \hs \mbox{and} \hs p_{i; c} = \Delta_{\omega_i, c\omega_i}|_{\Lcc}. 
\]
One checks directly that the extended cluster of $\Sigma^{\BFZ}(\Lcc)$ 
is given by 
\[
(x_{\omega_1}, \, \ldots, \, x_{\omega_r}, \, p_{1; c}, \, \ldots, p_{r; c})
\]
with $(p_{1; c},  \ldots, p_{r; c})$ as frozen variables, 
and the extended mutation matrix of $\Sigma^{\BFZ}(\Lcc)$  is $\widetilde{B}^{\BFZ}$ in \eqref{eq:MM}  (see  also 
\cite[Example 2.24]{BFZ:III}).

Recall from \eqref{eq:pi} that for $i \in [1, r]$, we have defined 
\[
p_i = z_i z_{r+i} -\left(\prod_{j = i+1}^r z_j^{-a_{j, i}}\right) \left(\prod_{j=1}^{i-1} z_j^{-a_{j, i}}\right) \in \CC[z_1, \ldots, z_{2r}].
\]
Recall from $\S$\ref{ss:Lie-intro} that we have set 
$\CC^{2r}_{{\bf p} \neq 0} = \{(z_1, \ldots, z_{2r}) \in \CC^{2r}:  \prod_{i=1}^r p_i(z_1, \ldots, z_{2r}) \neq 0\}$.
By \prref{:BFZ-1}, the seed $\Sigma^{\BFZ}$ in 
\eqref{eq:Sigma-BFZ} defines the cluster structure $\calS^{\BFZ}$ on $\CC^{2r}_{{\bf p} \neq 0}$.
Set
\[
\rho = (x_{\omega_1}, \, \ldots, \,x_{\omega_r}, \, x_{c\omega_1}, \, \ldots, \, x_{c\omega_r}): \;\;  \Lcc \longrightarrow  \CC^{2r}.
\]

\bth{:Lcc-Af}
The map  $\rho$
induces an isomorphism from $\Lcc$ to $\CC^{2r}_{{\bf p} \neq 0}$, which sends the cluster structure $\calS^{\BFZ}(\Lcc)$ on $\Lcc$ to the cluster structure 
$\calS^{\BFZ}$ on $\CC^{2r}_{{\bf p} \neq 0}$. 
\eth

\begin{proof} Let $i \in [1, r]$.  Using the fact that $s_j \omega_i =\omega_i$ for all $j \in [1, r]$
and $j \neq i$, one has $c\omega_i = s_1\cdots s_i\omega_i$. Applying \eqref{eq:master-equation} to $u = v = s_1 \cdots s_{i-1}$, one has
\begin{align*}
\Delta_{\omega_i, \,\omega_i} \Delta_{c\omega_i,\, c\omega_i} & = \Delta_{\omega_i,\,  c\omega_i} \Delta_{c\omega_i, \, \omega_i} + \prod_{j \neq i} \Delta_{s_1\cdots s_{i-1}\omega_j, \, s_1\cdots s_{i-1}\omega_j}^{-a_{j, i}}\\
& = \Delta_{\omega_i, \, c\omega_i} \Delta_{c\omega_i,\,  \omega_i} + \prod_{j =1}^{i-1} \Delta_{c\omega_j,\,  c\omega_j}^{-a_{j, i}} \prod_{j=i+1}^r 
\Delta_{\omega_j,\,  \omega_j}^{-a_{j, i}}.
\end{align*}
Using the fact that $\Delta_{c\omega_i, \omega_i}|_{\Lcc} = 1$,  one has 
\[
p_{i; c} = \Delta_{\omega_i, c\omega_i}|_{\Lcc} = x_{\omega_i}x_{c\omega_i} -\prod_{j=i+1}^r x_{\omega_j}^{-a_{j, i}} \prod_{j = 1}^{i-1} x_{c\omega_j}^{-a_{j, i}} = \rho^* p_i.
\]
As the seed $\Sigma^{\BFZ}(\Lcc)$ in $\CC(\Lcc)$ and the seed $\Sigma^{\BFZ}$ in $\CC(z_1, \ldots, z_{2r})$ given in \eqref{eq:Sigma-BFZ} 
have the same  extended mutation matrix, 
it remains to show that the map $\rho$ induces an isomorphism of varieties from $\Lcc$ to $\CC^{2r}_{{\bf p} \neq 0}$.

Recall that $N_c = N \cap (\oc N_- \oc^{\, -1})$. Let $N_c^\prime = N \cap (\oc N\oc^{\, -1})$. One then has the direct product decomposition
$N = N_c N_c^\prime$. Thus every element $g \in N \oc N$ can be uniquely written as 
\begin{equation}\label{eq:g}
g = n_1 \oc \,n_2 n^\prime  \hs \mbox{where} \;\; n_1, n_2 \in N_c, \, n' \in N_c^\prime, 
\end{equation}
and we have the morphism
\[
\pi: \; N \oc N \longrightarrow N_c \times N_c: \; \; g \longmapsto (n_1, \, n_2).
\]
Set $(N_c \times N_c)_0 = \{(n_1, \, n_2) \in N_c \times N_c: n_1 \oc \,n_2\, \oc \in B_-B\} \subset N_c \times N_c$. 

Let $g \in \Lcc \subset N \oc N$ and write $g = n_1 \oc \,n_2 n^\prime$ as in \eqref{eq:g}. Using the unique decomposition 
\[
B_- c^{-1} B_- = B_- \oc^{\, -1} (N_- \cap \oc N \oc^{\,-1}),
\]
 we have
$g= b_- \oc^{\, -1} n_-$ for some unique $b_- \in B_-$ and $n_- \in N_- \cap \oc N \oc^{\,-1}$. Then
\[
n_1 \oc \,n_2 \oc = g (n')^{-1} \oc = b_- \oc^{\, -1} n_- (n')^{-1} \oc = b_- (\oc^{\, -1} n_- \oc) \, (\oc^{\, -1} (n')^{-1} \oc) \in B_-B.
\]
Thus
$\pi(\Lcc) \subset (N_c \times N_c)_0$.  
Conversely, given $(n_1, \, n_2)\in (N_c \times N_c)_0$, since $n_1\oc\, n_2 \oc \in B_-B$, we can write
$n_1 \oc\, n_2 \oc = b_- n_1^\prime n_2^{\prime}$, 
where $b_- \in B_-$, $n_1^\prime
\in N \cap \oc^{\, -1}  N_- \oc$ and $n_2^\prime \in N \cap \oc^{\, -1}  N\oc$.  Let
\begin{equation}\label{eq:pi-prime}
\pi^\prime(n_1,n_2) = n_1 \oc\, n_2 \oc\, (n_2^\prime)^{-1} \oc^{\, -1} = b_- \oc^{\, -1} (\oc \,n_1^\prime \oc^{-1}).
\end{equation}
Since both $\oc\, (n_2^\prime)^{-1} \oc^{\, -1}$ and $\oc \,n_1^\prime \oc^{-1}$ are in $N$, one has
$\pi^\prime(n_1, n_2)\in \Lcc$.  It is then easy to see that $\pi^\prime: (N_c \times N_c)_0 \to \Lcc$ is the inverse of 
$\pi : \Lcc \longrightarrow (N_c \times N_c)_0$, showing that the latter is 
an isomorphism.

Since an element $g \in G$ lies in $B_-B$ if and only if $\Delta_{\omega_i, \omega_i}(g) \neq 0$ for every $i \in [1, r]$, we have
\[
(N_c \times N_c)_0= \{(n_1, n_2) \in N_c \times N_c: \Delta_{\omega_i, \omega_i}(n_1 \oc \,n_2 \oc)\neq 0, \forall \, i \in [1, r]\}.
\]
As $c = s_1\cdots s_r$ is reduced, we have the parametrization \cite[Proposition 2.11]{FZ:double}
\[
\CC^{2r} \longrightarrow N_c \oc \times N_c\oc, \;\; \nu(z) =(q_1(z_1) q_2(z_2) \cdots q_r(z_r), \,
q_1(z_{r+1}) q_2(z_{r+2}) \cdots q_r(z_{2r})).
\]
For $i \in [1, r]$, define $p_i^\prime \in \CC[z_1, \ldots, z_{2r}]$ by 
\begin{equation}\label{eq:fi-x}
p_i^\prime(z) = \Delta_{\omega_i, \omega_i}(q_1(z_1) q_2(z_2) \cdots q_r(z_r)q_1(z_{r+1}) q_2(z_{r+2}) \cdots q_r(z_{2r})).
\end{equation}
It now remains to show that $p_i^\prime = p_i(z_1, \ldots, z_{2r})$ for every $i \in [1, r]$. To this end, let $z = (z_1, \ldots, z_{2r}) \in \CC^{2r}$, let
$i \in [1, r]$, and let
\[
g_i = q_{i+1}(z_{i+1})  \cdots q_{r}(z_r) q_{r+1}(z_{r+1}) \cdots q_{i-1}(z_{r+i-1}).
\]
By \leref{:gig} and using the facts that $\Delta_{\omega_i, \omega_i}(q_i(a)) = a$ for $a \in \CC$ and
$\Delta_{\omega_i, \omega_i}(\overline{s_i}) =0$, 
one has
\begin{align*}
p_i^\prime & = \Delta_{\omega_i, \, \omega_i}(q_i(z_i) g_i q_i(z_{r+i})) = z_i \Delta_{\omega_i, \, \omega_i}(g_i q_i(z_{r+i}))-\Delta_{s_i\omega_i, \, \omega_i}(g_i q_i(z_{r+i}))\\
& = z_iz_{r+i} - \Delta_{s_i\omega_i, \, \omega_i}(g_i) - \Delta_{s_i \omega_i, \, s_i\omega_i}(g_i)= z_iz_{r+i} -\Delta_{s_i \omega_i, \, s_i\omega_i}(g_i).
\end{align*}
By \eqref{eq:master-equation},
\[
\Delta_{\omega_i, \omega_i}(g_i) \Delta_{s_i\omega_i, s_i\omega_i}(g_i) =
\Delta_{s_i\omega_i, \omega_i}(g_i) \Delta_{\omega_i, s_i\omega_i}(g_i) +
\prod_{j \neq i} (\Delta_{\omega_j, \omega_j}(g_i))^{-a_{j, i}}.
\]
 As $\Delta_{\omega_i, \omega_i}(g_i)=1$ and $\Delta_{s_i\omega_i, \omega_i}(g_i)=0$, one has 
\[
p_i^\prime = z_iz_{r+r}- \prod_{j \neq i} (\Delta_{\omega_j, \omega_j}(g_i))^{-a_{j, i}} = z_iz_{r+i} - 
\prod_{j=i+1}^r z_j^{-a_{j, i}} \prod_{j = 1}^{i-1} z_{r+j}^{-a_{j, i}} =  p_i(z_1, \ldots, z_{2r}).
\]
\end{proof}

\bre{:inv-rho}
{\rm
By \eqref{eq:pi-prime}, the inverse of the isomorphism $\rho: \Lcc \to \CC^{2r}_{{\bf p} \neq 0}$ is given by
\[
\rho^{-1}: \;\; \CC^{2r}_{{\bf p} \neq 0} \longrightarrow \Lcc,\;\; (z_1, \ldots, z_{2r}) \longmapsto 
q_1(z_1)  \cdots q_r(z_r)q_1(z_{r+1}) \cdots q_r(z_{2r})\oc^{\, -1} n,
\]
where $n = \oc\, (n_2^\prime)^{-1} \oc^{\, -1} \in N$, and $n_2^\prime \in N \cap \oc^{\, -1}  N\oc$ is the unique element in the decomposition
\[
q_1(z_1)  \cdots q_r(z_r)q_1(z_{r+1}) \cdots q_r(z_{2r})= 
b_- n_1^\prime n_2^{\prime},
\]
where $b_- \in B_-$ and $n_1^\prime \in N \cap \oc^{\, -1}  N_- \oc$.  
\hfill $\diamond$
}
\ere

 In view of \thref{:Lcc-Af}, we  will denote $\calS^{\BFZ}(\Lcc)$ by $\calS^\BFZ$ and call it the {\it cluster structure on $\Lcc$ with BFZ coefficients}.
Thus $\calS^\BFZ$ has the seed 
\begin{equation}\label{eq:seed-BFZ}
\calS^\BFZ = ({\bf x}_c = (x_{\omega_1}, \, \ldots, \, x_{\omega_r}), \, {\bf p}_c = (p_{1; c}, \, \ldots, p_{r; c}), \; \wt{B}^\BFZ).
\end{equation}
We have the following immediate consequence of \prref{:BFZ-test}.

\bco{:BFZ-in-minor}
A point $g \in \Lcc$ is an $\calS^{\BFZ}$-frieze point if and only if 
\[
\Delta_{\omega_i, \omega_i} \in \ZZ_{>0},\; \;\;\Delta_{c\omega_i, c\omega_i} \in \ZZ_{>0} \;\;\; \mbox{and} \;\; \;
\Delta_{\omega_i, c\omega_i} \in \ZZ_{>0}, \hs i \in [1, r].
\]
\eco

\subsection{The cluster structure on $\Lcc$ with principal coefficients}\label{ss:Lie-prin}
Introduce
\begin{equation}\label{eq:y-ic}
y_{i; c} = p_{i; c} \prod_{j=1}^{i-1} (p_{i; c})^{a_{j, i}} \in \CC[\Lcc], \hs i \in [1, r].
\end{equation}
Modifying the seed $\Sigma^{\BFZ}$ in $\CC(\Lcc)$ given in \eqref{eq:seed-BFZ}, we have the seed
\begin{equation}\label{eq:seed-prin}
\Sigma^{\rm prin} \stackrel{\rm def}{=} ({\bf x}_c = x_{\omega_1}, \, \ldots, \, x_{\omega_r}), \, {\bf y}_c = (y_{1; c}, \, \ldots, y_{r; c}), \, \wt{B}^{\rm prin})
\end{equation} 
in $\CC(\Lcc)$, where $\wt{B}^{\rm prin}$ is given in \eqref{eq:MM}.  As ${\bf y}_c = {\bf p}_c^{\UA}$, where $\UA \in GL(r, \ZZ)$ is  given in \eqref{eq:U0-00}, the seed $\Sigma^{\rm prin}$ defines a cluster structure on $\Lcc$, which, by abusing notation, will also be denoted by $\calS^{\rm prin}$ 
and called 
the {\it cluster structure on $\Lcc$ with principal coefficients}. 

The cluster structure $\calS^{\rm prin}$ on $\Lcc$ was studied by S.-W. Yang and A. Zelevinsky \cite{YZ:Lcc}. 
Recall that \cite[Proposition 1.3]{YZ:Lcc} for each $i \in [1, r]$ one has the integer $h(i; c) \geq 1$ such that
\begin{equation}\label{eq:hic}
\omega_i > c\omega_i >  \cdots > c^{h(i; c)} \omega_i = -\omega_{i^*},
\end{equation}
where the notation $\delta > \gamma$ for $\delta \neq \gamma$ in the weight lattice  ${\mathcal{P}}$ 
means that $\delta - \gamma$ is 
a non-negative  linear integral combination of the simple roots, and  $i^* \in [1, r]$ for $i \in [1, r]$ is determined via the longest element
$w_0 \in W$ by $w_0 \omega_i = -\omega_{i^*}$.
Following \cite{YZ:Lcc}, set  
\[
\Pi(c) = \{c^m\omega_i: \; i \in [1, r], \, 0 \leq m \leq h(i; c)\}, 
\]
and for $\gamma \in \Pi(c)$, set
$x_{\gamma} = \Delta_{\gamma, \gamma}|_{\Lcc}$.

\bth{YZ:Lcc}\cite[Theorem 1.4]{YZ:Lcc} The cluster variables of the cluster structure $\calS^{\rm prin}$on $\Lcc$  are precisely the functions $x_{\gamma}$ for $\gamma \in \Pi(c)$.
\eth

\bre{:1-c}
{\rm
In \cite{YZ:Lcc}, the $g$-vector of the cluster variable $x_\gamma$ with respect to the initial seed $\Sigma^{\rm prin}$ is identified with the coefficient vector of $\gamma$ in the 
$\ZZ$-basis of $(\omega_1, \ldots, \omega_r)$ of $\calP$. 
To give a geometric meaning of the 
$g$-vectors, consider the $T$-action on $\Lcc$ given by 
\[
T \times \Lcc \longrightarrow \Lcc, \;\; t \cdot g = t g (t^{-1})^c.
\]
The restriction to $\Lcc$ of a generalized minor $\Delta_{\delta, \gamma}$, if non-zero, is then 
a $T$-weight vector with
$T$-weight $\delta-c\gamma$. In particular,  $x_\gamma$ for $\gamma \in \Pi(c)$ has $T$-weight $\gamma - c\gamma$. 
For $i \in [1, r]$, let 
\begin{equation}\label{eq:beta-i}
\beta_i = \omega_i - c\omega_i = \omega_i - s_1\cdots s_i \omega_i  = s_1\cdots s_{i-1}\alpha_i.
\end{equation}
One checks directly that (see also \cite[(2.6)]{YZ:Lcc}) that
\begin{equation}\label{eq:beta-alpha-i}
(\beta_1, \ldots, \beta_r) \UA= (\alpha_1, \ldots, \alpha_r), 
\end{equation}
where $\UA \in SL(r, \ZZ)$ is given in \eqref{eq:U0-00}. 
It follows that $(\beta_1, \ldots, \beta_r)$ is a $\ZZ$-basis for the root lattice 
$\calQ = \ZZ \alpha_1 + \ldots + \ZZ \alpha_r$, so
\[
1-c: \;\; \calP \longrightarrow \calQ, \;\; \gamma \longmapsto \gamma - c\gamma,
\]
is a lattice isomorphism. We can thus identify  the $g$-vector of $x_\gamma$ with its $T$-weight $\gamma-c\gamma$.
\hfill $\diamond$
}
\ere

We now show that for each $i \in [1, r]$, the regular function $y_{i; c}$ on $\Lcc$ is the restriction to $\Lcc$ of a generalized minor. 

\ble{:y-ic}
For each $i \in [1, r]$ there is a unique $\nu_i \in \Pi(c)$ such that $c\nu_i \in \Pi(c)$ and 
$\nu_i - c\nu_i =\alpha_i$, and one has
$y_{i; c} = \Delta_{\nu_i, \, c\nu_i}|_{\Lcc}$.
\ele

\begin{proof} 
Recall that $\Phi_+$ is the set of all positive roots. 
By \cite[Lemma 2.2]{YZ:Lcc}, one has the bijection
\[
\Pi(c) \backslash \{-\omega_i: i \in [1, r]\} \longrightarrow \Phi_+, \;\; \gamma \longmapsto \gamma - c\gamma.
\]
Thus for each $i \in [1, r]$ there is a unique $\nu_i  \in \Pi(c)$ with $c\nu_i \in \Pi(c)$ such that $\nu_i - c\nu_i = \alpha_i$. 
By \cite[Lemma 3.4]{YZ:Lcc},
$y_{i; c} =  \left( \Delta_{c\nu_i, \, \nu_i}\Delta_{\nu_i, \, c\nu_i}\right)|_{\Lcc}$.
By the proof of \cite[Lemma 3.4]{YZ:Lcc} in \cite[Page 780]{YZ:Lcc}, $\Delta_{c\nu_i, \, \nu_i}|_{\Lcc} = 1$. 
Thus $y_{i; c} = \Delta_{\nu_i, \, c\nu_i}|_{\Lcc}$.
\end{proof}

For $i \in [1, r]$, we consider now the cluster variable $x_{\omega_i}^\prime$ in 
$\mu_i(\Sigma^{\rm prin})$ that is not in $\Sigma^{\rm prin}$.  In other words, let $x_{\omega_i}^\prime$ be  the unique regular function on $\Lcc$ determined by
\begin{equation}\label{eq:x-i-prime}
x_{\omega_i} x_{\omega_i}^\prime = y_{i; c} \prod_{j=1}^{i-1} x_{\omega_j}^{-a_{j, i}} + \prod_{j=i+1}^r x_{\omega_j}^{-a_{j, i}}.
\end{equation}

\ble{:xi-prime}
 One has 
$x_{\omega_i}^\prime = \Delta_{c\omega_i, \, s_i\omega_i}|_{\Lcc}$ for every $i \in [1, r]$.
\ele

\begin{proof} Let $g \in \Lcc$, and  let $x_i = x_{\omega_i}(g)$ and $y_i = y_{i; c}(g)$. 
By \cite[(3.26)]{YZ:Lcc}, for generic $g \in \Lcc$ in the sense that $x_i \neq 0$ for every $i \in [1, r]$, one has
\[
g = e_{-1}(a_1) \cdots e_{-r}(a_r) \alpha_1^\vee(x_1) \cdots \alpha_r^\vee(x_r) e_r(y_ra_r) \cdots e_1(y_1a_1),
\]
where $a_i = x_i^{-1} \prod_{j=1}^{i-1} x_j^{-a_{j, i}}$ for $i \in [1, r]$, so 
\[
\Delta_{\omega_i, \, s_i\omega_i}(g) = \left(\alpha_1^\vee(x_1) \cdots \alpha_r^\vee(x_r)\right)^{\omega_i} 
\Delta_{\omega_i, \, s_i\omega_i}(e_r(y_ra_r) \cdots e_1(y_1a_1)) = x_iy_ia_i = y_i\prod_{j=1}^{i-1} x_j^{-a_{j, i}}.
\]
Thus for each $i \in [1, r]$ one  has
\begin{equation}\label{eq:i-si}
\Delta_{\omega_i, \, s_i\omega_i}|_{\Lcc} = y_{i; c} \prod_{j=1}^{i-1} x_{\omega_j}^{-a_{j, i}}.
\end{equation}
On the other hand, applying \eqref{eq:master-equation} with $i \in [1, r]$ and $u = s_1 \cdots s_{i-1}$, $v = e$, one has
\[
	 \Delta_{c\omega_i,\,s_i\omega_i} \Delta_{\omega_i,\,\omega_i}= 
  \Delta_{c\omega_i,\,\omega_i}\Delta_{\omega_i,\,s_i\omega_i} + \prod_{j=1}^{i-1} \Delta_{c\omega_j,\,\omega_j}^{-a_{j,i}} \prod_{j =i+1}^r
  \Delta_{\omega_j,\,\omega_j}^{-a_{j,i}}.
\]
Restricting to $\Lcc$ and using $\Delta_{c\omega_j,\,\omega_j}|_{\Lcc} = 1$, one has 
\[
x_{\omega_i} \; \Delta_{c\omega_i,s_i\omega_i}\vert_{\Lcc} = \Delta_{\omega_i,s_i\omega_i}\vert_{\Lcc} + \prod_{j=i+1}^rx_{\omega_j}^{-a_{j,i}} =
y_{i; c} \prod_{j=1}^{i-1} x_{\omega_j}^{-a_{j, i}} + \prod_{j=i+1}^rx_{\omega_j}^{-a_{j,i}}.
\]
Thus $x_{\omega_i}^\prime = \Delta_{c\omega_i, \, s_i\omega_i}|_{\Lcc}$ for ever $i \in [1, r]$.
\end{proof}

\leref{:xi-prime} now leads to the following consequence of \prref{:prin-1}.

\bco{:prin-in-minor}
A point $g \in \Lcc$ is an $\calS^{\rm prin}$-frieze point if and only if 
\[
\Delta_{\omega_i, \,\omega_i} \in \ZZ_{>0},\; \;\;\Delta_{c\omega_i, \, s_i\omega_i} \in \ZZ_{>0} \;\;\; \mbox{and} \;\; \;
\Delta_{\nu_i, c\nu_i} \in \ZZ_{>0}, \hs i \in [1, r].
\]
\eco

In the remainder of $\S$\ref{ss:Lie-prin}, we prove some more facts on $\calS^{\rm prin}$ in order to relate $\calS^{\rm prin}$-frieze points in $\Lcc$ to 
the Lusztig total positivity in $\S$\ref{ss:frieze-lusz-pos}. To this end,
consider the reduced word ${\bf c} = (s_1, s_2, \ldots, s_r)$ of $c$.
Let
$m_c = {\rm max}\{h(i; c): i \in [1, r]\}$, and 
for $0 \leq m \leq m_c$, let
$I_m(c) = \{i \in [1, r]: h(i; c) \geq m\}$. Note that
$[1, r] =I_0(c) = I_1(c) \supset I_2(c) \supset \cdots \supset I_{m_c}(c)$.
For $1 \leq m \leq m_c$, if 
 $I_m(c) = \{j_{1}, \ldots, j_p\}$ with $j_{1} < \ldots < j_p$,  set 
\begin{equation}\label{eq:wm-1}
{\bf w}_m = (s_{j_1}, \ldots, s_{j_p}) \hs \mbox{and} \hs w_m = s_{j_{1}}\cdots s_{j_p} \in W.
\end{equation}
Note that ${\bf w}_1 = (s_1, \ldots, s_r) = {\bf c}$ and $w_1 = c$.  Define 
\[
{\bf w}_0({\bf c}) = ({\bf w}_1, {\bf w}_2, \ldots, {\bf w}_{m_c}) =({\bf c}, {\bf w}_2, \ldots, {\bf w}_{m_c}).
\]
It is
proved in \cite[$\S$7]{Stella:g-fan} that ${\bf w}_0({\bf c})$ is a reduced word of the longest elment $w_0$, called \cite{Reading:Cambrian}  the {\it ${\bf c}$-sorting word of $w_0$}, and that
\begin{equation}\label{eq:cmm-1}
c^m \omega_i= w_1w_2\cdots w_m \omega_i,\hs i \in [1, r], \; 1 \leq m \leq h(i; c).
\end{equation}
Let $d = l_0 + r$, where $l_0$ is the length of $w_0$, and set
\[
{\bf c} {\bf w}_0({\bf c}) = ({\bf c}, \, {\bf w}_0({\bf c})) = (s_{i_1}, \, \ldots, \,s_{i_r}, \,s_{i_{r+1}}, \ldots, \,s_{i_d}),
\]
where $(s_{i_1}, \ldots, s_{i_r}) = (s_1, \ldots, s_r)$.
For $k \in [1, d]$, define
$\gamma_{k}= s_{i_1}\cdots s_{i_{k-1}} \omega_{i_k} \in \calP$.
\ble{:mk}
One has $\{\gamma_k: k \in [1, d]\} = \Pi(c)$. 
\ele
\begin{proof}
Define the {\it block number} of  $k \in [1, d]$ to be $m(k) = 0$ if $k \in [1, r]$, and for $k \in [r+1, d]$, 
\begin{equation}\label{eq:block-mk}
m(k) = m\in [1, m_c] \;\;\; \mbox{if}\;\; r+l(w_1\cdots w_{m-1}) +1 \leq k \leq r + l(w_1\cdots w_{m-1}w_m).
\end{equation}
We prove that $\gamma_k = c^{m(k)}\omega_{i_k} \in \Pi(c)$ for every $k \in [1, d]$. Indeed, 
the statement holds for $k \in [1, r]$. Let $k \in [r+1, d]$, so that $m(k) \geq 1$, and let
 $k' =r+l(w_1\cdots w_{m(k)-1})$. Then
\[
\gamma_{k}= cw_1 \cdots w_{m(k)-1} s_{i_{k^{\prime}+1}} \cdots s_{i_{k-1}} \omega_{i_k} =
cw_1 \cdots w_{m(k)-1}  \omega_{i_k} =c^{m(k)} \omega_{i_k} \in \Pi(c).
\]
For every $i \in [1, r]$ and $m \in [0, h(i; c)]$ there is a unique $k \in [1, d]$ with $i_k = i$ and $m(k) = m$, so $c^m \omega_i = \gamma_k$. Thus the map $[1, d]\to\Pi(c), k \mapsto \gamma_k$, is bijective.
\end{proof}
Recall from \eqref{eq:ga} the element $q_i(z) \in G$ for $i \in [1, r]$ and $z \in \CC$.
\bpr{:BS-coordinates}
The map $g_{\cwc}: \CC^d \rightarrow N \oc N$ given by
\[
g_{\cwc}(x_1,\ldots, x_d) = q_{i_1}(x_1)  \cdots q_{i_r}(x_r)  q_{i_{r+1}}(x_{i+1})  \cdots q_{i_d}(x_d)\overline{w}_0^{\, -1}
\]
is an isomorphism, and for every $x = (x_1, \ldots, x_d) \in \CC^d$, one has
\[
x_k = \Delta_{\gamma_k, \gamma_k}(g_{\cwc}(x)), \hs k \in [1, d].
\]
\epr
\begin{proof}
By the isomorphism $N_c \times N \to N\oc N, (n, n') \mapsto n \oc n'$ and the isomorphisms
\begin{align*}
&\CC^r \longrightarrow N_c\oc,\;\; (x_1 \ldots, x_r) \longmapsto q_{i_1}(x_1)  \cdots q_{i_r}(x_r),\\
& \CC^{l_0} \longrightarrow N, \;\;  (x_{r+1},  \ldots, x_d) \longmapsto q_{i_{r+1}}(x_{i+1})  \cdots q_{i_d}(x_d) \overline{w}_0^{\, -1},
\end{align*}
we see that $g_{\cwc}$ is an isomorphism. Fix now $x = (x_1, \ldots, x_{d}) \in \CC^{d}$, and let 
\[
n(x) = q_1(x_1) \cdots q_r(x_r) \overline{c}^{\, -1}\in N_c \hs \mbox{and} \hs n^\prime(x)= 
q_{i_{r+1}}(x_{r+1}) \cdots q_{i_{d}}(x_{d})  \overline{w_0}^{\, -1}\in N,
\]
so that $g_{\cwc}(x) = n(x)\oc\, n'(x)$. For $k \in [1, r]$, one  has $\gamma_k = \omega_k$ and
\[
\Delta_{\omega_k, \, \omega_k}(g_{\cwc}(x)) = \Delta_{\omega_k, \, \omega_k}(n(x)\oc\, n'(x)) = 
\Delta_{\omega_k, \, \omega_k}(n(x)\oc) = 
\Delta_{\omega_k, \, \omega_k}(q_1(x_1) \cdots q_r(x_r)) = x_k,
\]
where in the last step we used \leref{:gig}. Let now $k \in [1, d]$, and let 
$m(k)\geq 1$ be the block number of $k$ given in \eqref{eq:block-mk}. Let $w^\prime_{m(k)} = w_2\cdots w_{m(k)}$ (or $e$ if $m(k) = 1)$.  By the proof of \leref{:mk}, 
$\gamma_{k} = c^{m(k)} \omega_{i_k} = cw^\prime_{m(k)} \omega_{i_k}$, and 
\begin{equation}\label{eq:u-omegak}
w^\prime_{m(k)}\omega_{i_k} = c^{m(k)-1} \omega_{i_k} =
w_1 w_2 \cdots w_{m(k)-1} \omega_{i_k}, 
\end{equation}
where recall that $w_1 = c$ and $w_1 w_2 \cdots w_{m(k)-1}= e$ if $m(k) = 1$. Thus
\begin{align*}
\Delta_{\gamma_{k}, \,\gamma_{k}}(g_{\cwc}(x)) & = \Delta_{cw^\prime_{m(k)}\omega_{i_k}, \,cw^\prime_{m(k)}\omega_{i_k}}(n(x)\oc n'(x)) 
 = \Delta_{w^\prime_{m(k)}\omega_{i_k}, \, cw^\prime_{m(k)}\omega_{i_k}}(n'(x))\\
&= \Delta_{w_1\cdots w_{m(k)-1}\omega_{i_k}, \,w_1 \cdots w_{m(k)}\omega_{i_k}}(n'(x)).
\end{align*}
Writing $n'(x)\overline{w_0} = (n_1(x) \overline{w_1}) \, (n_2(x) \overline{w_2}) \, \cdots \, (n_{m_c} (x)\overline{w_{m_c}})$,
where $n_m(x) \in N_{w_m}$ for $m \in [1, m_c]$. By \leref{le:Nww} and \leref{:gig}, one has
\[
\Delta_{\gamma_{k},\, \gamma_{k}}(g_{\cwc}(x))  = \Delta_{\omega_{i_k}, \,\omega_{i_k}}(n_{m(k)}\overline{w_{m(k)}}) = x_{k}.
\]
\end{proof}
In the setting of \prref{:BS-coordinates}, the functions $x_1, x_2, \ldots, x_d$ on $N\oc N$ via the isomorphism $g_{\cwc}$ will be called the {\it Bott-Samelson coordinates} on $N \oc N$ associated to $\cwc$. 
\bco{:BS-coor}
The restrictions to $\Lcc\subset N\oc N$ of the Bott-Samelson coordinates $x_1, \ldots, x_{d}$ on $N\oc N$ are precisely
all the cluster variables of the cluster structure $\calS^{\rm prin}$ on $\Lcc$.
\eco

\ble{:minor-in-BS}
The restriction  of every generalized minor to $N \oc N$ is a polynomial in the Bott-Samelson coordinates $x_1, \ldots, x_d$ on $N \oc N$ with integer coefficients.
\ele

\begin{proof}
Let $\Delta = \Delta_{\delta, \gamma}$ be any generalized minor of $G$. Recall from \eqref{eq:si} that  $\overline{s_i} = e_i(-1) e_{-i}(1) e_i(-1)$. Using any reduced word for 
$w_0$, one can also write $\overline{w_0}^{\, -1}$ as a product of $3l_0$ elements, each of which is of the form $e_i(\pm 1)$
or $e_{-i}(\pm 1)$ for some $i \in [1, r]$.
It follows from \leref{:minor-ui} that there exists a polynomial 
$P \in \ZZ_{\geq 0}[x_1, x_2, \ldots, x_{d}, x_1^\prime, x_2^\prime, \ldots, x^\prime_{3d+3l_0}]$
such that $\Delta_{\gamma, \delta}(g_{\cwc}(x_1, \ldots, x_{d}))$ is obtained from $P$ by setting each $x_k^\prime$, for $k \in [1, 3d+3l_0]$, 
to be equal to $1$ or $-1$. It follows that $\Delta_{\gamma, \delta}(g_{\cwc}(x_1, \ldots, x_{d})) \in \ZZ[x_1, \ldots, x_{d}]$.
\end{proof}

\subsection{$\calS^{\rm prin}$-frieze points and Lusztig total positivity}\label{ss:frieze-lusz-pos}
In \cite{Lusztig1994}, G. Lusztig defined the {\it totally non-negative part} of $G$ to be the multiplicative monoid $G_{\geq 0} \subset G$ generated by 
the set
\[
\{e_i(a), \; e_{-i}(b), \; \alpha_i^\vee(c): \; i \in [1,r], \; a, b, c \in \RR_{>0}\}.
\]
By \cite[Theorem 3.1]{FZ:Oscillatory},
$g \in G$ lies in $G_{\geq 0}$ if 
and only if $\Delta(g) \geq 0$ for every generalized minor $\Delta$ of $G$. Let $G(\ZZ)$ be the sub-group of $G$ generated by 
$\{e_i(a), \; e_{-i}(b): \; i \in [1,r], \; a, b \in \ZZ\}$. 

\bth{:frieze-pos}
For $g \in \Lcc$, the following are equivalent:

1) $g$ is an $\calS^{\rm prin}$-frieze point;

2) $\Delta(g) \in \ZZ_{\geq 0}$ for every generalized minor $\Delta$ of $G$;

3) $g \in G_{\geq 0} \cap G(\ZZ)$.
\eth

\begin{proof}
By \cite[Theorem 1.11]{FZ:double}, one has $\Lcc \cap G_{\geq 0} = L^{c, c^{-1}}_{>0}$, where
\begin{align*}
L^{c, c^{-1}}_{>0} &= \{g \in \Lcc: \; x_{\omega_i}(g) >0, \; p_{i; c}(g) >0, \; \forall \, i \in [1, r]\}\\
& = \{g \in \Lcc: \; x_{\omega_i}(g) >0, \; y_{i; c}(g) >0, \; \forall \, i \in [1, r]\}.
\end{align*}

Assume 1).  Then $g \in L^{c, c^{-1}}_{>0}$, so $g \in G_{\geq 0}$. 
For every generalized minor $\Delta$ of $G$, one then has $\Delta(g) \geq 0$. 
By \leref{:minor-in-BS}, $\Delta(g) \in \ZZ$. Thus $\Delta(g) \in \ZZ_{\geq 0}$. This shows that  1) implies 2).

Assume 2). Then $g \in G_{\geq 0}$ and by \coref{:BS-coor} all the Bott-Samelson coordinates of $g \in \Lcc \subset N \oc N$ are non-negative
integers. Since $\overline{s_i} \in G(\ZZ)$ for each $i \in [1, r]$, we have $g \in g_{\cwc}(\ZZ^d) \subset  G(\ZZ)$. Thus 2) implies 3).

Assume 3). Then $g \in \Lcc \cap G_{\geq 0} = L^{c, c^{-1}}_{>0}$, so 
$x_{\omega_i}(g) >0$ and  $y_{i; c}(g) >0$ for every $i \in [1, r]$.
By the Laurent Phenomenon, 
 $x_\gamma(g) >0$ for every $\gamma \in \Pi(c)$. As $g \in G(\ZZ)$, it follows from \leref{:minor-ui} that 
$x_{\gamma}(g) \in \ZZ_{>0}$ for every $\gamma \in \Pi(c)$. By \leref{:y-ic}, one also has $y_{i; c} (g) \in \ZZ_{>0}$. 
Thus $g$ is an $\calS^{\rm prin}$-frieze point. Hence 3) implies 1). 
\end{proof}

\bre{:int}
{\rm
Let $G^{\rm int}$ be the set of all $g \in G$ such that $\Delta(g) \in \ZZ$ for every generalized minor $\Delta$, and let $G_{\geq 0}^{\rm int}$ the set of all $g \in G$ such that $\Delta(g) \in \ZZ_{\geq 0}$ for every generalized minor $\Delta$.
Let  $G(\ZZ_{>0})$ be the sub-monoid of $G(\ZZ)$ generated by
the set $\{e_i(a), e_{-i}(b):\; i \in [1,r], \; a, b \in \ZZ_{>0}\}$. By \leref{:minor-ui}, one has
\[
G(\ZZ) \subset G^{\rm int} \hs \mbox{and} \hs 
G(\ZZ_{>0}) \subset G_{\geq 0}^{\rm int}.
\]
When $G = {\rm SL}(n, \CC)$ with the standard pinning, one has $G^{\rm int} ={\rm SL}(n, \ZZ)$, and since ${\rm SL}(n, \ZZ)$ is generated by the elementary matrices (see e.g. \cite[1.2.11]{Hahn-Omeara:k-theory}), one also has $G(\ZZ) = G^{\rm int}$. In general, it would be interesting to further study the sets $G^{\rm int}$ and $G(\ZZ)$, particularly in connection with frieze points in arbitrary double Bruhat cells.
\hfill $\diamond$ 
}
\ere

\subsection{Periodicity of cluster variables on the acyclic belt}\label{ss:belt-finite}
Continuing to assume that $A$ is of finite type, we further assume in this section that $A$ is 
indecomposable, so that  if $h$ is the Coxeter number of $A$, i.e., if 
$h$ is the order of $c \in W$, then  \cite[Proposition 1.7]{YZ:Lcc}, 
\begin{equation}\label{eq:hhic}
h(i; c) + h(i^*; c) = h, \hs i \in [1, r].
\end{equation}
We consider the acyclic belt $\{\Sigma_{t(i, m)}: (i, m) \in [1, r] \times \ZZ\}$ of $\calS^{\rm prin}$ with
$\Sigma_{t(1, 0)} = \Sigma^{\rm prin}$ as in \eqref{eq:seed-prin}. For this special case of 
principal coefficients,   we write, for $(i, m) \in [1, r] \times \ZZ$,
\[
\Sigma_{t(i, m)} = ({\bf x}_{t(i, m)}, \, {\bf y}_c, \; \wt{B}_{t(i, m)}) = (B_{t(i, m)}\bb C_{t(i, m)}),
\]
and denote the $i$th cluster variable  of ${\bf x}_{t(i, m)}$ by $x(i, m)$. We call the assignment 
\begin{equation}\label{eq:xim-pattern}
[1, r] \times \ZZ \longrightarrow \fX(\calS^{\rm prin}), \;\; (i, m) \longmapsto x(i, m),
\end{equation}
the {\it principal generic frieze pattern} associated to $A$ (see \deref{:generic-frieze}).
Let $c(i, m)$ be the $i$th column of $C_{t(i, m)}$, so by 
\eqref{eq:ex-0} each $(i, m) \in [1, r] \times \ZZ$ gives  the mutation relations
\begin{equation}\label{eq:xmi-cmi}
x(i, m)\, x(i, m+1) =  {\bf y}_c^{[c(i, m)]_+} + {\bf y}_c^{[-c(i, m)]_+}\prod_{j=i+1}^{r} x(j, m)^{-a_{j,i}} \prod_{j=1}^{i-1} x(j,m+1)^{-a_{j,i}}.
\end{equation} 
             In this section we state \thref{:finite-main}, which identifies the cluster variables $x(i, m)$ using weights in $\Pi(c)$, 
describes the mutation relations in \eqref{eq:xmi-cmi}, and determines the periodicity of the principal generic frieze pattern
\eqref{eq:xim-pattern}. The proof of \thref{:finite-main} will be given in $\S$\ref{ss:proof-finite-main}.

To prepare for the statement of \thref{:finite-main}, we introduce some notation and recall some facts from \cite{YZ:Lcc}. 
We first note that by \cite{Keller:periodicity} one has $(\mu_r \cdots \mu_2 \mu_1)^{h+2} \;\Sigma^{\rm prin}= \Sigma^{\rm prin}$, hence
\begin{equation}\label{eq:hplus2}
\Sigma_{t(i, m+h+2)} = \Sigma_{t(i, m)}, \hs \forall \; (i, m) \in [1, r] \times \ZZ.
\end{equation}
It thus suffices to consider the cluster variables $x(i, m)$ and the mutation relations in \eqref{eq:xmi-cmi} for
$(i, m) \in [1, r] \times [0, h+1]$. Introduce
\begin{equation}\label{eq:frakF}
\bfF: \; [1, r] \times \ZZ \longrightarrow [1, r] \times \ZZ, \;\; \bfF(i, m) = (i^*, \;m+h(i^*; c) +1).
\end{equation}
It follows from \eqref{eq:hhic} that 
\begin{equation}\label{eq:F-2}
\bfF^2(i, m) = (i, \,m+h+2), \hs (i, m) \in [1, r]\times \ZZ.
\end{equation}
For $i, j \in [1, r]$, write $i \prec_c j$ if $i < j$ and $a_{i, j} \neq 0$. By \cite[Proposition 1.6]{YZ:Lcc}, 
\begin{equation}\label{eq:h-ij}
h(i; c) - h(j; c) = \begin{cases} 1, & \hs \mbox{if}\;\;i \prec_c j \;\;\,\mbox{and}\; \;\,j^* \prec_c i^*,\\ 0, & \hs \mbox{if} \;\;i \prec_c j \;\;\,\mbox{and}\; \;\, i^* \prec_c j^*.\end{cases}
\end{equation}
For $i \in [1, r]$ and $m \in [0, h(i; c)-1]$, the following are proved in  \cite[(1.10) and (1.11)]{YZ:Lcc}:
\begin{align}\label{eq:110}
x_{-\omega_i} x_{\omega_i} &= 1 + y_{i;c} \prod_{i \prec_c j} x_{-\omega_j}^{-a_{j, i}} \prod_{j \prec_c i}x_{\omega_j}^{-a_{j, i}},\\
\label{eq:111}
x_{c^m\omega_i} x_{c^{m+1}\omega_i} &= {\bf y}_c^{[c^m\beta_i:\,\underline{\alpha}]} + \prod_{i \prec_c j} x_{c^m\omega_j}^{-a_{j, i}} 
\prod_{j \prec_c i}x_{c^{m+1}\omega_j}^{-a_{j, i}},
\end{align}
where $\beta_i$ for $i \in [1, r]$ is given in \eqref{eq:beta-i}, and
$[c^m\beta_i:\underline{\alpha}]$ is the 
integral (column) coefficient vector of the root $c^m\beta_i$ in its 
expansion in $\underline{\alpha} = (\alpha_1, \ldots, \alpha_r)$. 
Note that \eqref{eq:h-ij} implies that the weights appearing on the right hand side of \eqref{eq:111} are in $\Pi(c)$. 

For $(i, m) \in [1, r] \times [0, h+1]$, define
\begin{equation}\label{eq:gmi}
\gamma(i, m) = \begin{cases} c^m \omega_i, & \hs  m \in [0, \,h(i; c)],\\
c^{m-h(i; c)-1} \omega_{i^*}, & \hs m \in [h(i; c)+1,\, h+1].\end{cases}
\end{equation}

\bth{:finite-main} For every $(i, m) \in [1, r] \times [0, h+1]$, one has $x(i, m) = x_{\gamma(i, m)}$, and

1) if $m \in [0, h(i; c)-1]$, the mutation relation in \eqref{eq:xmi-cmi}  is \eqref{eq:111}; 

2) if $m = h(i; c)$, the mutation relation in \eqref{eq:xmi-cmi}  is \eqref{eq:110} with $i$ replaced by $i^*$; 

3) if $m \in [h(i; c) +1, \, h]$, the mutation relation in \eqref{eq:xmi-cmi}  is \eqref{eq:111} with $i$ replaced by $i^*$ and $m$ replaced by $m-h(i; c)-1$;

4) if $m = h+1$, the mutation relation in \eqref{eq:xmi-cmi}  is \eqref{eq:110}.

\noindent
Moreover, one has $x(\bfF(i, m)) = x(i, m)$ for all $(i, m) \in [1, r] \times \ZZ$.
\eth

\bre{:D}
{\rm
Note that we have the disjoint union $[1, r]\times [0, h+1] = \calD \sqcup \bfF(\calD)$, where
\begin{equation}\label{eq:D}
\calD = \{(i,m): \, i \in [1, r], \, m \in [0, h(i; c)]\}.
\end{equation}
Combining with \eqref{eq:F-2}, one sees that $\calD$ is a fundamental domain of $\bfF:
[1, r] \times \ZZ \to [1, r] \times \ZZ$. It follows from \thref{:finite-main} and the bijection
\[
\calD \longrightarrow \Pi(c):\;\; (i, m) \longmapsto \gamma(i, m) = c^m\omega_i,
\]
that each cluster variable of $\calS^{\rm prin}$ 
appears exactly once in the domain $\calD$ 
and exactly twice in the domain $[1, r] \times [0, h+1]$ for the principal generic frieze pattern $(i, m)\mapsto x(i, m)$.
\hfill $\diamond$
}
\ere

\bre{:primitive}
{\rm
Recall from \cite[$\S$12]{FZ:IV}  and \cite[Page 877]{YZ:Lcc} that 
a mutation relation \eqref{eq:xk} is said to be
{\it primitive} if one of the two products of cluster variables on the right-hand side is  equal to $1$. It is shown in \cite[Theorem 1.5]{YZ:Lcc} that
\eqref{eq:110} and \eqref{eq:111} are exactly all the primitive mutation relations in the cluster structure $\calS^{\rm prin}$ on $\Lcc$. While it is clear that
all the mutation relations in \eqref{eq:xmi-cmi}, coming from source mutations, are primitive,
 \thref{:finite-main} shows that all primitive mutation relations in $\calS^{\rm prin}$
appear in \eqref{eq:xmi-cmi}.
\hfill $\diamond$
}
\ere

\bco{:uim-F-inv}
Assume that the Cartan matrix $A$ is of finite type and indecomposable. For any $P \in M_{l, r}(\ZZ)$, the generic frieze pattern 
\[
[1, r]\times \ZZ \longrightarrow \fX(\BA \bb P), \;\; (i, m) \longmapsto u(i, m),
\]
is surjective, $\bfF$-invariant, and has $\calD \subset [1, r] \times \ZZ$ in \eqref{eq:D} as a fundamental domain. 
\eco

\begin{proof}
The statements are direct consequence of the $\bfF$-invariance of the principal generic frieze pattern \eqref{eq:xim-pattern} and 
the  Separation Formula \cite[Theorem 3.7]{FZ:IV}.
\end{proof}


\noindent
{\it Proof of \thref{:F-inv-all}}. Assume that the Cartan matrix $A$ is of finite type and indecomposable. Let 
$P \in M_{l, r}(\ZZ)$ be arbitrary, and let
$f: [1, r] \times \ZZ \to \ZZ$ be a frieze pattern associated to $(\BA \bb P)$.
By \coref{:finite-equivalent}, $f$ is obtained by evaluation of a frieze $h: \calA_\ZZ((\BA \bb P), \emptyset) \to \ZZ$ on the generic frieze pattern associated to $(\BA \bb P)$. By \coref{:uim-F-inv}, $f$ is $\bfF$-invariant.

\bre{:F-w0=-1}
{\rm 
When the Cartan matrix $A$ is of type $A_1, B,C,D_{2n}, n \geq 2,  E_7,E_8, F_4$ or $G_2$, 
the longest element $w_0$ in $W$ acts as the negative of the identity operator on $\t$, so $i = i^*$ and 
$h(i; c) = \frac{1}{2} h$ for every $i \in [1, r]$. 
\hfill $\diamond$
}
\ere

Suppose now that $A$ is of type $A_n, D_{2n+1}$, $n \geq 2$, or $E_6$. 
Let $Q$ be the Dynkin quiver defined by $\BA$, i.e., $Q$ has vertex set $[1, r]$, and for $i \neq j$ there is one arrow from $i$ to $j$ if and only if 
$i < j$ and $a_{i, j} =-1$. The map  $i \mapsto i^*$ is the permutation of the vertex set $[1,r]$ induced by the unique graph automorphism of the Dynkin diagram underlying $Q$. For
 $i \in [1,r]$, let 
 \[
 \pi(i, i^*): \; \;i = i_0 \mathdash \;i_1 \mathdash 
 \;i_2 \mathdash  \cdots \mathdash \; i_{l-1}  \mathdash \;i_l = i^*
 \]
 be the {\it unique path in the Dynkin diagram connecting $i$ to $i^*$}, and write $l = l(i, i^*) \geq 0$, the number of edges in the path $\pi(i, i^*)$. If $i = i^*$, we have $l(i, i^*) = 0$ and set $a_i = b_i = 0$. If $i \neq i^*$, 
  let $a_i$ be the 
number of arrows in $Q$, whose underlying edges lie in $\pi(i, i^*)$, that are pointing towards $i^*$, and let $b_i$ be the number of arrows in $Q$, whose underlying edges lie in $\pi(i, i^*)$, that are pointing towards $i$. In other words,
\begin{align*}
    a_i  =\#\{j \in [0, l-1]: \, i_j \prec_c i_{j+1}\} \hs \mbox{and} \hs 
    b_i = \#\{j \in [0, l-1]: \, i_{j+1} \prec_c i_{j}\}.
\end{align*}
In \cite{Bedard:commu-classes} R. B\'{e}dard introduced the numbers
\[
m_i = \frac{1}{2}(h + a_i - b_i), \hs i \in [1, r].
\]

\ble{:hic-mi} For $A$ of type $A_n, D_{2n+1}$, $n \geq 2$, or $E_6$, one has 
$h(i;c) = m_i$ for every $i \in [1, r]$.
\ele

\begin{proof}
If $i = i^*$, then $a_i = b_i = 0$ and the statement holds.
Suppose that $l(i, i^*) = 1$ and $i \prec_c i^*$. By \eqref{eq:hhic} and \eqref{eq:h-ij}, $h(i; c) = \frac{1}{2}(h+1)$. On the other hand, $a_i = 1$ and $b_i = 0$, so the statement holds. Similarly the statement holds if $l(i, i^*) = 1$ and $i^* \prec_c i$. 

Suppose now that $l(i, i^*) \geq 2$, and assume that the statement holds for all $k \in [1, r]$ such that $l(k, k^*) \leq  l(i, i^*)-2$. Applying $-w_0$ to the vertices of the path $\pi(i, i^*)$ and using the uniqueness of the path $\pi(i, i^*)$, one sees that $i_j^* = i_{l-j}$ for every $j \in [0, l]$. In particular, we have
\[
\pi(i_1, i_1^*):\;\; i_1 \mathdash \; 
 i_2 \mathdash \; \cdots \; \mathdash\; i_{l-1}.
 \]
 
{\bf Case 1:} $i_0 \prec_c i_1$ and $i_{l-1} \prec_c i_l$. By \eqref{eq:h-ij} and the induction assumption, one has
\[
h(i; c) = h(i_1; c)+1 =\frac{1}{2}(h + a_i-2 -b_i)+1 = m_i;
\]

{\bf Case 2:} $i_0 \prec_c i_1$ and $i_l \prec_c i_{l-1}$. By \eqref{eq:h-ij} and the induction assumption, one has
\[
h(i; c) = h(i_1; c) =\frac{1}{2}(h + a_i-1 -(b_i-1)) = m_i;
\]

{\bf Case 3:} $i_1 \prec_c i_0$ and $i_{l-1} \prec_c i_l$. By \eqref{eq:h-ij} and the induction assumption, one has
\[
h(i; c) = h(i_1; c) =\frac{1}{2}(h + a_i-1 -(b_i-1)) = m_i;
\]

{\bf Case 4:} $i_1 \prec_c i_0$ and $i_l \prec_c i_{l-1}$. By \eqref{eq:h-ij} and induction assumption, one has
\[
h(i; c) = h(i_1; c)-1 =\frac{1}{2}(h + a_i -(b_i-2)) -1= m_i.
\]
This shows that $h(i; c) = m_i$ for every $i \in [1, r]$.
\end{proof}

\bre{:ADE-F} 
{\rm When $A$ is of type $A_n$,  $D_n$, $E_7$ or $E_8$ with {\it standard labelling} (see e.g. \cite[p.43]{Kac:inf-dim-Lie-alg}), our formulas for $\bfF$ recover those in \cite[p.909]{Sophie-M:survey}. When $A$ is of type $E_6$, we have  
    \[
    \bfF(i,m) = \begin{cases}
        (5-i,\, m + i+4), & \hs  i = 1,2, \ldots , 5, \\
        (6,\, m + 7), & \hs i = 6,
    \end{cases}
    \]
which is different from the formula for $\bfF$ in \cite[p.909]{Sophie-M:survey}: the latter is based on an incorrect formula for the Nakayama permutation written down in \cite[p.48]{dlab-gabriel:rep1}.
\hfill $\diamond$
}
\ere

\subsection{Proof of \thref{:finite-main}}\label{ss:proof-finite-main}
To prepare  for the proof of  \thref{:finite-main}, we first observe that 
by comparing the $T$-weights of the functions on the left hand side of \eqref{eq:111} and of the second term of the right hand side of \eqref{eq:111}, one has,
for all $i \in [1, r]$ and $m \in [0, h(i; c)-1]$,
\begin{align}\label{eq:T-wt-0}
c^m \beta_i + c^{m+1}\beta_i &= \sum_{i \prec_c j} (-a_{j, i}) c^m \beta_j + \sum_{j\prec_c i} (-a_{j, i}) c^{m+1}\beta_j\\
\label{eq:T-wt-1}
& = \sum_{j \in J_1(i) \sqcup J_2(i)} (-a_{j, i}) c^m \beta_j + \sum_{j \in J_3(i) \sqcup J_4(i)} (-a_{j, i}) c^{m+1}\beta_j,
\end{align}
where for $i \in [1, r]$, 
\begin{align}\label{eq:JJ-0}
J_1(i) &= \{j \in [1, r]: \, i \prec_c j, \, i^* \prec_c j^*\}, \hs J_2(i) = \{j \in [1, r]: \, i \prec_c j, \, j^* \prec_c i^*\},\\
\label{eq:JJ-1}J_3(i) &= \{j \in [1, r]: \, j \prec_c i, \, j^* \prec_c i^*\}, \hs J_4(i) = \{j \in [1, r]: \, j \prec_c i, \, i^* \prec_c j^*\}.
\end{align}
By \eqref{eq:h-ij}, for $j \in [1, r]$ and $a_{i, j} \neq 0$, one has
\begin{equation}\label{eq:h-ij-1}
h(j; c) = \begin{cases} h(i; c), &\hs  j \in J_1(i) \sqcup J_3(i),\\
h(i; c)-1, & \hs j \in J_2(i),\\
h(i; c)+1, & \hs j \in J_4(i).\end{cases}
\end{equation}

We next determine $c(i, m) \in \ZZ^r$ in \eqref{eq:xmi-cmi} for $(i, m) \in [1, r]\times \ZZ$. By \eqref{eq:hplus2}, 
$c(i, m+h+2) = c(i, m)$ for all  $(i,m) \in [1, r] \times \ZZ$.
For $(i, m) \in [1, r] \times [0, h+1]$, define $\alpha(i, m) \in \Phi$ by
\[
\alpha(i, m)  = \begin{cases} c^m \beta_i, & \hs m \in [0, \,h(i; c)-1],\\
 -\alpha_{i^*}, & \hs m = h(i; c),\\
c^{m-h(i; c)-1}\beta_{i^*}, & \hs m \in [h(i; c)+1, \,h],\\
-\alpha_i, & \hs m = h+1.\end{cases}
\]
Let again $\underline{\alpha} = (\alpha_1, \ldots, \alpha_r)$, and for a root $\beta$ let
$[\beta:\underline{\alpha}] \in \ZZ^r$ be such that $\beta = \underline{\alpha} [\beta:\underline{\alpha}]$.

\ble{:cim-0}
One has $c(i, m) = [\alpha(i, m):\underline{\alpha}]$ for $(i, m) \in  [1, r] \times [0, h+1]$.
\ele

\begin{proof} 
Note first that for any $i \in [1, r]$, we have $c(i, h+1) = c(i, -1)$, and by taking $P = I_r$ and $m = -1$ in \eqref{eq:pmi-1}, we have
$c(i, -1) = (0, \ldots, 0, -1, 0, \ldots, 0)^T$, where $-1$ is at the $i$th place. As $\alpha(i, h+1) = -\alpha_i$, we have $c(i, m) = [\alpha(i, m):\underline{\alpha}]$
for $i \in [1, r]$ and $m = h+1$.

We now prove \leref{:cim-0} for $(i, m) \in [1, r] \times [0, h]$ using
induction on $(i, m) \in [1, r] \times [0, h]$ in the total order on $[1, r]\times \ZZ$ defined in 
\eqref{eq:order-0}.

Consider first the case of $m = 0$ and recall the matrix $\UA$ from \eqref{eq:U0-00}. 
Taking $P = I_r$ in \eqref{eq:pmi-1} and noting that $\UA^{-1} \geq 0$,  one sees that $(c(1, 0), \ldots, c(r, 0)) = \UA^{-1}$.
 By \eqref{eq:beta-alpha-i}, 
$(\beta_1, \ldots, \beta_r) = \underline{\alpha} \UA^{-1}$. Since $\alpha(i, 0) = \beta_i$, one has
$c(i, 0) = [\beta_i:\underline{\alpha}] = [\alpha(i, 0):\underline{\alpha}]$ for every $i \in [1, r]$.

Let now $(i, m) \in [1, r] \times [0, h-1]$ and assume that $c(j, m') =[\alpha(j, m'):\underline{\alpha}]$ for all 
$j \in [1, r]$ and $m' \in [0, h]$ such that 
$(j, m') < (i,m+1)$.
By \prref{:additive}, 
\begin{equation}\label{eq:cmi-11}
c(i,m) + c(i, m+1) = \sum_{i \prec_c j} (-a_{j, i})[c(j, m)]_+ + \sum_{j \prec_c i} (-a_{j, i}) [c(j, m+1)]_+.
\end{equation}
As $(j, m) < (i, m+1)$ for all $j \in [i, r]$, and $(j, m+1) < (i, m+1)$ for all $j \in [1, i-1]$, we can 
apply the induction assumption on $c(i, m)$ and on the terms appearing on the right hand side of \eqref{eq:cmi-11} to determine $c(i,m+1)$. 
We consider four cases.

{\bf Case 1:} $m \in [0, h(i; c)-2]$. By \eqref{eq:h-ij-1}, for $j \in [1, r]$, 
if $i \prec_c j$  then  $m \leq h(j; c)-1$, so $\alpha(j, m) = c^m\beta_j \in \Phi_+$, and if
$j \prec_c i$, then $m+1 \leq h(j; c)-1$, so $\alpha(j, m+1) = c^{m+1}\beta_j \in \Phi_+$.
Applying the induction assumption to the terms in \eqref{eq:cmi-11}, one gets
\[
[c^m \beta_i:\underline{\alpha}] + c(i, m+1) = \sum_{i \prec_c j} (-a_{j, i})[c^m\beta_j:\underline{\alpha}] + \sum_{j \prec_c i} (-a_{j, i}) [c^{m+1}\beta_j:\underline{\alpha}].
\]
Comparing with \eqref{eq:T-wt-0}, one sees that $c(i, m+1) = [c^{m+1}\beta_i:\underline{\alpha}] = [\alpha(i, m+1):\underline{\alpha}]$.

{\bf Case 2:} $m = h(i; c)-1$. In this case, for $j \in [1, r]$, by \eqref{eq:h-ij-1} one has 
\[
\begin{cases}
m = h(j; c)-1 \hs \mbox{and} \hs \alpha(j, m) = c^m \beta_j \in \Phi_+, & \hs j \in J_1(i),\\
m = h(j; c) \hs \mbox{and} \hs \alpha(j, m) = -\beta_{j^*} \in -\Phi_+, & \hs j \in J_2(i),\\
m +1= h(j; c) \hs \mbox{and} \hs \alpha(j, m+1) = -\beta_{j^*} \in -\Phi_+, & \hs j \in J_3(i),\\
m +1= h(j; c)-1 \hs \mbox{and} \hs \alpha(j, m+1) = c^{m+1} \beta_j \in \Phi_+, & \hs j \in J_4(i).\end{cases}
\]
Applying the induction assumption to the terms in \eqref{eq:cmi-11}, one has
\begin{align*}
[c^m\beta_i:\underline{\alpha}] + c(i, m+1)& = \sum_{j \in J_1(i) \sqcup J_2(i)} (-a_{j, i})[c(m, j)]_+ 
\sum_{j \in J_3(i) \sqcup J_4(i)} (-a_{j, i}) [c(m+1), j]_+ \\
& =  \sum_{j \in J_1(i)} (-a_{j, i})[c^m\beta_j:\underline{\alpha}] + \sum_{j \in J_4(i)} (-a_{j, i})[c^{m+1}\beta_j:\underline{\alpha}].
\end{align*}
Comparing with \eqref{eq:T-wt-1} and using $a_{j^*, i^*} = a_{j, i}$ for $j \in J_2(i) \sqcup J_3(i)$ and \eqref{eq:beta-alpha-i},  one has
\begin{align*}
c(i, m+1) &= -[\beta_{i^*}:\underline{\alpha}] +\sum_{j \in J_2(i) \sqcup J_3(i)} (-a_{j, i})([\beta_{j^*}:\underline{\alpha}] \\
& = -[\beta_{i^*}:\underline{\alpha}] -\sum_{j^* \prec_c i^*} a_{j^*, i^*}([\beta_{j^*}:\underline{\alpha}] = -[\alpha_{i^*}:\underline{\alpha}] = 
[\alpha(i, m+1):\underline{\alpha}].
\end{align*}

{\bf Case 3:} $m = h(i; c)$. In this case,  for $j \in [1, r]$, by \eqref{eq:h-ij-1} one has 
\[
\begin{cases}
m = h(j; c) \hs \mbox{and} \hs \alpha(j, m) = -\beta_{j^*} \in -\Phi_+, & \hs j \in J_1(i),\\
m = h(j; c)+1 \hs \mbox{and} \hs \alpha(j, m) = \beta_{j^*} \in \Phi_+, & \hs j \in J_2(i),\\
m +1= h(j; c)+1 \hs \mbox{and} \hs \alpha(j, m+1) = \beta_{j^*} \in \Phi_+, & \hs j \in J_3(i),\\
m +1= h(j; c) \hs \mbox{and} \hs \alpha(j, m+1) = -\beta_{j^*} \in -\Phi_+, & \hs j \in J_4(i).\end{cases}
\]
Applying the induction assumption to the terms in \eqref{eq:cmi-11}, one has
\[
-[\alpha_{i^*}:\underline{\alpha}] +  c(i, m+1)  = \sum_{j \in J_2(i) \sqcup J_3(i)} (-a_{j, i}) [\beta_{j^*}:\underline{\alpha}] =\sum_{j^* \prec_c i^*} (-a_{j^*, i^*}) [\beta_{j^*}:\underline{\alpha}].
\]
It follows  again from \eqref{eq:beta-alpha-i} that $c(i, m+1) = [\beta_{i^*}:\underline{\alpha}] = [\alpha(i, m+1):\underline{\alpha}]$.

{\bf Case 4:} $m \in [h(i; c)+1, h-1]$. The the proof of this case is similar to that of {\bf Case 1}. 
\end{proof}

\smallskip
\noindent
{\it Proof of \thref{:finite-main}}.
We again use induction on $(i,m) \in [1, r] \times [0, h+1]$ in the total order in 
 \eqref{eq:order-0} to prove that  $x(i, m) = x_{\gamma(i, m)}$ and that the the mutation relation \eqref{eq:xmi-cmi} is as described.

By definition, $x(i, 0) = x_{\omega_i} = x_{\gamma(i, 0)}$ for every $i \in [1, r]$. Let $(i, m) \in [1, r] \times [0, h]$ and suppose that 
$x(j, m') = x_{\gamma(j, m')}$ for all $(j, m') < (i, m+1)$. Consider again four cases.

{\bf Case 1:} $m \in [0, \, h(i; c)-1]$. By \eqref{eq:xmi-cmi}, \leref{:cim-0}, \eqref{eq:h-ij-1}, and the induction assumption, 
\begin{align*}
x_{c^m\omega_i} x(m+1, i) &= {\bf y}_c^{[c(i, m)]_+} +{\bf y}_c^{[-c(i, m)]_+} \prod_{j=i+1}^r x(j,m)^{-a_{j, i}} \prod_{j=1}^{i-1} x(j,m+1)^{-a_{j, i}}\\
&={\bf y}_c^{[c^m\beta_i:\underline{\alpha}]} + \prod_{j=i+1}^r x_{c^m \omega_j}^{-a_{j, i}} \prod_{j=1}^{i-1} x_{c^{m+1} \omega_j}^{-a_{j, i}}.
\end{align*}
Comparing with \eqref{eq:111}, one sees that $x(i, m+1) = x_{c^{m+1}\omega_i} = x_{\gamma(i, m+1)}$ and the mutation relation in \eqref{eq:xmi-cmi}  is \eqref{eq:111};

{\bf Case 2:} $m =h(i; c)$. By \eqref{eq:xmi-cmi}, \leref{:cim-0}, \eqref{eq:h-ij-1}, and the induction assumption, 
\begin{align*}
x_{-\omega_{i^*}} x(i, m+1) & = 1 + y_{i^*; c} \prod_{j=i+1}^r x(j,m)^{-a_{j, i}} \prod_{j=1}^{i-1} x(j, m+1)^{-a_{j, i}}\\
& = 1 + y_{i^*;c}\prod_{j \in J_1(i)\sqcup J_4(i)}x_{\omega_{-j^*}}^{-a_{j, i}} \prod_{j \in J_2(i) \sqcup J_3(i)} x_{\omega_{j^*}}^{-a_{j, i}} \\
& = 1 + y_{i^*;c}\prod_{i^* \prec_c j^*} x_{\omega_{-j^*}}^{-a_{j^*, i^*}} \prod_{j^* \prec_c i^*} x_{\omega_{j^*}}^{-a_{j^*, i^*}}.
\end{align*}
Comparing with \eqref{eq:110}, one sees that $x(i,m+1) = x_{\omega_{i^*}} = x_{\gamma(i,m+1)}$, 
and the mutation relation in \eqref{eq:xmi-cmi}  is \eqref{eq:110} with $i$ replaced by $i^*$;

{\bf Case 3:} $m \in [h(i; c)+1, h]$. The proof is similar to that of {\bf Case 1} and we omit the details;

{\bf Case 4:} $m =h+1$. The proof is similar to that of {\bf Case 2} and we omit the details.

Looking at the above four cases separately, one checks directly that $x(\bfF(i, m)) = x(i, m)$ for all $(i, m) \in [1, r] \times \ZZ$.
This finishes the proof of \thref{:finite-main}.

\bre{:cim}
{\rm Using \leref{:cim-0} and by looking at the four cases in the above  proof of \thref{:finite-main}, one can also check that $c(\bfF(i, m)) = c(i, m)$ for all $(i, m) \in [1, r] \times \ZZ$.
\hfill $\diamond$
}
\ere

\section{Appendix: Proof of \prref{:calM}}\label{appen:calM}
Keeping the notation from $\S$\ref{ss:prin}, we now prove 
\prref{:calM} which says that the set ${\mathcal{M}} = \{\bfx^{\langle \bfm \rangle}: \bfm \in \ZZ^r\}$ of standard monomials is a basis of 
$\calA_{\ZZ}(\bsig^{\rm prin}, \emptyset)$ 
as a module over the polynomial ring $\ZZ[\bfy]$. We prepare some lemmas by modifying some arguments from \cite[$\S$6]{BFZ:III}.

By the Laurent Phenomenon, we have the inclusion 
\[
I: \;\; \calA_\ZZ(\bspr, \emptyset) \hookrightarrow \ZZ[\bfx^{\pm 1}, \bfy].
\]
For $i \in [1, r]$, let
$\calM_{\lan i \ran}$ be the set of all $\bfx^{\langle \bfm \rangle} \in \calM$ with no factor of $x_i$ or $x_i^\prime$. If $S \in \ZZ[\bfy][\calM_{\lan i \ran}]$, then 
$I(S) \in \ZZ[\bfx^{\pm 1}, \bfy]$ is polynomial in $x_i$, so by setting $x_i = 0$ we have the well-defined
\[
I(S)|_{{}_{x_i=0}} \in \ZZ[\bfx_{\lan i \ran}^{\pm 1}, \bfy], 
\]
where 
$\bfx_{\lan i \ran} = (x_1, \ldots, x_{i-1}, x_{i+1}, \ldots, x_r)$. For any $\calM^\prime \subset \calM$ and any sub-ring $R$ of $\ZZ[\bfy]$, let $R[\calM']$ be the $R$-span of $\calM'$ in 
$\calA_\ZZ(\bspr, \emptyset)$. 

\ble{:IS}
For any non-zero $S \in \ZZ[\bfy] (\calM_{\lan i \ran})$, one has $I(S)|_{{}_{x_i=0}} \neq 0$;
\ele

\begin{proof}
The case of $i = 1$ is part of \cite[Lemma 6.2]{BFZ:III}. We now modify the proof of 
\cite[Lemma 6.2]{BFZ:III} for an arbitrary $i \in [1, r]$. Let $\sigma$ be the permutation on $[1, r]$ given by 
\[
(\sigma(1), \, \ldots, \, \sigma(r)) = (i, \, i+1, \, \ldots,\,  r,\,  i-1,\,  \ldots, \, 2,\,  1).
\]
 For
${\bfm}, {\bfm}^\prime \in \ZZ^r$ and 
 ${\bfm}' \neq {\bfm}$, define 
${\bfm} \prec_i {\bfm}^\prime$ if the first non-zero entry of $\sigma({\bfm}) -\sigma({\bfm}^\prime)$ is negative, where (in this proof we write
elements in $\ZZ^r$ as row vectors)
\[
\sigma(m_1,\,   \ldots, \, m_r) = 
(m_{\sigma(1)}, \ldots, m_{\sigma(r)}) = (m_i, \,m_{i+1}, \,\ldots, \,m_r, \,m_{i-1},\, \ldots, \, m_2, \, m_1).
\]
Let $\prec_i$ also denote the corresponding total orders on $\calM$
and on the set $\calM^o$ of ordinary Laurent
monomials of $\bfx$. In particular we have
\begin{equation}\label{eq:prec-x}
 1\prec_i x_{1} \prec_i x_2 \prec_i \cdots \prec_i x_{i-1} \prec_i x_r \prec_i \cdots \prec_i x_{i+1} \prec_i x_i.
\end{equation}
For a non-zero $L \in \ZZ[\bfx^{\pm 1}, \bfy]$,
let ${\rm lm}^o(L) =\bfx^{d(L)}$, where $d(L) \in \ZZ^r$, be 
the {\it lowest} Laurent monomial of $\bfx$ in $L$ (without the coefficient in $\ZZ[{\bf y}]$) with respect to the total order $\prec_i$ on $\calM^o$.
Using \eqref{eq:prec-x} and \eqref{eq:ppi}, one checks directly that
for any $j \in [1, r]$,
\[
{\rm lm}^o(I((x_j^\prime)) = \begin{cases} x_j^{-1} q_j^+, & \hs j \leq i, \\ x_j^{-1} q_j^-, & \hs j > i,\end{cases}
\]
where $q_j^+$ and $q_j^-$ are given in \eqref{eq:ppi}. Thus for any $j \in [1, r]$ there exist $c_{\sigma^{-1}(j), k} \in \ZZ_{\geq 0}$ for $k \in [\sigma^{-1}(j)+1, r]$ such that for any $m_j \in \ZZ$, 
\[
\sigma(d(I(x_j^{\lan m_j \ran}))) = \left(0,\,  \ldots, \, 0, \, m_j,\,  [-m_j]_+c_{\sigma^{-1}(j), \sigma^{-1}(j)+1}, \, \ldots, \, [-m_j]_+c_{\sigma^{-1}(j), r}\right).
\]
where $m_j$ is at the $\sigma^{-1}(j)$th entry. 
For $\bfm \in \ZZ^r$, writing  $\sigma(\bfm) =(\sigma(\bfm)_1, \ldots, \sigma(\bfm)_r)$, then
\[
\left(\!\sigma\!\left(d\left(I(\bfx^{\lan \bfm \ran})\right)\!\right)\!\right)_k  = \sigma(\bfm)_k + \sum_{l=1}^{k-1} [-\sigma(\bfm)_l]_+ c_{l, k}, \hs k \in [1, r].
\]
It follows that for any $\bfm, \bfm' \in \ZZ^r$, 
\begin{equation}\label{eq:bfm-bfm}
{\bfm}\prec_i {\bfm}^\prime \hs \mbox{implies} \hs  d(I(\bfx^{\lan {\bfm}\ran})) \prec_i d(I(\bfx^{\lan {\bfm}^\prime\ran})).
\end{equation}

Let $S \in \ZZ[\bfy](\calM_{\lan i \ran})$ be non-zero, and suppose that $I(S)|_{x_i=0} = 0$. Then every monomial term in $I(S)$, so in particular ${\rm lm}^o(I(S))$, contains a positive power of $x_i$. Let $\bfx^{\langle \bfm \rangle}$
be the {\it lowest standard monomial} in $S$ as a $\ZZ[\bfy]$-linear combination of 
standard monomials in the total order $\prec_i$ on $\calM$. By \eqref{eq:bfm-bfm}, 
${\rm lm}^o(I(S)) = {\rm lm}^o(I(\bfx^{\lan {\bfm} \ran}))$.
Since
$\bfx^{\lan {\bfm} \ran} \in \calM_{\lan i \ran}$, one knows that
$I(\bfx^{\lan {\bfm} \ran})|_{{}_{x_i=0}}$ is a (non-zero) monomial in $\bfx_{\lan i \ran}$, so 
${\rm lm}^o(I(\bfx^{\lan {\bfm} \ran}))$ contains no positive power of $x_i$,  a contradiction.
\end{proof}

For $i \in [1, r]$, let $\bfy_{\lan i \ran} = (y_1, \ldots, y_{i-1}, y_{i+1}, \ldots, y_r)$ and consider the ring homomorphism
\[
\ZZ[\bfx^{\pm 1},\, \bfy] \longrightarrow \ZZ[\bfx^{\pm 1},\, \bfy_{\lan i \ran}], \;\; Z \longmapsto  Z|_{y_i = 0}.
\]

\ble{:xiyi-0}
For every $i \in [1, r]$ and any
non-zero $T \in \ZZ[\bfy_{\lan i \ran}](\calM)$, one has  $I(T)|_{{}_{y_i=0}} \neq 0$.
\ele

\begin{proof} 
 Let $T \in \ZZ[\bfy_{\lan i \ran}][\calM]$ be non-zero and write  
$T = \sum_{k =0}^l T_k (x_i^\prime)^k$,
where $T_0 \in \ZZ[x_i, \bfy_{\lan i \ran}][\calM_{\lan i\ran}]$, 
and $T_k \in  \ZZ[\bfy_{\lan i \ran}][\calM_{\lan i\ran}]$  for $k > 0$.
With again  $q_i^+$ and $q_i^-$ given in \eqref{eq:ppi}, one has
\[
I(T) =  \sum_{k =0}^l I(T_k) \left(\frac{y_i q_i^+ + q_i^-}{x_i}\right)^k.
\]
Since $I(x_j^\prime)|_{{}_{y_i=0}}= I(x_j^\prime)$ for every $j \neq i$, we have  $I(T_k)|_{{}_{y_i=0}} = I(T_k)$ for every $k \geq 0$. 
If $l  =0$, then $T = T_0 \neq 0$, and $I(T)|_{{}_{y_i=0}} = I(T_0)  \neq 0$ and we are done. Assume that $l > 0$, so 
\[
I(T)|_{{}_{y_i=0}}= \sum_{k =0}^{l} I(T_k) \left(\frac{q_i^-}{x_i}\right)^k.
\]
By \leref{:IS}, $I(T_l)|_{x_i = 0} \neq 0$, so $I(T)|_{{}_{y_i=0}}$, when
expressed as a Laurent polynomial in $x_i$ with coefficients in $\ZZ[\bfx_{\lan i \ran}^{\pm 1}, \bfy_{\lan i \ran}]$,
 has $x_i^{-l}$ as the lowest order term. Thus $I(T)|_{{}_{y_i=0}} \neq 0$.
\end{proof}

\medskip
\noindent
{\it Proof of \prref{:calM}}.  
Let $x$ be any cluster variable of $\bspr$. As  ${\mathcal{M}}$ is a basis of $\calA_{\ZZ}(\bsig^{\rm prin})$ as a module over the Laurent polynomial ring $\ZZ[\bfy^{\pm 1}]$,
we can write $x$ uniquely as
\[
x = \sum_{{\bfm} \in {\bf M}} \lambda_{{\bfm}} \bfx^{\lan \bfm \ran},
\]
where ${\bf M} \subset \ZZ^r$ is a finite set, and $\lambda_\bfm \in \ZZ[\bfy^{\pm 1}]$ is non-zero for each $\bfm \in {\bf M}$. 
We want to show that $\lambda_\bfm \in \ZZ[\bfy]$ for every $\bfm \in {\bf M}$. Suppose that this is not the case. Then, by {\it clearing the 
denominators} in the Laurent expansions of the $\lambda_{\bfm}$'s in $\bfy$, we can find ${\bf k} = (k_1, \ldots, k_r) \in (\ZZ_{\geq 0})^r$ with
$k_i > 0$ for some $i \in [1, r]$ and a non-empty subset ${\bf M}'$ of ${\bf M}$ such that
\begin{equation}\label{eq:k-x}
\bfy^{\bf k} x = \sum_{{\bfm} \in {\bf M}'} \mu_{{\bfm}} \bfx^{\lan \bfm \ran} + y_i X,
\end{equation}
where $\mu_{{\bfm}} \in \ZZ[\bfy_{\lan i \ran}]$ for  each $\bfm \in {\bf M}'$ and $X \in \ZZ[\bfy][\calM]$. Applying $I$ to both sides of \eqref{eq:k-x} and setting $y_i = 0$, we get $I(T)|_{{}_{y_i=0}} = 0$ where
$T=  \sum_{{\bfm} \in {\bf M}'} \mu_{{\bfm}} \bfx^{\lan \bfm \ran}$ is a non-zero element in $\ZZ[\bfy_{\lan i \ran}][\calM]$, a contradiction by \leref{:xiyi-0}. This finishes the proof of \prref{:calM}.

\bibliographystyle{amsalpha}
\bibliography{ref}
\end{document}